\documentclass[11pt,oneside]{amsart}

\usepackage{amsmath}

\usepackage{tikz-cd}

\usepackage{fix-cm}

\DeclareMathAlphabet{\mathcal}{OMS}{cmsy}{m}{n}

\usepackage[arrow,curve,matrix,tips,2cell]{xy}

\usepackage[hmargin=2cm,vmargin=1.5cm]{geometry}

\usepackage{anyfontsize}
\usepackage{mathrsfs}

\usepackage{bm}

\usepackage{bbm}

\usepackage{fancyhdr}
\usepackage[obeyspaces]{url}

\usepackage{tensor}

\usepackage{mathtools}

\usepackage{accents}

\usepackage{amscd}
\usepackage{amssymb}
\usepackage{latexsym}
\usepackage{url}

\usepackage{graphicx}

\newtheorem{theorem}{Theorem}[section]
\newtheorem*{theorem*}{Theorem}
\newtheorem{lemma}[theorem]{Lemma}
\newtheorem*{lemma*}{Lemma}
\newtheorem{corollary}[theorem]{Corollary}
\newtheorem{proposition}[theorem]{Proposition}

\newtheorem{remark}[theorem]{Remark}
\newtheorem{definition}[theorem]{Definition}
\newtheorem*{definition*}{Definition}

\newtheorem{question}[theorem]{Question}
\newtheorem*{question*}{Question}

\newtheorem{example}[theorem]{Example}
\newtheorem{examples}[theorem]{Examples}

\makeatletter
\def\revddots{\mathinner{\mkern1mu\raise\p@
\vbox{\kern7\p@\hbox{.}}\mkern2mu
\raise4\p@\hbox{.}\mkern2mu\raise7\p@\hbox{.}\mkern1mu}}
\makeatother 
\newcommand{\bgl}{\begin{equation}} 
\newcommand{\egl}{\end{equation}}
\newcommand{\bgloz}{\begin{equation*}} 
\newcommand{\egloz}{\end{equation*}}
\newcommand{\bgln}{\begin{eqnarray}} 
\newcommand{\egln}{\end{eqnarray}}
\newcommand{\bglnoz}{\begin{eqnarray*}} 
\newcommand{\eglnoz}{\end{eqnarray*}}
\newcommand{\btheo}{\begin{theorem}}
\newcommand{\etheo}{\end{theorem}}
\newcommand{\btheooz}{\begin{theorem*}}
\newcommand{\etheooz}{\end{theorem*}}

\newcommand{\blemma}{\begin{lemma}}
\newcommand{\elemma}{\end{lemma}}
\newcommand{\blemmaoz}{\begin{lemma*}}
\newcommand{\elemmaoz}{\end{lemma*}}
\newcommand{\bproof}{\begin{proof}}
\newcommand{\eproof}{\end{proof}}
\newcommand{\bbew}{\begin{beweis}}
\newcommand{\ebew}{\end{beweis}}
\newcommand{\bremark}{\begin{remark}\em}
\newcommand{\eremark}{\end{remark}}
\newcommand{\bdefin}{\begin{definition}}
\newcommand{\edefin}{\end{definition}}
\newcommand{\bdefinoz}{\begin{definition*}}
\newcommand{\edefinoz}{\end{definition*}}
\newcommand{\bex}{\begin{example}\em}
\newcommand{\eex}{\end{example}}
\newcommand{\bexs}{\begin{examples}}
\newcommand{\eexs}{\end{examples}}

\newcommand{\bprop}{\begin{proposition}}
\newcommand{\eprop}{\end{proposition}}
\newcommand{\bcor}{\begin{corollary}}
\newcommand{\ecor}{\end{corollary}}
\newcommand{\bfa}{\begin{cases}} 
\newcommand{\efa}{\end{cases}}

\newcommand{\bquestion}{\begin{question}}
\newcommand{\equestion}{\end{question}}
\newcommand{\bquestionoz}{\begin{question*}}
\newcommand{\equestionoz}{\end{question*}}
\newtheorem{introtheorem}{Theorem}

\newtheorem{introcor}[introtheorem]{Corollary}

\newcommand{\cC}{\mathcal C}

\newcommand{\cG}{\mathcal G}

\newcommand{\cM}{\mathcal M}
\newcommand{\cN}{\mathcal N}
\newcommand{\cO}{\mathcal O}
\newcommand{\cP}{\mathcal P}

\newcommand{\cR}{\mathcal R}

\newcommand{\cU}{\mathcal U}
\newcommand{\cV}{\mathcal V}

\def\Kz{\mathbb{K}}

\def\Nz{\mathbb{N}}

\def\Qz{\mathbb{Q}}
\def\Rz{\mathbb{R}}

\def\Zz{\mathbb{Z}}
\def\Sz{\mathbb{S}}
\newcommand{\fB}{\mathfrak B}
\newcommand{\fC}{\mathfrak C}

\newcommand{\fF}{\mathfrak F}
\newcommand{\fG}{\mathfrak G}
\newcommand{\fH}{\mathfrak H}

\newcommand{\fL}{\mathfrak L}

\newcommand{\fN}{\mathfrak N}
\newcommand{\fP}{\mathfrak P}

\newcommand{\mfd}{\mathfrak d}

\newcommand{\mfl}{\mathfrak l}

\newcommand{\mft}{\mathfrak t}

\newcommand{\an}[1]{``#1''} 
\newcommand{\ti}{\tilde}

\newcommand{\ma}{\mapsto} 
\newcommand{\into}{\hookrightarrow} 
\newcommand{\Rarr}{\Rightarrow} 

\def\SEMI{\mbox{$\times\kern-2pt\vrule height5pt width.6pt \kern3pt $}}

\newcommand{\id}{{\rm id}}

\renewcommand{\dim}{{\rm dim}\,}

\renewcommand{\ker}{{\rm ker}\,}
\newcommand{\coker}{{\rm coker}\,}
\newcommand{\reg}{^\times} 
\newcommand{\lspan}{{\rm span}} 
\newcommand{\abs}[1]{\lvert#1\rvert} 
\newcommand{\defeq}{\mathrel{:=}} 

\newcommand{\dop}{\text{: }} 
\newcommand{\falls}{\text{ if }} 
\newcommand{\sonst}{\text{ else}} 
\newcommand{\ilim}{\varinjlim} 
\newcommand{\plim}{\varprojlim} 
\newcommand{\lge}{\left\{} 
\newcommand{\rge}{\right\}} 
\newcommand{\lru}{\left(} 
\newcommand{\rru}{\right)} 
\newcommand{\lsp}{\left\langle} 
\newcommand{\rsp}{\right\rangle} 
\newcommand{\rukl}[1]{\lru #1 \rru} 
\newcommand{\gekl}[1]{\lge #1 \rge} 
\newcommand{\spkl}[1]{\lsp #1 \rsp} 
\newcommand{\menge}[2]{\gekl{ #1 \dop #2 }} 
\def\bf1{\mathbf{1}}
\newcommand{\bma}{\bm{a}}

\newcommand{\bmf}{\bm{f}}

\newcommand{\bmA}{\bm{A}}
\newcommand{\bmD}{\bm{D}}

\newcommand{\bmF}{\bm{F}}
\newcommand{\bmK}{\bm{K}}
\newcommand{\bmS}{\bm{S}}

\newcommand{\oset}[2]{%
  \mathop{#2}\limits^{\vbox to -1.66ex{%
  \kern -1.4ex\hbox{$#1$}\vss}}}

\newcommand{\C}{\mathscr{C}}

\newcommand{\I}{\mathscr{I}}

\newcommand{\crO}{\mathscr{O}}

\newcommand{\Gn}{{G^{(0)}}}
\newcommand{\Gnu}{{G^{(\nu)}}}
\newcommand{\Gnum}{{G^{(\nu - 1)}}}
\newcommand{\Gnup}{{G^{(\nu + 1)}}}
\newcommand{\fBnu}{{\fB^{(\nu)}}}
\newcommand{\fBnum}{{\fB^{(\nu - 1)}}}
\newcommand{\fBnup}{{\fB^{(\nu + 1)}}}
\newcommand{\rmr}{\ensuremath{\mathrm{r}}}
\newcommand{\rms}{\ensuremath{\mathrm{s}}}

\newcommand{\bbs}{\mathbbm{s}}
\newcommand{\bbt}{\mathbbm{t}}

\newcommand{\obj}{{\rm obj}\,}
\newcommand{\mor}{{\rm mor}\,}

\newcommand{\sfA}{\mathsf{A}}
\newcommand{\sfC}{\mathsf{C}}

\newcommand{\Tor}{{\rm Tor\,}}

\newcommand{\fBZ}{\fB_{G \curvearrowright Z}}
\newcommand{\acts}{\curvearrowright}

\DeclareFontFamily{U}{mathb}{\hyphenchar\font45}
\DeclareFontShape{U}{mathb}{m}{n}{
      <5> <6> <7> <8> <9> <10>
      <10.95> <12> <14.4> <17.28> <20.74> <24.88>
      mathb10
      }{}
\DeclareSymbolFont{mathb}{U}{mathb}{m}{n}

\DeclareMathSymbol{\sqbullet}{1}{mathb}{"0D}
\usepackage{newtxtext}
\usepackage{newtxmath}

\begin{document}

\title{Ample groupoids, topological full groups,\\algebraic K-theory spectra and infinite loop spaces}

\thispagestyle{fancy}

\author{Xin Li}

\address{Xin Li, School of Mathematics and Statistics, University of Glasgow, University Place, Glasgow G12 8QQ, United Kingdom}
\email{Xin.Li@glasgow.ac.uk}


\thanks{This project has received funding from the European Research Council (ERC) under the European Union's Horizon 2020 research and innovation programme (grant agreement No. 817597).}

\begin{abstract}
Inspired by work of Szymik and Wahl on the homology of Higman-Thompson groups, we establish a general connection between ample groupoids, topological full groups, algebraic K-theory spectra and infinite loop spaces, based on the construction of small permutative categories of compact open bisections. This allows us to analyse homological invariants of topological full groups in terms of homology for ample groupoids. 

Applications include complete rational computations, general vanishing and acyclicity results for group homology of topological full groups as well as a proof of Matui's AH-conjecture for all minimal, ample groupoids with comparison.
\end{abstract}

\maketitle


\setlength{\parindent}{0cm} \setlength{\parskip}{0.5cm}

\section{Introduction}

The construction of topological full groups has recently attracted attention because it led to solutions of several outstanding open problems in group theory. For instance, it gave rise to first examples of finitely generated, infinite simple amenable groups (see \cite{JM} and also \cite{JNS}). Topological full groups also led to the first examples of finitely generated, infinite simple groups with intermediate growth \cite{Nek18b}. Moreover, this construction also produces new families of infinite simple groups with prescribed finiteness properties \cite{SWZ}. 

Topological full groups arise from generalized dynamical systems in the form of topological groupoids, which describe orbit structures of dynamical systems in situations where the actual space of orbits might be very badly behaved. Topological groupoids only capture local symmetries arising in dynamical systems, which is enough to determine their orbit structures. Roughly speaking, elements of topological full groups are global symmetries which are pieced together from local symmetries encoded by topological groupoids.

Topological groupoids and their topological full groups arise in a variety of settings, for instance from topological dynamical systems given by actions of groups on topological spaces by homeomorphisms. Indeed, the first examples of topological full groups were studied in \cite{Kri,GPS99} in the setting of Cantor minimal systems and the closely related context of Bratteli diagrams. They also arise from shifts of finite type, or more generally, from graphs (see for instance \cite{Mat15}). Further examples have been constructed from self-similar groups or actions and higher rank graphs (see for instance \cite{Mat16,Nek18a}). In this context, there is an interesting connection to C*-algebra theory because topological groupoids serve as models for C*-algebras (see \cite{Ren}) such as Cuntz algebras, Cuntz-Krieger algebras, graph C*-algebras or higher rank graph C*-algebras, many of which play distinguished roles in the classification programme for C*-algebras. There is also an interesting link to group theory because Thompson's group $V$ and many of its generalizations and variations \cite{Hig,Ste92,Bri04} can be described as topological full groups of corresponding topological groupoids. In the case of $V$ this observation goes back to \cite{Nek04}. This gives a dynamical perspective on Thompson-like groups, which have been popular and important objects of study in group theory ever since the introduction of $V$ by Thompson (see for instance \cite{CFP}).

While general structural properties \cite{Mat12,Mat15,Mat16,Nek19,MB} and rigidity results have been developed \cite{Rub,Med11,Mat15}, and several deep results have been established for particular examples of topological full groups \cite{JM,JNS,Nek18b,SWZ,SW}, it would be desirable to create a dictionary between dynamical properties and invariants of topological groupoids on the one hand and group-theoretic properties and invariants of topological full groups on the other hand. This would allow us to study topological full groups -- which are very interesting but in many aspects still remain mysterious -- through the underlying topological groupoids which are often much more accessible. The goal of this paper is to develop this programme in the context of homological invariants by establishing a link between groupoid homology and group homology of topological full groups. This leads to a better understanding of this class of groups arising from dynamics, not only because it allows us to compute group homology, which is a fundamental invariant, but also because the methods developed to accomplish this task reveal interesting connections.

For the particular example class of Thompson's group $V$ and its generalizations, the study of homological invariants and properties has a long history \cite{BG,Bro92}. It was shown in \cite{Bro92} that $V$ is rationally acyclic. Only recently it was established in \cite{SW} that $V$ is even integrally acyclic. The new approach in \cite{SW} also allows for many more homology computations for Higman-Thompson groups. However, for other classes of topological full groups, very little is known about homological invariants. In degree one, Matui has formulated the AH-conjecture, which describes $H_1$ (i.e., the abelianization) of topological full groups in terms of groupoid homology of the underlying topological groupoids. This AH-conjecture has been verified for general classes of topological full groups (for instance for almost finite, principal groupoids, see \cite{Mat12,Mat15,Mat16}) as well as for several example classes (for instance for groupoids arising from shifts of finite type \cite{Mat15,Mat16}, graphs \cite{NO21a} or self-similar actions \cite{NO21b}, as well as transformation groupoids of odometers \cite{Sca19} and Cantor minimal dihedral systems \cite{Sca22}). However, no general results of this nature are known concerning homology groups in higher degree.

In this paper, we develop a new approach to homological invariants of topological full groups. The key novelties are the construction of small permutative categories of bisections for all ample groupoids and the realization of groupoid homology as (reduced) stable homology of the associated algebraic K-theory spectra. Another key ingredient is the identification of homology of the corresponding infinite loop spaces with group homology of the topological full groups we are interested in. This last identification is inspired by \cite{SW} and at the same time vastly generalizes corresponding results on the particular example class of Higman-Thompson groups in \cite{SW}. Our new insights allow us to apply powerful tools from algebraic topology to the study of homology of topological full groups, bringing together group theory, topological dynamics, algebraic topology as well as ideas from operator algebras. Among other things, our insights lead to
\setlength{\parindent}{0cm} \setlength{\parskip}{0cm}

\begin{itemize}
\item a complete description of rational group homology for large classes of topological full groups,
\item general vanishing and acyclicity results, explaining and generalizing the result that $V$ is acyclic in \cite{SW},
\item a verification of Matui's AH-conjecture for a general class of ample groupoids, including all purely infinite and minimal ones.
\end{itemize}
We establish these results under very mild assumptions, i.e., for all ample groupoids which are minimal, whose unit spaces do not have isolated points, and which have comparison. The first two conditions are necessary, as we explain below in our discussion of Theorem~\ref{thm:IntroB}. Comparison appears naturally and is in itself an interesting property. It has been verified in many situations \cite{DZ,KN,GGKN}, and there is the conjecture that comparison holds in great generality, as we explain below. For topological full groups of amplified groupoids, we prove analogous results in complete generality, i.e., for all ample groupoids. The present work is a significant step forward in our understanding of homological invariants of ample groupoids and topological full groups, both at the conceptual level as well as concerning concrete applications. Indeed, our results on rational group homology are the first explicit computations of that kind which work in all degrees. The acyclicity results imply, for example, that all of Brin's groups $nV$ are integrally acyclic, and that all Brin-Higman-Thompson groups $nV_{k,r}$ are rationally acyclic. In addition, we are able to construct continuum many pairwise non-isomorphic infinite simple groups which are all integrally acyclic. Moreover, our work leads to a conceptual explanation and strengthening of Matui's AH-conjecture as we obtain precise obstructions for the strong AH-conjecture and establish that the amplified version of the AH-conjecture is always true (i.e., for all ample groupoids).
\setlength{\parindent}{0cm} \setlength{\parskip}{0.5cm}

Let us now formulate our main results. Let $G$ be a topological groupoid, i.e., a topological space which is at the same time a small category with invertible morphisms, such that all operations (range, source, multiplication and inversion maps) are continuous. We always assume the unit space $\Gn$ consisting of the objects of $G$ to be locally compact and Hausdorff. In addition, suppose that $G$ is ample, in the sense that it has a basis for its topology given by compact open bisections (see \S~\ref{ss:GPD}). If $\Gn$ is compact, then the topological full group $\bmF(G)$ is defined as the group of global compact open bisections. In the general case, $\bmF(G)$ is the inductive limit of topological full groups of restrictions of $G$ to compact open subspaces of $\Gn$. The new examples of infinite simple groups mentioned at the beginning (see \cite{JM,JNS,Nek18b,SWZ}) are given by commutator subgroups $\bmD(G)$ of $\bmF(G)$. Given an ample groupoid $G$ as above, we construct a small permutative category $\fB_G$ of compact open bisections of $G$ (see \S~\ref{s:smallpermcat}). Let $\Kz(\fB_G)$ be the algebraic K-theory spectrum of $\fB_G$ and $\Omega^{\infty} \Kz(\fB_G)$ the associated infinite loop space (see \S~\ref{ss:AlgKSmallPermCat}). 

Our first main result identifies the (reduced) stable homology of $\Kz(\fB_G)$ with the groupoid homology of $G$ as introduced in \cite{CM} and studied in \cite{Mat12}.
\begin{introtheorem}[see Theorem~\ref{thm:HKB=HG}]
\label{thm:IntroA}
Let $G$ be an ample groupoid with locally compact Hausdorff unit space. Then we have
$$
 \ti{H}_*(\Kz(\fB_G)) \cong H_*(G).
$$
\end{introtheorem}
For the second main result, we need the assumption that $G$ is minimal, i.e., every $G$-orbit is dense in $\Gn$, and that the unit space of $G$ does not have isolated points. These two conditions are necessary for Theorem~\ref{thm:IntroB}. We also require $G$ to have comparison, which roughly means that $G$-invariant measures on $\Gn$ control when one compact open subspace of $\Gn$ can be transported into another by compact open bisections of $G$. Comparison appears naturally and is needed for the key ingredient, Morita invariance (see Theorem~\ref{introthm:Morita}), which allows us to compare compact open subspaces of the unit space and the corresponding topological full groups of the restricted groupoids in homology. Under these assumptions, we can identify group homology of the topological full group $\bmF(G)$ with the homology of $\Omega^{\infty}_0 \Kz(\fB_G)$, the connected component of the base point in $\Omega^{\infty} \Kz(\fB_G)$.
\begin{introtheorem}[see Theorem~\ref{thm:HFG=HOKB}]
\label{thm:IntroB}
Let $G$ be an ample groupoid whose unit space is locally compact Hausdorff without isolated points. Assume that $G$ is minimal and has comparison. Then we have
$$
 H_*(\bmF(G)) \cong H_*(\Omega^{\infty}_0 \Kz(\fB_G)).
$$
\end{introtheorem}
Note that the group completion theorem (see \cite{McDS,R-W}) gives us a similar isomorphism in homology, where the left-hand side is replaced by homology of the amplified version of the topological full group. However, that alone is not enough to derive Theorem~\ref{thm:IntroB}. Indeed, the key (and also most demanding) step is to show that the amplified version of the topological full group and the topological full group itself have the same group homology. This is achieved by Morita invariance (Theorem~\ref{introthm:Morita}), which plays the role of homological stability in \cite{SW} (see also \cite{RWW}).

Groupoid homology is much more accessible than group homology of topological full groups because there are many tools to compute groupoid homology, and several computations have been produced for various example classes of ample groupoids (see \S~\ref{sss:GPDHomEx}). Thus the point of our two main results is that they enable us to study group homology of topological full groups in terms of groupoid homology, provided we understand how to relate $H_*(\Omega^{\infty}_0 \Kz(\fB_G))$ to $\ti{H}_*(\Kz(\fB_G))$. This problem has been studied in algebraic topology, where powerful tools have been developed. The precise relation between homology of infinite loop spaces and the corresponding spectra is not easy to understand. But we can still derive several consequences. In the following, let us present a selection of such consequences.

In order to present our results on rational group homology, we need the following notation:
$$
 H^{\rm odd}_*(G,\Qz) \defeq
 \bfa
 H_*(G,\Qz) & \falls * > 0 \text{ odd},\\
 \gekl{0} & \sonst,
 \efa
 \qquad \text{and} \quad
 H^{\rm odd}_{*>1}(G,\Qz) \defeq
 \bfa
 H_*(G,\Qz) & \falls * > 1 \text{ odd},\\
 \gekl{0} & \sonst,
 \efa
$$
as well as 
$$
 H^{\rm even}_*(G,\Qz) \defeq
 \bfa
 H_*(G,\Qz) & \falls * > 0 \text{ even},\\
 \gekl{0} & \sonst.
 \efa
$$
\begin{introcor}[see Corollaries~\ref{cor:RatHom} and \ref{cor:DRatHom}]
Let $G$ be an ample groupoid, with locally compact Hausdorff unit space without isolated points. Assume that $G$ is minimal and has comparison. Then, as graded vector spaces over $\Qz$,
$$
 H_*(\bmF(G),\Qz) \cong {\rm Ext}(H^{\rm odd}_*(G,\Qz)) \otimes {\rm Sym}(H^{\rm even}_*(G,\Qz)).
$$

In particular, $\bmF(G)$ is rationally acyclic (i.e., $H_*(\bmF(G),\Qz) \cong \gekl{0}$ for all $* > 0$) if and only if $H_*(G,\Qz) \cong \gekl{0}$ for all $*>0$. 

For the commutator subgroup $\bmD(G)$ of $\bmF(G)$, we obtain, again as graded vector spaces over $\Qz$,
$$
 H_*(\bmD(G),\Qz) \cong {\rm Ext}(H^{\rm odd}_{*>1}(G,\Qz)) \otimes {\rm Sym}(H^{\rm even}_*(G,\Qz)).
$$
\end{introcor}
\setlength{\parindent}{0cm} \setlength{\parskip}{0cm}

Here ${\rm Ext}$ stands for exterior algebra (see for instance \cite[\S~5]{Gre}) and ${\rm Sym}$ stands for symmetric algebra (see for instance \cite[\S~9]{Gre}).
\setlength{\parindent}{0cm} \setlength{\parskip}{0.5cm}

Next, we present vanishing results which generalize and provide a conceptual explanation for the result in \cite{SW} that $V$ is acyclic.

\begin{introcor}[see Corollaries~\ref{cor:HVanish}, \ref{cor:Acyclic} and \ref{cor:HDVanish}]
\label{introcor:Vanish}
Let $G$ be an ample groupoid whose unit space is locally compact Hausdorff and does not have isolated points. Assume that $G$ is minimal and has comparison. 

Suppose that $k \in \Zz$ with $k > 0$. If $H_*(G) \cong \gekl{0}$ for all $* < k$, then $H_*(\bmF(G)) \cong \gekl{0}$ for all $0 < * < k$ and $H_k(\bmF(G)) \cong H_k(G)$. If $k \geq 2$, then this implies $\bmF(G) = \bmD(G)$. In particular, if $H_*(G) \cong \gekl{0}$ for all $* \geq 0$, then $\bmF(G)$ is integrally acyclic, i.e., $H_*(\bmF(G)) \cong \gekl{0}$ for all $* > 0$, and $\bmF(G) = \bmD(G)$. 

For the commutator subgroup, we always have $H_1(\bmD(G)) \cong \gekl{0}$.
\end{introcor}
\setlength{\parindent}{0cm} \setlength{\parskip}{0cm}

Concrete examples where Corollary~\ref{introcor:Vanish} applies can be found in \S~\ref{ss:Examples}. In particular, we construct continuum many pairwise non-isomorphic infinite simple groups which are all integrally acyclic (see Remark~\ref{rem:ManyAcyclic}).
\setlength{\parindent}{0cm} \setlength{\parskip}{0.5cm}

In low degrees, we obtain the following exact sequence.
\begin{introcor}[see Corollary~\ref{cor:AHConj}]
Let $G$ be an ample groupoid, with locally compact Hausdorff unit space without isolated points. Assume that $G$ is minimal and has comparison. Then there is an exact sequence
$$
 \xymatrix{
 H_2(\bmD(G)) \ar[r] & H_2(G) \ar[r] & H_0(G,\Zz/2) \ar[r]^{\zeta} & H_1(\bmF(G)) \ar[r]^{\hspace{0.25cm} \eta} & H_1(G) \ar[r] & 0.
 }
$$
The maps $\eta$ and $\zeta$ coincide with the ones in \cite[\S~2.3]{Mat16} and \cite[\S~7]{Nek19}. 

In particular, Matui's AH-conjecture is true for ample groupoids $G$ which are minimal, have comparison and whose unit spaces are locally compact Hausdorff without isolated points.
\end{introcor}
\setlength{\parindent}{0cm} \setlength{\parskip}{0cm}

Note that this in particular verifies Matui's AH-conjecture for all purely infinite minimal ample groupoids, which was not known before. Our result also verifies the AH-conjecture for all minimal ample groupoids which are $\sigma$-compact, Hausdorff and almost finite, and whose unit spaces are compact without isolated points. Previously, this was only known under the additional assumption of principality \cite{Mat12}.
\setlength{\parindent}{0cm} \setlength{\parskip}{0.5cm}

Our results also lead to several new concrete homology computations. For instance, if $G$ is the transformation groupoid of a Cantor minimal $\Zz$-system, then the commutator subgroup $\bmD(G)$ of $\bmF(G)$ is always rationally acyclic. In particular, this covers the class of infinite simple amenable groups found in \cite{JM}. Moreover, we obtain that all Brin-Higman-Thompson groups $nV_{k,r}$ are rationally acyclic, and that $nV_{2,r}$ is integrally acyclic for all $n$ and $r$ (this is the case $k=2$). We also obtain concrete computations of rational group homology for topological full groups of certain tiling groupoids, graph groupoids and groupoids attached to self-similar actions such as Katsura-Exel-Pardo groupoids. Moreover, we compute rational group homology for classes of Thompson-like groups introduced by Stein in \cite{Ste92}, including irrational slope versions of Thompson's group $V$. The reader may consult \S~\ref{ss:Examples} for results concerning homology for concrete example classes of topological full groups.

Let us now explain the main ideas. The construction of $\fB_G$ is a key ingredient. The underlying category of $\fB_G$ actually already appeared in \cite[\S~2]{Li21b} in the context of finiteness properties, but the extra structure making $\fB_G$ a small permutative category has not been exploited before. For us, the key insight is that for general ample groupoids, we can replace the small permutative categories constructed from Cantor algebras in \cite{SW} in the setting of Higman-Thompson groups by the categories of bisections $\fB_G$. This allows us to treat general ample groupoids, for which the notion of Cantor algebras is not available. Actually, from the point of view of groupoids, $\fB_G$ is more natural because it takes into account all compact open subspaces of the unit space of our groupoid, whereas for groupoids giving rise to Higman-Thompson groups, the small permutative categories constructed from Cantor algebras in \cite{SW} are only subcategories of our $\fB_G$. To create the structure of a small permutative category starting with compact open bisections of our ample groupoid $G$, the idea of amplification is crucial, i.e., we pass from $G$ to $\cR \times G$. Here $\cR$ is the full equivalence relation on $\Nz = \gekl{1, 2, 3, \dotsc}$, i.e., $\cR = \Nz \times \Nz$ with the discrete topology. On the C*-algebraic level, this corresponds to passing to matrix algebras, an idea which is at the heart of K-theory. Actually, if we are willing to replace $\bmF(G)$ by $\bmF(\cR \times G)$, then our results above are unconditionally true, i.e., they do not need the assumptions that $\Gn$ has no isolated points, and that $G$ is minimal and has comparison (see Thereom~\ref{thm:GroupCompl} and the results in \S~\ref{s:App}). In particular, we obtain a proof, for general ample groupoids, of a modified AH-conjecture with $\bmF(\cR \times G)$ in place of $\bmF(G)$ (see Theorem~\ref{thm:H2H0H1H1} and Remark~\ref{rem:AmplifiedAHConj}).

To go back from $\bmF(\cR \times G)$ to $\bmF(G)$, at least in homology, we need to establish Morita invariance, which plays the role of homological stability in \cite{SW} (see also \cite{RWW}).
\begin{introtheorem}[see Theorem~\ref{thm:GUUGVV} and Remark~\ref{rem:DMorita}]
\label{introthm:Morita}
Suppose that $G$ is an ample groupoid which is minimal, has comparison, and whose unit space $\Gn$ is locally compact Hausdorff without isolated points. Then for all non-empty compact open subspaces $U \subseteq V$ of $\Gn$, the canonical maps $\bmF(G_U^U) \to \bmF(G_V^V)$ and $\bmD(G_U^U) \to \bmD(G_V^V)$ induce isomorphisms in homology in all degrees.
\end{introtheorem}
\setlength{\parindent}{0cm} \setlength{\parskip}{0cm}

Here $G_U^U$ and $G_V^V$ are the restrictions of $G$ to $U$ and $V$, respectively. Theorem~\ref{introthm:Morita} implies that the homology of Brin-Higman-Thompson groups $n V_{k,r}$ does not depend on $r$ (see \S~\ref{ss:Examples}), just as in the case of Higman-Thompson groups \cite{SW}. Moreover, Theorem~\ref{introthm:Morita} implies that homology of topological full groups and their commutator subgroups is invariant under (Morita) equivalence of groupoids (see Corollary~\ref{cor:Morita} and Remark~\ref{rem:DMorita}).
\setlength{\parindent}{0cm} \setlength{\parskip}{0.5cm}

Interestingly, the notion of comparison also appears in the classification programme of C*-algebras \cite{Ker,KS,KN,MW}. For instance, transformation groupoids of free minimal actions of groups with subexponential growth and elementary amenable groups on the Cantor space have comparison \cite{DZ,KN}. Moreover, all purely infinite minimal groupoids have comparison. This includes transformation groupoids arising from amenable, minimal actions of many non-amenable groups on the Cantor space \cite{GGKN}. At the moment, there is no example of a minimal ample groupoid known which does not have comparison, and there is the conjecture that all transformation groupoids of free minimal group actions on the Cantor space have comparison. Our proof of Morita invariance splits naturally into the case of purely infinite minimal groupoids, where no non-zero invariant measures exist, and the case where non-zero invariant measures do exist (almost finite minimal groupoids, for instance). The cases covered in \cite{SW} belong to the purely infinite setting. In the setting where non-zero invariant measures do exist, no Morita invariance results were known and we had to develop new ideas. In both cases, we analyse connectivity of certain simplicial complexes constructed out of bisections, following the general criterion for homological stability formulated in \cite{RWW}.

Our result identifying (reduced) stable homology of the algebraic K-theory spectrum $\Kz(\fB_G)$ with groupoid homology of $G$ is a completely new insight which does not appear in \cite{SW}. This result is interesting on its own right because it gives a new perspective on groupoid homology, which is a fundamental invariant in topological dynamics. For instance, this invariant plays a key role in the classification of Cantor minimal systems up to topological orbit equivalence \cite{GPS95, GMPS08, GMPS10}. This new insight leads to a conceptual explanation why Thompson's group $V$ is acyclic (as proven in \cite{SW}). This is because, as observed in \cite{Nek04}, $V$ can be identified with the topological full group of the ample groupoid $G_2$, which is the Deaconu-Renault groupoid for the one-sided full shift on two symbols (see \S~\ref{sss:GraphGPD}). And it is known that the homology of $G_2$ vanishes (see \S~\ref{sss:GPDHomEx}). From the point of view of C*-algebras, this can be explained using Matui's HK-conjecture \cite[Conjecture~2.6]{Mat16}, because $G_2$ is a groupoid model for the Cuntz algebra $\cO_2$, whose K-theory vanishes.

Our main results also lead to a better understanding of Matui's AH-conjecture by relating it to the Atiyah-Hirzebruch spectral sequence. Our work on the AH-conjecture demonstrates that our proof of Theorem~\ref{thm:IntroA} reveals more information about the isomorphism $\ti{H}_*(\Kz(\fB_G)) \cong H_*(G)$, allowing us to identify the maps in the Atiyah-Hirzebruch spectral sequence with the ones appearing in Matui's AH-conjecture.

I would like to thank E. Scarparo, O. Tanner and M. Yamashita for very helpful comments and discussions.

\section{Preliminaries}

\subsection{Groupoids}
\label{ss:GPD}

A groupoid is a small category whose morphisms are all invertible. As usual, we identify the groupoid with its set of morphisms, say $G$, and view its set of objects (also called units) $\Gn$ as a subset of $G$ by identifying objects with the corresponding identity morphisms. By definition, our groupoid $G$ comes with range and source maps $\rmr: \: G \to \Gn$, $\rms: \: G \to \Gn$, a multiplication map
$$
 G \tensor[_{\rms}]{\times}{_{\rmr}} G = \menge{(g_1,g_2)}{\rms(g_1) = \rmr(g_2)} \to G, \, (g_1, g_2) \ma g_1 g_2
$$
and an inversion map $G \to G, \, g \ma g^{-1}$ such that $\rmr(g^{-1}) = \rms(g)$, $\rms(g^{-1}) = \rmr(g)$, $g g^{-1} = \rmr(g)$ and $g^{-1} g = \rms(g)$. These structure maps satisfy a list of conditions so that $G$ becomes a small category (see for instance \cite[Chapter~I, Section~1]{Ren}).

We are interested in the case of topological groupoids, i.e., our groupoid $G$ is endowed with a topology such that range, source, multiplication and inversion maps are all continuous. We do not assume that $G$ is Hausdorff, but $\Gn$ is always assumed to be Hausdorff in the subspace topology. We call $\Gn$ the unit space. We will also always assume that $\Gn$ is locally compact. A topological groupoid is called {\'e}tale if the range map (and hence also the source map) is a local homeomorphism. It follows that $\Gn$ is an open subspace of $G$ in that case. An open subspace $U \subseteq G$ is called an open bisection if the restricted range and source maps $\rmr \vert_U: \: U \to \rmr(U), \, g \ma \rmr(g)$, $\rms \vert_U: \: U \to \rms(U), \, g \ma \rms(g)$ are bijections (and hence homeomorphisms). If $G$ is {\'e}tale, then $G$ has a basis for its topology consisting of open bisections. Note that open bisections are always locally compact and Hausdorff because they are homeomorphic to open subspaces of the unit space. A topological groupoid $G$ is called ample if it is {\'e}tale and its unit space $\Gn$ is totally disconnected. An {\'e}tale groupoid is ample if and only if it has a basis for its topology consisting of compact open bisections.

A topological groupoid $G$ is called minimal if for all $x \in \Gn$, the orbit $G.x \defeq \menge{\rmr(g)}{g \in \rms^{-1}(x)}$ is dense in $\Gn$. Let $M(G)$ be the set of all non-zero Radon measures $\mu$ on $\Gn$ which are invariant, i.e., $\mu(\rmr(U)) = \mu(\rms(U))$ for all open bisections $U \subseteq G$. An ample groupoid $G$ is said to have groupoid strict comparison for compact open sets (abbreviated by comparison in the following) if for all non-empty compact open sets $U, V \subseteq \Gn$ with $\mu(U) < \mu(V)$ for all $\mu \in M(G)$, there exists a compact open bisection $\sigma \subseteq G$ with $\rms(\sigma) = U$ and $\rmr(\sigma) \subseteq V$ (see for instance \cite[\S~6]{MW}). Note that we restrict ourselves to non-empty open sets $U$ and $V$ because we want our definition of comparison to cover purely infinite groupoids, where $M(G) = \emptyset$ (see below).

Examples of groupoids with comparison include locally compact $\sigma$-compact Hausdorff ample groupoids with compact unit spaces which are almost finite in the sense of \cite[\S~6]{Mat12} (see also \cite[Proposition~7.2]{MW}). This covers many examples, for instance AF groupoids, classes of transformation groupoids and tiling groupoids, as we explain below in \S~\ref{ss:ExGPDs}. Another class of groupoids with comparison is given by purely infinite minimal groupoids, in the following sense: An ample groupoid $G$ is purely infinite minimal if and only if for all compact open subspaces $U, V \subseteq \Gn$ with $V \neq \emptyset$, there exists a compact open bisection $\sigma \subseteq G$ such that $\rms(\sigma) = U$ and $\rmr(\sigma) \subseteq V$ (compare \cite[\S~4.2]{Mat15}, but we do not require essential principality or Hausdorffness). Concrete examples are discussed in \S~\ref{ss:ExGPDs}. By definition, it is clear that purely infinite minimal groupoids have comparison.

\subsection{Examples of groupoids}
\label{ss:ExGPDs}

Let us discuss several classes of examples of topological groupoids. 

\subsubsection{AF groupoids}
\label{sss:AF}

AF groupoids are inductive limits of elementary groupoids. Here an elementary groupoid is a disjoint union of groupoids of the form $\cR \times X$, where $\cR$ is the full equivalence relation on a finite set and $X$ is a totally disconnected, locally compact Hausdorff space. AF groupoids are represented by Bratteli diagrams. We refer to \cite[Chapter~III, \S~1]{Ren}, \cite[\S~3]{GPS04} and \cite[\S~2]{Mat12} for details.

\subsubsection{Transformation groupoids}
\label{sss:Trafo}

Let $\Gamma$ be a discrete group acting on a locally compact Hausdorff space $X$ via $\Gamma \times X \to X, \, (\gamma, x) \ma \gamma.x$. We form the transformation groupoid $\Gamma \ltimes X \defeq \Gamma \times X$, equipped with the product topology. The unit space of $G = \Gamma \ltimes X$ is given by $\Gn = \gekl{e} \times X \cong X$ (where $e$ is the identity of $\Gamma$), with source and range maps $\rms(\gamma,x) = x$ and $\rmr(\gamma,x) = \gamma.x$. Multiplication is given by $(\gamma', \gamma.x) (\gamma, x) \defeq (\gamma' \gamma,x)$. Such a transformation groupoid is always {\'e}tale. It is ample if and only if $X$ is totally disconnected. Moreover, $G$ is minimal if and only if $\Gamma$ acts minimally on $X$, i.e., for all $x \in X$, the orbit $\menge{\gamma.x}{\gamma \in \Gamma}$ is dense in $X$. For our transformation groupoid, $M(G)$ coincides with the $\Gamma$-invariant non-zero Radon measures on $X$.

Suppose that $\Gamma$ is countably infinite, that $X$ is compact, metrizable and totally disconnected, and that the $\Gamma$-action on $X$ is free. Then the transformation groupoid $G = \Gamma \ltimes X$ has comparison if all finitely generated subgroups of $\Gamma$ have subexponential growth \cite{DZ} (see also \cite{KS}) or if $\Gamma$ is elementary amenable \cite{KN}.

Concrete examples are given by Cantor minimal systems, i.e., the case when $\Gamma = \Zz$ and $X$ is homeomorphic to the Cantor space (see \cite{GPS95}), or by Cantor minimal $\Zz^d$-systems, i.e., the case when $\Gamma = \Zz^d$ and $X$ is homeomorphic to the Cantor space (see \cite{GMPS08,GMPS10}). A class of concrete examples is given by interval exchange transformations (see for instance \cite{CJN}). Transformation groupoids of Cantor minimal dihedral systems also have comparison by \cite{OS}. 

Another concrete class of examples is given by odometers: Let $\Gamma_i$ be a decreasing sequence of finite index subgroups of a group $\Gamma$. Then the left multiplication action of $\Gamma$ on $\Gamma / \Gamma_i$ induces an action of $\Gamma$ on $X \defeq \plim_i \Gamma / \Gamma_i$. $X$ is always totally disconnected, the action is always minimal, and the corresponding transformation groupoid always has comparison. 

Furthermore, it was shown in \cite{GGKN} that transformation groupoids of amenable, minimal actions of many non-amenable groups on the Cantor space have comparison.

\subsubsection{Tiling groupoids}
\label{sss:Tiling}

Groupoids associated with tilings have been constructed in \cite{Kel}. For aperiodic, repetitive tilings with finite local complexity, the corresponding tiling groupoids are {\'e}tale, minimal, have unit spaces homeomorphic to the Cantor space, and are almost finite by \cite{IWZ}, hence have comparison.

\subsubsection{Graph groupoids}
\label{sss:GraphGPD}

Groupoids attached to graphs have been constructed in \cite{Ren,Dea} (see also \cite{KPRR,Pat,NO21a}, for instance). We will refer to these as graph groupoids, and remark that they are special cases of Deaconu-Renault groupoids. The reader will find criteria when graph groupoids are purely infinite minimal in \cite{NO21a}.

Let us describe groupoids associated with shifts of finite type (abbreviated by SFT groupoids), which are special cases of graph groupoids. Consider a shift of finite type encoded by a finite directed graph with vertices $E^0$ and edges $E^1$. Let $A$ be the corresponding adjacency matrix, i.e., $A = (A(j,i))_{j, i \in E^0}$ where $A(j,i)$ is the number of edges from $i$ to $j$. Assume that $A$ is irreducible, in the sense that for all $j, i$, there exists $n$ such that $A^n(j,i) > 0$, and that $A$ is not a permutation matrix. Let $X_A$ be the infinite path space of our graph, i.e., $X_A$ consists of infinite sequences $(x_k)_k$ such that the target of $x_{k+1}$ is the domain of $x_k$. Equip $X_A$ with the product topology. Define the one-sided shift $\sigma_A: \: X_A \to X_A$ by setting $(\sigma_A(x_k))_k = x_{k+1}$. The groupoid attached to our shift of finite type is given by
$$
 G_A \defeq \menge{(x,n,y) \in X_A \times \Zz \times X_A}{\exists \, l, m \in \Zz 	\text{ with } l, m 	\geq 0 \text{ such that } n = l-m \text{ and } \sigma_A^l(x) = \sigma_A^m(y)}.
$$
The topology of $G_A$ is generated by sets of the form
$$
 \menge{(x,l-m,y) \in G_A}{x \in U, \, y \in V, \, \sigma_A^l(x) = \sigma_A^m(y)},
$$
where $l, m \in \Zz$ with $l, m \geq 0$, and $U$, $V$ are open subspaces of $X_A$ such that $\sigma_A^l$ and $\sigma_A^m$ induce homeomorphisms
$$
 \xymatrix{
 U \ar[r]^{\hspace{-1cm} \sigma_A^l}_{\hspace{-1cm} \cong} & \sigma_A^l(U) = \sigma_A^m(V) & \ar[l]_{\hspace{1cm} \sigma_A^m}^{\hspace{1cm} \cong} V.
 }
$$
The unit space of $G_A$ is given by $\menge{(x,0,x) \in G_A}{x \in X_A}$, which is canonically homeomorphic to $X_A$. Source and range maps are given by $\rms(x,n,y) = y$, $\rmr(x,n,y) = x$, and multiplication is given by $(x,n,y)(y,n',z) = (x,n+n',z)$. In this setting, our groupoid $G_A$ is always purely infinite minimal, with unit space homeomorphic to the Cantor space. Note that compared to the convention in \cite{Mat15}, the direction of our arrows is reversed.

\subsubsection{Higher rank graph groupoids}
\label{sss:HigherRankGraphGPD}

Higher rank graphs are small categories which come with a functor to $\Nz^k$ satisfying a certain factorisation property (see \cite{KuPa}). Groupoids attached to higher rank graphs have been introduced and studied in \cite{KuPa,FMY}. These groupoids can be identified with boundary groupoids arising from left regular representations of higher rank graphs, so that \cite[Proposition~5.21]{Li21a} gives sufficient conditions when these groupoids are purely infinite minimal. Particular examples are given by products of SFT groupoids, which are always purely infinite minimal and have been studied in \cite{Mat16}.

\subsubsection{Groupoids arising from self-similar actions}
\label{sss:SelfSimGPD}

Groupoids associated with self-similar actions on trees have been introduced and studied in \cite[Example~6.5]{Nek18a} (see also \cite{NO21b}, for instance). These are always {\'e}tale and purely infinite minimal, with unit space homeomorphic to the Cantor space. Note, however, that these groupoids may be non-Hausdorff.

{\'E}tale groupoids attached to self-similar actions on graphs have been studied in \cite{EP} (see also \cite{NO21b,Ort}). A special case is given by Katsura-Exel-Pardo groupoids $G_{A,B}$ (in the language of \cite{NO21b}), where $A$ and $B$ are two $N \times N$ row-finite matrices with integer entries, where $N \in \Nz \cup \gekl{\infty}$, and all entries of $A$ are non-negative. If $A$ is irreducible and not a permutation matrix, then these groupoids are purely infinite minimal, with unit space homeomorphic to the Cantor space.

\subsubsection{Groupoids arising from piecewise affine transformations}
\label{sss:AffineGPD}

For fixed $\lambda \in (0,1)$, groupoids arising from piecewise affine transformations on the unit interval, which on sub-intervals of the form $[a,b)$, for $a, b \in \Zz[\lambda,\lambda^{-1}]$, are given by $t \ma \lambda^i t + c$ for some $i \in \Zz$ and $c \in \Zz[\lambda,\lambda^{-1}]$, have been studied in \cite{Li15} (where they are denoted by $G \ltimes O_{P \subseteq G} \vert_{N(P)}^{N(P)}$). These groupoids are {\'e}tale, minimal, with unit space homeomorphic to the Cantor space, and a similar argument as in \cite[Proposition~4.1]{Li15} shows that the groupoids are purely infinite.

\bremark
The groupoids in \S~\ref{sss:GraphGPD}, \S~\ref{sss:HigherRankGraphGPD} and \S~\ref{sss:SelfSimGPD} are special cases of boundary groupoids arising from left regular representations of left cancellative small categories (see 	\cite{Spi20,Li21a}). Actually, in all these cases, the underlying categories have natural Garside structures, which allow for a detailed analysis of structural properties of the corresponding groupoids (see \cite{Li21b}).
\eremark

\subsection{Groupoid homology}

Let us discuss groupoid homology in the general setting of non-Hausdorff groupoids. We refer the reader to \cite{CM,Mat12} for more information about groupoid homology.

\subsubsection{Functions with compact open support and definition of groupoid homology}
\label{sss:C,GPDHom}

Let $Z$ be a topological space and $\crO$ a family of subspaces $O \subseteq Z$ which are Hausdorff, open, locally compact and totally disconnected in the subspace topology, such that $Z = \bigcup_{O \in \crO} O$. This implies that $\crO$ determines the topology of $Z$ because a subset of $Z$ is open if and only if its intersection with every $O \in \crO$ is open.

Let $\sfC$ be a $\Zz$-module, i.e., an abelian group. Given $c \in \sfC$ and a subset $U \subseteq Z$, let $c_U$ denote the function $Z \to \sfC$ with $c_U \equiv c$ on $U$ and $c_U \equiv 0$ on $Z \setminus U$. Define
$$
 \C(Z,\sfC) \defeq \lspan \menge{c_U}{U \text{ compact open subspace of some } O \in \crO, \, c \in \sfC}.
$$
By construction, $\C(Z,\sfC)$ consists of functions $Z \to \sfC$. Clearly, $\C(Z,\sfC)$ is an abelian group. As observed in \cite[Proposition~4.3]{Ste10}, $\C(Z,\sfC)$ is also the linear span of all functions of the form $c_K$ where $c$ runs through all $c \in \sfC$ and $K$ runs through all subspaces of $Z$ which are compact, open and Hausdorff.

If $Z$ is Hausdorff, then $\C(Z,\sfC)$ is the set of continuous $\sfC$-valued functions on $Z$ with compact (open) support. In that case, disjointification is a key technique in the analysis of $\C(Z,\sfC)$. In the non-Hausdorff setting, disjointification is not possible in general because intersections of compact sets might not be compact. Instead, the result below (Lemma~\ref{lem:C=Sum/I}) serves as a replacement. We include it because similar proof techniques will appear frequently in the non-Hausdorff setting.

Given $U \subseteq Z$, let $\sfC_U \defeq \menge{c_U}{c \in \sfC}$. Consider $\bigoplus_U \sfC_U$, where the sum runs over all compact open subsets $U$ of some $O \in \crO$, and let $\I$ be the subgroup of $\bigoplus_U \sfC_U$ generated by elements of the form $c_{U \amalg V} - c_U - c_V$, where $U, V$ are disjoint compact open subspaces of some $O \in \crO$.

\blemma
\label{lem:C=Sum/I}
The kernel of the canonical projection map $\pi: \: \bigoplus_U \sfC_U \to \C(Z,\sfC), \, c_U \ma c_U$ coincides with $\I$.
\elemma
\setlength{\parindent}{0cm} \setlength{\parskip}{0cm}

\bproof
Suppose that $f = \sum_{i \in I} (c_i)_{U_i}$ satisfies $\pi(f) = 0$, where $I$ is a finite index set. Suppose that $\gekl{O_1, \dotsc, O_n}$ is a finite subset of $\crO$ such that for every $i \in I$ there exists $1 \leq m \leq n$ with $U_i \subseteq O_m$. We proceed inductively on $n$. 
\setlength{\parindent}{0cm} \setlength{\parskip}{0.5cm}

If $n = 1$, then all $U_i$ are contained in some $O \in \crO$. Then we can disjointify $U_i$ in $O$, i.e., we let $\gekl{V_j}$ be the set of non-empty subspaces of $O$ of the form $\bigcap_{i \in I'} U_i \cap \bigcap_{i' \in I \setminus I'} U_i^c$, where $I'$ runs through all non-empty subsets of $I$ and $U_i^c = Z \setminus U_i$. By construction, $\gekl{V_j}$ is a family of pairwise disjoint subsets, and because $O$ is Hausdorff, every $V_j$ is compact open. Moreover, every $U_i$ is a disjoint union of $V_j$ because $U_i = \bigcup_{i \in I'} (\bigcap_{i \in I'} U_i \cap \bigcap_{i' \in I \setminus I'} U_i^c)$. So we can write $U_i = \coprod_{j_i} V_{j_i}$. Hence it follows that $(c_i)_{U_i} \equiv \sum_{j_i} (c_i)_{V_{j_i}} \ {\rm mod} \ \I$. Hence we obtain $f \equiv \sum_j (\ti{c}_j)_{V_j}  \ {\rm mod} \ \I$. But now $\pi(f) = 0$ implies that every $\ti{c}_j$ must be zero because the $V_j$ are pairwise disjoint. Hence $f \equiv \sum_j (\ti{c}_j)_{V_j} = 0 \ {\rm mod} \ \I$, i.e., $f \in \I$.

Now suppose that $n > 1$. By disjointifying, we may assume that all the $U_i$ which are contained in a single $O_m$ are pairwise disjoint. Now fix $i \in I$ with $U_i \subseteq O_n$. If there exists $z \in Z$ with $z \in U_i$, $z \notin U_{i'}$ for all $i' \neq i$, then $c_i = 0$ as $0 = \pi(f)(z) = c_i$. Hence we may assume $U_i \subseteq \bigcup_{i \neq i' \in I} U_{i'}$. Actually, we even have $U_i \subseteq \bigcup_{i' \in I'} U_{i'}$ where $I'$ is a subset of $I \setminus \gekl{i}$ such that for every $i' \in I'$, we have $U_{i'} \subseteq O_m$ for some $1 \leq m < n$. Here we are using that all $U_{i'}$ with $U_{i'} \subseteq O_n$ are pairwise disjoint. Therefore, for every $x \in U_i$ there exists a compact open neighbourhood $V_x$ of $x$ with $V_x \subseteq U_{i'}$ for some $i' \in I'$. As $U_i$ is compact, we can write $U_i$ as a finite union $U_i = \bigcup_{j \in J} V_j$, where $V_j = V_{x_j}$ for some $x_j \in U_i$. As $U_i$ is contained in $O_n$, all $V_j$ are also contained in $O_n$, so that we can (after disjointifying) assume that the $V_j$ are pairwise disjoint. Thus $(c_i)_{U_i} \equiv \sum_{j \in J} (c_i)_{V_j} \ {\rm mod} \ \I$ and thus $f \equiv \sum_{j \in J} (c_i)_{V_j} + \sum_{i \neq i' \in I} (c_{i'})_{U_{i'}} \ {\rm mod} \ \I$. Now run this procedure for all $i \in I$ such that $U_i \subseteq O_n$. In this way, we are able to replace $\gekl{O_1, \dotsc, O_n}$ by $\gekl{O_1, \dotsc, O_{n-1}}$ and then apply induction hypothesis.
\eproof
\setlength{\parindent}{0cm} \setlength{\parskip}{0.5cm}

We derive the following immediate consequence.
\bcor
$\C(Z,\sfC) \cong \C(Z,\Zz) \otimes_{\Zz} \sfC$.
\ecor

Let $G$ be an ample groupoid with locally compact Hausdorff unit space $\Gn$. Let $\sfC$ be as above. Define 
$$
 \C(G,\Zz) \defeq \lspan \menge{1_U}{U \subseteq G \text{ compact open bisection}}.
$$
Note that since every compact open bisection is automatically Hausdorff, $\C(G,\Zz)$ coincides with $\C(Z,\Zz)$ as defined above, for $Z = G$ and $\crO$ given by the collection of all open bisections. $\C(G,\Zz)$ becomes a $\Zz$-algebra with respect to convolution given by $(f_1 f_2) (g) = \sum_{h_1 h_2 = g} f_1(h_1) f_2(h_2)$ for $f_1, f_2 \in \C(G,\Zz)$. Algebras of this form have for instance been studied in \cite{Ste10, CEPSS}.

Now consider 
$$
 \C(\Gn,\sfC) \defeq \lspan \menge{c_U}{U \subseteq \Gn \text{ compact open}, \, c \in \sfC}.
$$
$\C(\Gn,\sfC)$ is a left- and right-$\C(G,\Zz)$-module via
$$
 (f m) (x) = \sum_{g^{-1} g = x} f(g^{-1}) m(\rmr(g)), \qquad (m f) (x) = \sum_{g g^{-1} = x} m(\rms(g)) f(g^{-1}),
$$
for $f \in \C(G,\Zz)$ and $m \in \C(\Gn,\sfC)$. Let us now define groupoid homology in terms of the bar resolution, and then explain an alternative approach using $\Tor$.

Let $\Gnu \defeq \menge{(g_1, \dotsc, g_{\nu}) \in G^{\nu}}{\rms(g_{\mu + 1}) = \rmr(g_{\mu})}$, equipped with the subspace topology coming from the product topology on $G^{\nu}$. Let $\crO^{(\nu)}$ be the collection of subsets of $\Gnu$ of the form 
$$
 O_1 \tensor[_{\rms}]{\times}{_{\rmr}} \dotso \tensor[_{\rms}]{\times}{_{\rmr}} O_{\nu} \defeq \menge{(g_1, \dotsc, g_{\nu}) \in \Gnu}{g_{\mu} \in O_{\mu}},
$$
where $O_{\mu}$ are open bisections with $\rms(O_{\mu + 1}) = \rmr(O_{\mu})$. Let $\C(\Gnu,\sfC)$ be defined as above (with $Z = \Gnu$, $\crO = \crO^{(\nu)}$). Consider the maps $\ti{d}_{\nu}^{\mu}: \: \Gnu \to \Gnum$ given by $\ti{d}_1^0 = \rms$, $\ti{d}_1^1 = \rmr$ and 
$$
 \ti{d}_{\nu}^{\mu}(g_1, \dotsc, g_{\nu}) 
 =
 \bfa
 (g_2, \dotsc, g_{\nu}) & \falls \mu = 0,\\
 (g_1, \dotsc, g_{\mu} g_{\mu + 1}, \dotsc, g_{\nu}) & \falls 0 < \mu < \nu,\\
 (g_1, \dotsc, g_{\nu - 1}) & \falls \mu = \nu.
 \efa
$$
Since $\ti{d}_{\nu}^{\mu}$ are local homeomorphisms, they induce homomorphisms $(\ti{d}_{\nu}^{\mu})_*: \: \C(\Gnu,\sfC) \to \C(\Gnum,\sfC)$ given by $(\ti{d}_{\nu}^{\mu})_*(f)(z) = \sum_{y \ \in \ (d_{\nu}^{\mu})^{-1}(z)} f(y)$. Now define 
\begin{equation}
\label{e:tipartial}
\ti{\partial}_{\nu} \defeq \sum_{\mu = 0}^{\nu} (-1)^{\mu} (\ti{d}_{\nu}^{\mu})_*.
\end{equation}
It is straightforward to check that $B_*(G,\sfC) \defeq (\C(\Gnu,\sfC), \ti{\partial}_{\nu})_{\nu}$ is a chain complex. Groupoid homology is defined as the homology of this chain complex, i.e.,
$$
 H_*(G,\sfC) \defeq H_*(B_*(G,\sfC)).
$$

\subsubsection{Groupoid homology for examples}
\label{sss:GPDHomEx}

Given an AF groupoid $G$ as in \S~\ref{sss:AF}, the $0$-th homology $H_0(G,\sfC)$ is given by the dimension group, with coefficients in $\sfC$, of a Bratteli diagram describing $G$ (the dimension group is independent of the choice of the diagram), and all higher homology groups vanish, i.e., $H_*(G,\sfC) \cong \gekl{0}$ for all $* > 0$ (see for instance \cite{Ren, Kri, FKPS, Mat12}).

For a transformation groupoid $G = \Gamma \ltimes X$ as in \S~\ref{sss:Trafo}, it follows from the definitions that groupoid homology is canonically isomorphic to group homology with coefficients in the $\Gamma$-module $C_c(X,\sfC)$, i.e., $H_*(G,\sfC) \cong H_*(\Gamma,C_c(X,\sfC))$ (see \cite{Bro} for the definition of group homology). Here $C_c(X,\sfC)$ denotes the set of compactly supported continuous functions on $X$ with values in $\sfC$, where $\sfC$ is equipped with the discrete topology. Note that since $X$ is Hausdorff, $C_c(X,\sfC)$ coincides with $\C(X,\sfC)$ from \S~\ref{sss:C,GPDHom}. Let us describe groupoid homology more explicitly in the case of Cantor minimal systems, following \cite{GPS95}. In that case, $\Gamma = \Zz$, $X$ is homeomorphic to the Cantor space, and $C_c(X,\sfC) = C(X,\sfC)$ because $X$ is compact. Let $\varphi \in {\rm Homeo}(X)$ be the homeomorphism corresponding to the canonical generator $1$ of $\Zz$. Then
$$
 H_0(\Zz \ltimes X,\sfC) \cong C(X,\sfC) / \menge{f - f \circ \varphi^{-1}}{f \in C(X,\sfC)},
$$
$H_1(\Zz \ltimes X,\sfC) \cong \sfC$ and all higher homology groups vanish, i.e., $H_*(\Zz \ltimes X,\sfC) \cong \gekl{0}$ for all $*>1$.

Let us now consider tiling groupoids as in \S~\ref{sss:Tiling}. Given a tiling of $\Rz^d$, let $G$ be its tiling groupoid and $\Omega$ the hull space of our tiling. As observed in \cite[\S~5.2]{PY}, groupoid cohomology of $G$ can be identified with sheaf cohomology of $\Omega$. Using the description of groupoid homology of $G$ in terms of group homology (see for instance \cite[\S~5.2]{PY}) and Poincar{\'e} duality, we obtain an identification of groupoid homology $H_*(G)$ with the $(d-*)$-th {\v{C}}ech cohomology $\check{H}^{d-*}(\Omega)$ of $\Omega$. Explicit homology computations can be found in \cite{GK,FHK01,FHK02,GHK13}.

For an SFT groupoid $G_A$ as in \S~\ref{sss:GraphGPD}, it was shown in \cite[Theorem~4.14]{Mat12} that
$$
 H_*(G_A,\sfC) \cong
 \bfa
 \coker(\id - A^t: \: \bigoplus_{E^0} \sfC \to \bigoplus_{E^0} \sfC) & \falls *=0,\\
 \ker(\id - A^t: \: \bigoplus_{E^0} \sfC \to \bigoplus_{E^0} \sfC) & \falls *=1,\\
 \gekl{0} & \sonst.
 \efa
$$
For homology computations for more general graph groupoids, we refer to \cite{NO21a} and the references therein.

For groupoids of higher rank graphs as in \S~\ref{sss:HigherRankGraphGPD}, groupoid homology has been computed for some cases in \cite{FKPS}. Let us briefly summarize the result from \cite{FKPS} in the one vertex case. Let $\Lambda$ be a one vertex $k$-graph and $G_{\Lambda}$ the corresponding groupoid. Let $\Lambda^{\varepsilon_i}$ be the elements of $\Lambda$ with degree $\varepsilon_i$, where $\varepsilon_i$ are the standard generators of $\Nz^k$. Write $N_i \defeq \abs{\Lambda^{\epsilon_i}} - 1$. If $\Lambda$ is row-finite and $N_i \geq 1$ for all $i$, then
$$
 H_*(G_{\Lambda}) \cong
 \bfa
 (\Zz/\gcd(N_1, \dotsc, N_k) )^{\binom{k-1}{*}} & \falls 0 \leq * \leq k-1,\\
 \gekl{0} & \sonst.
 \efa
$$
For products of SFT groupoids, which is another particular case of groupoids of higher rank graphs, a complete computation of groupoid homology has been established in \cite[Proposition~5.4]{Mat16}.

Consider Katsura-Exel-Pardo groupoids $G_{A,B}$, which are special cases of groupoids attached to self-similar actions on graphs (see \S~\ref{sss:SelfSimGPD}). Let us present the groupoid homology computation in \cite{Ort} (see also \cite{NO21b}). We use the same notation as in \S~\ref{sss:SelfSimGPD}. Assume that $A$ and $B$ are row-finite matrices with integer entries, and all entries of $A$ are non-negative. Suppose that for all $1 \leq i, j \leq N$, $B_{i,j} = 0$ if and only if $A_{i,j} = 0$. Then
$$
 H_*(G_{A,B}) \cong
 \bfa
 \coker(\id - A) & \falls * = 0,\\
 \ker(\id - A) \oplus \coker(\id - B) & \falls * = 1,\\
 \ker(\id - B) & \falls * = 2,\\
 \gekl{0} & \sonst.
 \efa
$$

For groupoids arising from piecewise affine transformations as in \S~\ref{sss:AffineGPD}, groupoid homology computations for classes of examples can be found in \cite{Li15}.

\subsubsection{Description of groupoid homology using derived functors}
\label{sss:GPDHom=Tor}

Our goal now is to show that 
$$
 H_*(G,\sfC) \cong \Tor_*^{\C(G,\Zz)}(\C(\Gn,\Zz),\C(\Gn,\sfC)).
$$
For Hausdorff groupoids, this is shown in \cite{BDGW,Mil}. We treat the case of non-Hausdorff groupoids. The results of \S~\ref{sss:GPDHom=Tor} are merely included for completeness; they are not needed in the sequel. However, some of the ideas will appear again in \S~\ref{s:HK=HG}.

First of all, note that, following for instance the approach in \cite[\S~4.1]{Mil}, we will be able to use standard results in homological algebra, even though they are usually formulated for unital rings whereas our ring $\C(G,\Zz)$ is in general not unital, only locally unital. 

First of all, the inversion map induces an involution on $\C(G,\Zz)$ which flips the order of multiplication, which in turn allows us to interchange left-$\C(G,\Zz)$-modules and right-$\C(G,\Zz)$-modules and thus leads to the identification $ \Tor_*^{\C(G,\Zz)}(\C(\Gn,\Zz),\C(\Gn,\sfC)) \cong  \Tor_*^{\C(G,\Zz)}(\C(\Gn,\sfC),\C(\Gn,\Zz))$. So it suffices to show that $H_*(G,\sfC) \cong \Tor_*^{\C(G,\Zz)}(\C(\Gn,\sfC),\C(\Gn,\Zz))$. 

Next, we define another chain complex $E_*(G,\Zz) \defeq (\C(\Gnup,\Zz),\partial_{\nu + 1})_{\nu}$, where $\partial_{\nu + 1}$ is given as follows: The maps $d_{\nu + 1}^{\mu}: \: \Gnup \to \Gnu$ given by
$$
 d_{\nu + 1}^{\mu}(g_0, \dotsc, g_{\nu}) 
 =
 \bfa
 (g_0, \dotsc, g_{\mu} g_{\mu + 1}, \dotsc, g_{\nu}) & \falls 0 \leq \mu < \nu,\\
 (g_0, \dotsc, g_{\nu - 1}) & \falls \mu = \nu,
 \efa
$$
are local homeomorphisms, hence induce homomorphisms $(d_{\nu + 1}^{\mu})_*: \: \C(\Gnup,\Zz) \to \C(\Gnu,\Zz)$ given by $(d_{\nu + 1}^{\mu})_*(f)(z) = \sum_{y \ \in \ (d_{\nu + 1}^{\mu})^{-1}(z)} f(y)$. Define $\partial_{\nu + 1} \defeq \sum_{\mu = 0}^{\nu} (d_{\nu + 1}^{\mu})_*$. Now consider the left $G$-action on $\Gnup$ with respect to the anchor map $\rho: \: \Gnup \to \Gn, \, (g_0, \dotsc, g_{\nu}) \ma \rmr(g_0)$ and the action $g.(g_0, \dotsc, g_{\nu}) \defeq (gg_0, g_1, \dotsc, g_{\nu})$ for all $g \in G$ and $(g_0, \dotsc, g_{\nu}) \in \Gnup$ with $\rho(g_0, \dotsc, g_{\nu}) = \rms(g)$. This $G$-action induces a left-$\C(G,\Zz)$-module structure on $\C(\Gnup,\Zz)$ via $(fm)(z) \defeq \sum_{g \in G, \, y \in \Gnup, \, g.y \, = \, z} f(g) m(y)$ for $f \in \C(G,\Zz)$, $m \in \C(\Gnup,\Zz)$. It is straightforward to check that $E_*(G,\Zz)$ is a chain complex in the category of $\C(G,\Zz)$-modules.

Observe that $B_*(G,\sfC) \cong \C(\Gn,\sfC) \otimes_{\C(G,\Zz)} E_*(G,\sfC)$. This is because the identification $\Gn \times_G \Gnup \cong \Gnu, \, (\rmr(g_0), (g_0, \dotsc, g_{\nu}))) \ma (g_1, \dotsc, g_{\nu})$ induces an isomorphism 
$$
 \C(\Gn,\sfC) \otimes_{\C(G,\Zz)} \C(\Gnup,\Zz) \cong \C(\Gnu,\sfC)
$$
sending $c_{\rmr(U_0)} \otimes 1_{U_0 \tensor[_{\rms}]{\times}{_{\rmr}} \dotso \tensor[_{\rms}]{\times}{_{\rmr}} U_{\nu}}$ to $c_{U_1 \tensor[_{\rms}]{\times}{_{\rmr}} \dotso \tensor[_{\rms}]{\times}{_{\rmr}} U_{\nu}}$.

Moreover, $E_*(G,\sfC)$ is exact. The corresponding chain homotopy is induced by the maps $h_{\nu}: \: \Gnu \to \Gnup, \, (g_0, \dotsc, g_{\nu - 1}) \ma (\rmr(g_0), g_0, \dotsc, g_{\nu - 1})$ for $\nu \geq 1$ and the inclusion $h_0: \: \Gn \to G$ for $\nu = 0$ (see for instance \cite{BDGW,Mil}, but note that $EG_{\bullet}$ in \cite[\S~2.3]{BDGW} does not coincide with our $G^{(\bullet)}$; instead, use the identification $EG_{\bullet} \ni (g_0, g_1, g_2, \dotsc) \ma (g_0, g_0^{-1} g_1, g_1^{-1} g_2, \dotsc) \in G^{(\bullet)}$).

Therefore, once we show that $\C(\Gnup,\Zz)$ are flat left-$\C(G,\Zz)$-modules, then we conclude that $H_*(G,\sfC) \cong \Tor_*^{\C(G,\Zz)}(\C(\Gn,\Zz),\C(\Gn,\sfC))$.

As a first step, observe that we have an isomorphism
\begin{equation}
\label{e:CxCnu}
 \C(G,\Zz) \otimes_{\C(\Gn,\Zz)} \C(\Gnu,\Zz) \cong \C(\Gnup,\Zz)
\end{equation}
sending $a \otimes f$ to the function $(g_0, g_1, \dotsc, g_{\nu}) \ma a(g_0) f(g_1, \dotsc, g_{\nu})$. The proof is similar as the one for Lemma~\ref{lem:C=Sum/I}.

Now suppose that $M$ is a right-$\C(G,\Zz)$-module. By \cite{Ste14}, there is a sheaf $\cM$ of $\Zz$-modules over $\Gn$ together with a $G$-action such that $M \cong \Gamma_c(\Gn,\cM)$ as right-$\C(G,\Zz)$-modules. Here $\Gamma_c$ stands for continuous sections with compact support. Using the anchor map $\rho: \: \Gnu \to \Gn$, define the pullback $\rho^* \cM$ as a sheaf over $\Gnu$ with fibre $(\rho^* \cM)_z = \cM_{\rho(z)}$ for $z \in \Gnu$. Let $\crO^{(\nu)}$ be the collection of subspaces of the form $O_1 \tensor[_{\rms}]{\times}{_{\rmr}} \dotso \tensor[_{\rms}]{\times}{_{\rmr}} O_{\nu}$ as above. Note that $\rho$ restricts to a homeomorphism on these subspaces $O_1 \tensor[_{\rms}]{\times}{_{\rmr}} \dotso \tensor[_{\rms}]{\times}{_{\rmr}} O_{\nu}$. For a compact open subspace $U \subseteq O \in \crO^{(\nu)}$ and $m \in M$ define $(\rho^*m)_U$ as the composite
$$
 \xymatrix{
 (\rho \vert_O)^{-1}(U) \ar[r]^{\hspace{0.5cm} \rho} & \rho(U) \ar[r]^{m 1_{\rho(U)}} & \cM.
 }
$$
Set
$$
 \Gamma_{\C}(\Gnu,\rho^*\cM) \defeq \lspan \menge{(\rho^*m)_U}{U \text{ compact open subspace of some } O \in \crO^{(\nu)}, \, m \in M}.
$$
By construction, $\Gamma_{\C}$ consists of sections $\Gnu \to \rho^*\cM$. Note that since $\Gnu$ is not Hausdorff in general, $\Gamma_{\C}$ does not coincide with $\Gamma_c$. Now a similar argument as for Lemma~\ref{lem:C=Sum/I} implies that the following map,
\begin{equation}
\label{e:MxCnu}
 M \otimes_{\C(\Gn,\Zz)} \C(\Gnu,\Zz) \cong \Gamma_{\C}(\Gnu,\rho^*\cM),
\end{equation}
sending $m \otimes f$ to the function $z \ma m(	\rho(z)) f(z)$, is an isomorphism.

Now we arrive at the desired conclusion.
\bprop
For all $\nu \geq 0$, $\C(\Gnup,R)$ is a flat left-$\C(G,\Zz)$-module.
\eprop
\bproof
Suppose that
$
 \xymatrix{
 0 \ar[r] & M \ar[r]^{\iota} & N
 }
$
is an exact sequence of right-$\C(G,\Zz)$-modules. By \cite{Ste14}, we obtain corresponding sheaves $\cM$ and $\cN$. Moreover, $\iota$ induces injective homomorphisms $\iota_x: \: \cM_x \into \cN_x$ on the fibres, for all $x \in \Gn$. We want to show that
$$
 \iota \otimes \id: \: M \otimes_{\C(G,\Zz)} \C(\Gnup,\Zz) \to N \otimes_{\C(G,\Zz)} \C(\Gnup,\Zz)
$$
is injective. Using \eqref{e:CxCnu} and \eqref{e:MxCnu}, we obtain the following identification:
\begin{align*}
 & M \otimes_{\C(G,\Zz)} \C(\Gnup,\Zz) 
 \cong M \otimes_{\C(G,\Zz)} \C(G,\Zz) \otimes_{\C(\Gn,\Zz)} \C(\Gnu,\Zz)
 \cong M \otimes_{\C(\Gn,\Zz)} \C(\Gnu,\Zz)\\
 \cong \ & \Gamma_{\C}(\Gnu,\rho^*\cM).
\end{align*}
Similarly, $N \otimes_{\C(G,\Zz)} \C(\Gnup,\Zz) \cong \Gamma_{\C}(\Gnu,\rho^*\cN)$. Identifying elements of $\Gamma_{\C}(\Gnu,\rho^*\cM)$ and $\Gamma_{\C}(\Gnu,\rho^*\cN)$ as functions on $\Gnu$ with values in $\cM$ and $\cN$, respectively, we see that, under the identifications above, $\iota \otimes \id$ sends $f \in \Gamma_{\C}(\Gnu,\rho^*\cM)$ to the map $z \ma \iota_{\rho(z)} (f(z))$. And since $\iota_x$ is injective for all $x \in \Gn$, we deduce that $\iota \otimes \id$ must be injective as well, as desired.
\eproof

As explained above, using \cite[Proposition~4.21]{Mil}, this leads to the desired description of groupoid homology in terms of $\Tor$.
\btheo
$H_*(G,\sfC) \cong \Tor_*^{\C(G,\Zz)}(\C(\Gn,\Zz),\C(\Gn,\sfC))$.
\etheo

\subsection{Topological full groups}
\label{ss:TFG}

In the following, let $G$ be an ample groupoid with locally compact Hausdorff unit space $\Gn$. 
\bdefin
If $\Gn$ is compact, then the topological full group $\bmF(G)$ is the group of compact open bisections $\sigma \subseteq G$ with $\rmr(\sigma) = \Gn = \rms(\sigma)$. Multiplication in $\bmF(G)$ is given by multiplication of bisections, i.e., $\sigma \tau \defeq \menge{gh}{g \in \sigma, \, h \in \tau, \, \rms(g) = \rmr(h)}$. 

In the general case where $\Gn$ is not necessarily compact, we set $\bmF(G) \defeq \ilim_U \bmF(G_U^U) = \bigcup_U \bmF(G_U^U)$. Here the limit is taken over all compact open subspaces $U \subseteq \Gn$, ordered by inclusion, and $G_U^U = \menge{g \in G}{\rms(g), \, \rmr(g) \in U}$. Given two compact open subspaces $U \subseteq V$ of $\Gn$, we view $\bmF(G_U^U)$ as a subgroup of $\bmF(G_V^V)$ via the embedding $\bmF(G_U^U) \into \bmF(G_V^V), \, \sigma \ma \sigma \amalg (V \setminus U)$.
\edefin
If $G$ is effective, i.e., when the interior of the isotropy subgroupoid $\menge{g \in G}{\rmr(g) = \rms(g)}$ coincides with $\Gn$, then the map sending $\sigma \in \bmF(G_U^U) \subseteq \bmF(G)$ to the homeomorphism $\Gn \to \Gn$ given by $x \ma \sigma.x$ on $U$ and identity on $\Gn \setminus U$ is injective, so that we may view $\bmF(G)$ as a subgroup of ${\rm Homeo}(\Gn)$. Here we use the notation that $\sigma.x$ denotes $\rmr(g)$ for the unique element $g \in \sigma$ with $\rms(g) = x$.

Topological full groups first appeared in \cite{Kri,GPS99} (for special classes of groupoids). Several subgroups of $\bmF(G)$ have been constructed, for instance the alternating full group $\bmA(G)$ (see \cite{Nek19}). It is known that for almost finite or purely infinite groupoids $G$ which are minimal, effective and Hausdorff, with unit space $\Gn$ homeomorphic to the Cantor space, the alternating full group coincides with the commutator subgroup $\bmD(G)$ of $\bmF(G)$ (see \cite{Mat15,Nek19}).

Nekrashevych showed in \cite{Nek19} that for every minimal, effective groupoid $G$ whose unit space $\Gn$ is homeomorphic to the Cantor space, the alternating full group $\bmA(G)$ is simple. Moreover, again for minimal, effective groupoids $G$ with unit space $\Gn$ homeomorphic to the Cantor space, it is possible to reconstruct the groupoid $G$ from the topological full group $\bmF(G)$ (see \cite{Mat15,Nek19}). A far-reaching generalization of these results has been obtained in \cite{MB}.

Matui formulated the AH-conjecture, which describes the first homology group $H_1(\bmF(G))$ in terms of groupoid homology of $G$. He constructed an index map $I: \: H_1(\bmF(G)) \to H_1(G)$ and conjectured for every minimal, effective groupoid $G$ whose unit space $\Gn$ is homeomorphic to the Cantor space, there is an exact sequence
$$
 \xymatrix{
 H_0(G) \otimes \Zz/2 \ar[r] & H_1(\bmF(G)) \ar[r]^{\hspace{0.25cm} I} & H_1(G) \ar[r] & 0.
 }
$$
Note that Matui restricts his discussion to second countable Hausdorff groupoids. Moreover, he formulated the AH-conjecture in terms of the abelianization $\bmF(G)^{\rm ab}$ of $\bmF(G)$, which is isomorphic to $H_1(\bmF(G))$. The AH-conjecture has been verified for all principal, almost finite, second countable Hausdorff groupoids as well as groupoids arising from shifts of finite type, products of groupoids from shifts of finite type, graph groupoids, Katsura-Exel-Pardo groupoids as well as transformation groupoids of of odometers and Cantor minimal dihedral systems (see \cite{Mat12,Mat15,Mat16,NO21a,NO21b,Sca19,Sca22}).

\subsection{Examples of topological full groups}
\label{ss:ExTopFullGrp}

For an AF groupoid $G$, the topological full group $\bmF(G)$ is the increasing union of finite direct sums of finite symmetric groups.

For a transformation groupoid $G$ of a Cantor minimal $\Zz$-system, Juschenko and Monod showed that $\bmF(G)$ is amenable \cite{JM}. By taking alternating full groups $\bmA(G)$ in the case of minimal subshifts (in this case $\bmA(G)$ coincides with the commutator subgroup $\bmD(G)$ of $\bmF(G)$), this produces the first examples of finitely generated amenable infinite simple groups, answering an open problem in group theory.

Now let $D_{\infty}$ be the infinite dihedral group $\Zz \rtimes (\Zz / 2) \cong (\Zz / 2) * (\Zz / 2)$. Starting with a Cantor minimal $D_{\infty}$-system, Nekrashevych constructed another Cantor minimal $\Gamma$-system (for some new group $\Gamma$) and shows that, under certain conditions, the alternating full group of the groupoid of germs for the $\Gamma$-action is a finitely generated simple periodic group of intermediate growth (see \cite{Nek18b}). This again answers an open problem in group theory.

Other examples include topological full groups of transformation groupoids of Cantor minimal $\Zz^d$-systems, tiling groupoids or topological full groups of interval exchange transformations (the latter have been studied in \cite{CJN}).

The groupoids discussed so far are almost finite (with the exception of Nekrashevych's examples, where the groupoids could be non-Hausdorff and almost finiteness is not known). Let us now turn to topological full groups of purely infinite groupoids.

Let $G_2$ be the groupoid attached to the one-sided full shift on two symbols. In the language of \S~\ref{sss:GraphGPD}, the graph we consider consists of one vertex and two edges (which must then be loops), the corresponding adjacency matrix $A$ is given by the $1 \times 1$-matrix with entry $2$, and we set $G_2 \defeq G_A$. Then $\bmF(G_2) \cong V$, where $V$ is Thompson's group $V$ (see \cite{CFP}). This was first observed in \cite{Nek04} (see also \cite{Mat15}). More generally, if we consider the one-sided full shift on $k$ symbols, its graph given by one vertex and $k$ edges, the adjacency matrix given by the $1 \times 1$-matrix with entry $k$ and let the corresponding groupoid be $G_k$, then $\bmF(\cR_r \times G_k^n) \cong n V_{k,r}$. Here $\cR_r$ is the groupoid given by the full equivalence relation on the finite set $\gekl{1, \dotsc, r}$, and $n V_{k,r}$ are the Brin-Higman-Thompson groups (see for instance \cite{Hig,Bri04} and also \cite{Mat16}). 

Szymik and Wahl show that the group homology $H_*(V_{k,r})$ does not depend on $r$, and that $V_{k,r}$ is acyclic for $k=2$. They also produce further computations of parts of $H_*(V_{k,r})$ (see \cite{SW}). Their work is based on Cantor algebras, which lead to the construction of a small permutative category and hence an algebraic K-theory spectrum $\Kz$ such that $H_*(V_{n,r}) \cong H_*(\Omega^{\infty}_0 \Kz)$. Here $\Omega^{\infty} \Kz$ is the infinite loop space corresponding to $\Kz$, and $\Omega^{\infty}_0 \Kz$ denotes the path component of the base point of $\Omega^{\infty} \Kz$ (these notions are introduced in \S~\ref{ss:AlgKSmallPermCat}). The results in \cite{SW} are then obtained by analysing the homotopy groups of $\Kz$. It is not immediate how to carry over the constructions in \cite{SW} to more general groupoids because the notion of Cantor algebras is tailored to the situation of Higman-Thompson groups.

Topological full groups for SFT groupoids and products of SFT groupoids have been studied in detail in \cite{Mat15,Mat16}.

For groupoids $G$ arising from self-similar actions on trees, the topological full groups $\bmF(G)$ are isomorphic to R{\"o}ver-Nekrashevych groups (see \cite{Rov,Nek04,Nek18a} as well as \cite{SWZ}). These groups have interesting finiteness properties. We say that a group is of type ${\rm F}_n$ if it admits a classifying space with a compact $n$-skeleton. These finiteness properties play an important role in group homology and specialize to familiar notions in low dimensions (a group is of type ${\rm F}_1$ if and only if it is finitely generated and of type ${\rm F}_2$ if and only if it is finitely presented). \cite{SWZ} shows that topological full groups arising from classes of self-similar actions give rise to first examples of infinite simple groups which are of type ${\rm F}_{n-1}$ but not of type ${\rm F}_n$, for each $n$.

For groupoids $G$ arising from piecewise affine transformations \cite{Li15}, the topological full groups $\bmF(G)$ are isomorphic to groups considered in \cite{Ste92} (for the parameters $l=1$, $A = \Zz[\lambda,\lambda^{-1}]$, $P=\spkl{\lambda}$ in the terminology of \cite{Ste92}). Moreover, a similar construction as in \cite{Li15} leads to ample groupoids whose topological full groups coincide with the groups denoted by $G(l,A,P)$ in \cite{Ste92}, for arbitrary parameters $l, A, P$ (see \cite{Tan}).

As these examples show, topological full groups and the closely related notion of alternating full groups lead to new examples in group theory with interesting properties. They provide a rich supply of infinite simple groups. However, even though we have seen much progress regarding particular example classes and spectacular advances have been made in our understanding of these groups, general results about topological full groups are rare and seem to be difficult to obtain. For instance, not much is known regarding analytic properties of topological full groups in general. All in all, it is a fascinating yet challenging problem to develop a better understanding of the interplay between group-theoretic properties of topological full groups and dynamical properties of the underlying topological groupoids.

\subsection{Algebraic K-theory spectra of small permutative categories}
\label{ss:AlgKSmallPermCat}

Let us now describe the construction of algebraic K-theory spectra from small permutative categories as in \cite{Seg,Tho}. We follow the exposition in \cite{EM}. We will use the language of simplicial sets (see for instance \cite{GZ,GJ}) and of spectra in the sense of algebraic topology (see for instance \cite{Swi,Ada78}). 

\bdefin
A small permutative category is a small category $\fB$ with object set $\obj \fB$, morphism set $\mor \fB$, together with a functor $\oplus: \: \fB \times \fB \to \fB$, an object $0 \in \obj \fB$ and natural isomorphisms
$$
 \gekl{\pi_{u,u'} \in \mor \fB: \: \pi_{u,u'}: \: u \oplus u' \cong u' \oplus u}_{u,u' \, \in \, \obj \fB},
$$
such that $\oplus$ is associative with unit $0$, and we have $\pi_{0,u} = \pi_{u,0} = \id_u$, $\pi_{u',u} \pi_{u,u'} = \id_{u \oplus u'}$, i.e., the diagram
$$
 \xymatrix{
 u \oplus u' \ar[dr]_{\id_{u \oplus u'}} \ar[r]^{\pi_{u,u'}} & u' \oplus u \ar[d]^{\pi_{u',u}}\\
 & u \oplus u'
 }
$$
commutes, and $(\pi_{u,u''} \oplus \id_{u'}) (\id_u \oplus \pi_{u',u''}) = \pi_{u \oplus u',u''}$, i.e., the diagram
$$
 \xymatrix{
 u \oplus u' \oplus u'' \ar[dr]_{\pi_{u \oplus u',u''}} \ar[r]^{\id \oplus \pi_{u',u''}} & u \oplus u'' \oplus u' \ar[d]^{\pi_{u,u''} \oplus \id_{u'}}\\
 & u'' \oplus u \oplus u'
 }
$$
commutes, for all objects $u, u', u''$ of $\fB$.
\edefin
Given $\sigma \in \mor \fB$, we write $\mft(\sigma)$ for its target and $\mfd(\sigma)$ for its domain, and we denote by $\fB(v,u)$ the set of morphisms of $\fB$ from $u$ to $v$, i.e., $\fB(v,u) = \menge{\sigma \in \mor \fB}{\mft(\sigma) = v, \, \mfd(\sigma) = u}$.

Now let $A$ be a finite based set, i.e., a finite set with a choice of an element called the base point.
\bdefin
Given a small permutative category $\fB$ and a finite based set $A$, let $\fB(A)$ be the category with objects of the form $\gekl{u_S, \varphi_{T,T'}}_{S,T,T'}$, where $S, T, T'$ run through all subsets of $A$ not containing the base point with $T \cap T' = \emptyset$, $u_S \in \obj \fB$ for all $S$ and $\varphi_{T,T'} \in \fB(u_{T \cup T'}, u_T \oplus u_{T'})$ are isomorphisms for all $T, T'$. We require that for $S = \emptyset$, $u_{\emptyset} = 0$, and for $T = \emptyset$, $\varphi_{\emptyset,T'} = \id_{u_{T'}}$. Moreover, for all pairwise disjoint $T, T', T''$, the following diagrams should commute:
$$
 \xymatrix{
 u_T \oplus u_{T'} \ar[d]_{\pi_{u_T,u_{T'}}} \ar[r]^{\varphi_{T,T'}} & u_{T \cup T'} \ar[d]^{\id}\\
 u_{T'} \oplus u_T \ar[r]_{\varphi_{T',T}} & u_{T' \cup T}
 }
$$
$$
 \xymatrix{
 u_T \oplus u_{T'} \oplus u_{T''} \ar[d]_{\id_{u_T} \oplus \varphi_{T',T''}} \ar[r]^{\varphi_{T,T'} \oplus \id_{u_{T''}}} & u_{T \cup T'} \oplus u_{T''} \ar[d]^{\varphi_{T \cup T', T''}}\\
 u_T \oplus u_{T' \cup T''} \ar[r]_{\varphi_{T,T' \cup T''}} & u_{T \cup T' \cup T''}
 }
$$
A morphism $f: \: \gekl{u_S, \varphi_{T,T'}} \to \gekl{\ti{u}_S, \ti{\varphi}_{T,T'}}$ consists of $f_S \in \fB(\ti{u}_S,u_S)$ for all $S$ such that $f_{\emptyset} = \id_0$ and the following diagram commutes for all disjoint $T, T'$:
$$
 \xymatrix{
 u_T \oplus u_{T'} \ar[d]_{f_T \oplus f_{T'}} \ar[r]^{\varphi_{T,T'}} & u_{T \cup T'} \ar[d]^{f_{T \cup T'}}\\
 \ti{u}_T \oplus \ti{u}_{T'} \ar[r]_{\ti{\varphi}_{T,T'}} & \ti{u}_{T \cup T'}
 }
$$
\edefin

The following result is for instance explained in \cite[Theorem~4.2]{EM}.
\btheo
$A \ma \fB(A)$ defines a functor from the category of finite based sets to the category of small categories, where a map of based sets $\alpha: \: A \to \bar{A}$ induces the functor $\gekl{u_S, \varphi_{T,T'}} \ma \gekl{u^{\alpha}_S, \varphi^{\alpha}_{T,T'}}$, with $u^{\alpha}_S \defeq u_{\alpha^{-1} S}$, $\varphi^{\alpha}_{T,T'} = \varphi_{\alpha^{-1} T,\alpha^{-1} T'}$, and a morphism $\gekl{f_S}$ is mapped to $\gekl{f^{\alpha}_S}$, with $f^{\alpha}_S = f_{\alpha^{-1} S}$.
\etheo

Next, we construct a $\Gamma$-space in the sense of \cite{Seg} (see also \cite{BF}), i.e., a functor from the category of finite based sets to the category of simplicial sets, sending the trivial based set to the simplicial set which is constantly given by one point. Given a finite based set $A$, let $\fN \fB (A)$ be the simplicial set with $p$-simplices $\fN_p \fB(A)$ consisting of elements of the form $(f_1, \dotsc, f_p)$, where $f_{\sqbullet}$ are morphisms in $\fB(A)$ such that $\mfd(f_{\sqbullet}) = \mft(f_{\sqbullet + 1})$. The face maps $\delta^p_{\sqbullet}: \: \fN_p \fB(A) \to \fN_{p-1} \fB(A)$ are given by 
$\delta^1_0(f) = \mfd(f)$, $\delta^1_1(f) = \mft(f)$ and 
$$
 \delta^p_{\sqbullet}(f_1, \dotsc, f_p)
 =
 \bfa
 (f_2, \dotsc, f_p) & \falls \sqbullet = 0,	\\
 (f_1, \dotsc, f_{\sqbullet} f_{\sqbullet + 1}, \dotsc, f_p) & \falls 0 < \sqbullet < p,\\
 (f_1, \dotsc, f_{p-1}) & \falls \sqbullet = p. 
 \efa
$$
Degeneracy maps are also part of the structure of a simplicial set, but since these are not needed for homology, we do not recall their definition here (see for instance \cite{GZ,GJ}).

Given a map $\alpha: \: A \to \bar{A}$ of finite based sets, define a map of simplicial sets $\fN \fB(\alpha): \: \fN \fB(A) \to \fN \fB(\bar{A})$ by setting $\fN_p \fB(\alpha)(f_1, \dotsc, f_p) \defeq (f^{\alpha}_1, \dotsc, f^{\alpha}_p)$.

The topological space $S^1$ is modelled by the simplicial set, also denoted by $S^1$, given by $S^1_q = \gekl{0, \dotsc, q}$ with base point $0$ and face maps 
$$
 d^q_{\bullet}: \: S^1_q \to S^1_{q-1}, \, a \ma
 \bfa
 a & \falls a \leq \bullet,\\
 a-1 & \falls a > \bullet.
 \efa
$$
For $\bullet = q$, $d^q_q$ sends $q$ to $0$. Again, we do not need the precise form of the degeneracy maps. Note that we have $S^1 \cong \Delta^1 / \partial \Delta^1$ (see for instance \cite{GZ,GJ}).

To obtain simplicial sets describing $S^n$, set $S^n = S^1 \wedge \dotso \wedge S^1$ ($n$ factors), as simplicial sets. Here the smash product $X \wedge S^1$ (where $X$ is some simplicial set) is given by 
$$
 (X \wedge S^1)_q = (X_q \times S^1_q) / ((\gekl{0_{X_q}} \times S^1_q) \cup (X_q \times \gekl{0_{S^1_q}})) \cong X_q\reg \times (S^1_q)\reg \cup \gekl{0},
$$
where $0$ stands for base point and, for a based set $A$ with base point $0$, $A\reg$ denotes $A \setminus \gekl{0}$. The face maps are induced from the face maps of $X$ and $S^1$. Concretely, 
$$
 S^n_q = \menge{(a_1, \dotsc, a_n)}{a_m \in \gekl{1, \dotsc, q}} \cup \gekl{0},
$$
and $d^q_{\bullet}: \: S^n_q \to S^n_{q-1}, \, (a_1, \dotsc, a_n) \ma (b_1, \dotsc, b_n)$ is given by
$$
 b_m = 
 \bfa
 a_m & \falls a_m \leq \bullet,\\
 a_m - 1 & \falls a_m > \bullet.
 \efa
$$
For $\bullet = q$, $b_m$ is defined to be $0$ if $a_m = q$. For $n = 0$, we set $S^0_q \defeq \gekl{0,1}$ for all $q$ and $d^q_{\bullet} = \id_{\gekl{0,1}}$.

Now consider the bisimplicial set $(p,q) \ma \fN_p \fB(S^n_q)$, with face maps $\delta_{\sqbullet}: \: \fN_p \fB(S^n_q) \to \fN_{p-1} \fB(S^n_q)$ and $\fN_q \fB(d_{\bullet}): \: \fN_p \fB(S^n_q) \to \fN_p \fB(S^n_{q-1})$. Furthermore, form the diagonal, i.e., the simplicial set $\fN \fB(S^n)$ given by $\fN \fB(S^n)_q \defeq \fN_q \fB(S^n_q)$ and the face maps, for $0 \leq \bullet \leq q$, given by the composites
$$
 \xymatrix{
 \fN_q \fB(S^n_q) \ar[r]^{\delta_{\bullet}} & \fN_{q-1} \fB(S^n_q) \ar[rr]^{\fN_{q-1} \fB(d_{\bullet})} & & \fN_{q-1} \fB(S^n_{q-1}).
 }
$$

Now we are ready for the definition of the algebraic K-theory spectrum $\Kz(\fB)$. The $n$-th simplicial set is given by $X_n = \fN \fB(S^n)$. We need to define the structure maps $\varsigma_n: \: \Sigma X_n \to X_{n+1}$. $\Sigma X_n$ is given by the smash product $\fN \fB(S^n) \wedge S^1$. As explained above, we have $(\fN \fB(S^n) \wedge S^1)_q \cong (\fN \fB(S^n)_q)\reg \times (S^1_q)\reg \cup \gekl{0}$, and the face maps are induced by the ones of $\fN \fB(S^n)$ and $S^1$. Every $b \in (S^1_q)\reg$ induces the map 
$$
 \iota_b: \: S^n_q \to S^{n+1}_q, \, \bma \ma (\bma,b), \, 0 \ma 0.
$$
This in turn induces $\fN_q \fB(\iota_b): \: \fN_q \fB(S^n_q) \to \fN_q \fB(S^{n+1}_q)$. Therefore, we obtain the based maps
$$
 (\fN \fB(S^n) \wedge S^1)_q \to \fN \fB(S^{n+1})_q, \, (\bmf,b) \ma \fN_q \fB(\iota_b)(\bmf), \, 0 \ma 0.
$$
In this way, we obtain the simplicial map $\varsigma_n: \: \fN \fB(S^n) \wedge S^1 \to \fN \fB(S^{n+1})$, as desired.

We have constructed a spectrum $\Kz(\fB)$ of simplicial sets. By taking geometric realizations, we also obtain a spectrum consisting of topological spaces.

It turns out that $\Kz(\fB)$ is a symmetric spectrum which is also a connective positive $\Omega$-spectrum, i.e., the adjoint maps $X_n \to \Omega X_{n+1}$ of $\varsigma_n$ are homotopy equivalences for all $n \geq 1$. The infinite loop space attached to $\Kz(\fB)$ is given by $\Omega^{\infty} \Kz(\fB) \defeq \Omega X_1$ (see \cite{Seg}). Note that $\Omega^{\infty} \Kz(\fB)$ coincides up to homotopy with $\Omega B \vert \fB \vert$. Here $\vert \fB \vert$ is the nerve or classifying space of $\fB$, and $B \vert \fB \vert$ is the bar construction of the monoid $\vert \fB \vert$, where the monoid structure is induced by the operation $\oplus$. We refer the reader to \cite{Ada78} for more information about infinite loop space theory.

Next, we briefly recall the definition of homology groups for simplicial sets and spectra. Let $\sfC$ be an abelian group as above. Let $X$ be a simplicial set with face maps $d_{\bullet}$. Define a chain complex $(C_* X, C_* d)$ by $C_q X \defeq \bigoplus_{X_q} \sfC$, setting $C_q d \defeq \sum_{\bullet = 0}^q (-1)^{\bullet} C_q d_{\bullet}$, where $C_q d_{\bullet}$ is the homomorphism $C_q X \to C_{q-1} X$ induced by $d_{\bullet}$. The homology $H_*(X,\sfC)$ is by definition the homology of the chain complex $(C_* X, C_* d)$. 

Note that by the Eilenberg-Zilber Theorem (see for instance \cite[Chapter~IV, \S~2.2]{GJ}), given a bisimplicial set like $(p,q) \ma \fN_p \fB(S^n_q)$, the homology of the diagonal (in our case $\fN \fB(S^n)$) is naturally isomorphic to the homology of the total complex (denoted by $C_{p,q} \fN_p \fB(S^n_q)$ in our case) associated to the bisimplicial set. Moreover, $H_*(X,\sfC)$ is canonically isomorphic to the singular homology with coefficients in $\sfC$ of the geometric realization of $X$ (see for instance \cite[Appendix Two, \S~1]{GZ}). 

Let us now define the homology of $\Kz(\fB)$. Applying the above definition of homology groups to $X_n$, we obtain the homology groups $H_{*+n}(X_n,\sfC)$.
\bdefin
$H_*(\Kz(\fB),\sfC) \defeq \ilim_n H_{*+n}(X_n,\sfC)$, where the inductive limit is taken with respect to the connecting maps
$$
 \xymatrix{
 H_{*+n}(X_n,\sfC) \ar[r]^{\hspace{-0.5cm} \cong} & H_{*+n+1}(\Sigma X_n,\sfC) \ar[rr]^{H_{*+n+1}(\varsigma_n)} & & H_{*+n+1}(X_{n+1},\sfC).
 }
$$
Here the first map is the suspension isomorphism.
\edefin

Let us also introduce the (stable) homotopy groups of $\Kz(\fB)$.
\bdefin
$\pi_*(\Kz(\fB)) \defeq \ilim_n \pi_{*+n}(X_n)$, where the inductive limit is taken with respect to the connecting maps
$$
 \xymatrix{
 \pi_{*+n}(X_n) \ar[r]^{\hspace{-0.5cm} \Sigma} & \pi_{*+n+1}(\Sigma X_n) \ar[r]^{(\varsigma_n)_*} & \pi_{*+n+1}(X_{n+1}).
 }
$$
Here the first map is the suspension homomorphism.
\edefin
Note that $\pi_*(\Kz(\fB)) \cong \pi_*(\Omega^{\infty} \Kz(\fB))$ (see for instance \cite[Chapter~I, \S~1]{Sch}).

The construction of algebraic K-theory spectra for small permutative spectra is functorial with respect to permutative functors, i.e., functors $\Phi: \: \fB \to \fC$ between small permutative categories which are compatible with all the structures, i.e., $\Phi(0) = 0$, $\Phi(u \oplus v) = \Phi(u) \oplus \Phi(v)$, similarly for morphisms, and $\Phi(\pi_{u,v}) = \pi_{\Phi(u),\Phi(v)}$. Indeed, given such a functor $\Phi: \: \fB \to \fC$ and a finite based set $A$, then we obtain a functor $\Phi(A): \: \fB(A) \to \fC(A)$ sending $\gekl{u_S, \varphi_{T,T'}}$ to $\gekl{\Phi(u_S), \Phi(\varphi_{T,T'})}$ and a morphism $\gekl{f_S}$ to $\gekl{\Phi(f_S)}$. Moreover, a based map $\alpha: \: A \to \bar{A}$ between finite based sets induces functors $\fB(\alpha): \: \fB(A) \to \fB(\bar{A})$ and $\fC(\alpha): \: \fC(A) \to \fC(\bar{A})$ such that the following diagram commutes:
$$
 \xymatrix{
 \fB(A) \ar[d]_{\fB(\alpha)} \ar[r]^{\Phi(A)} & \fC(A) \ar[d]^{\fC(\alpha)}\\
 \fB(\bar{A}) \ar[r]^{\Phi(\bar{A})} & \fC(\bar{A})
 }
$$
So $\Phi$ induces a map of bisimplicial sets between $(p,q) \ma \fN_p \fB(S^n_q)$ and $(p,q) \ma \fN_p \fC(S^n_q)$. It is straightforward to check that these maps are compatible with the structure maps defining the spectra $\Kz(\fB)$ and $\Kz(\fC)$.

\section{Small permutative categories of bisections}
\label{s:smallpermcat}

In this section, we construct small permutative categories from ample groupoids. Together with the construction of algebraic K-theory spectra from \S~\ref{ss:AlgKSmallPermCat}, this produces algebraic K-theory spectra for ample groupoids. We work in the more general setting of groupoid dynamical systems. This extra level of generality will be needed in \S~\ref{s:HK=HG}.

Let $G$ be an ample groupoid with locally compact Hausdorff unit space $\Gn$. Let $Z$ be a $G$-space with anchor map $\rho: \: Z \to \Gn$, i.e., the $G$-action $G \tensor[_{\rms}]{\times}{_{\rho}} Z \to Z, \, (g,z) \ma g.z$ is defined on $G \tensor[_{\rms}]{\times}{_{\rho}} Z = \menge{(g,z) \in G \times Z}{\rms(g) = \rho(z)}$. Assume that $\rho$ is a local homeomorphism. Let $\crO$ be a family of open Hausdorff subspaces of $Z$ covering $Z$, i.e., $Z = \bigcup_{O \in \crO} O$. If $Z$ is Hausdorff, then we can always take $\crO = \gekl{Z}$. Let $\cC \cO \defeq \menge{U \subseteq Z}{U \text{ compact open subspace of some } O \in \crO}$. Note that every $U \in \cC \cO$ is Hausdorff. Let $\cR$ be the full equivalence relation on $\Nz = \gekl{1, 2, 3, \dotsc}$, i.e., $\cR = \Nz \times \Nz$, $\cR^{(0)} = \Nz$ and we view an element $(j,i) \in \cR$ as a morphism from $i$ to $j$. Equip $\cR$ with the discrete topology. Equivalently, $\cR = \bigcup_N \cR_N$, where $\cR_N$ is the full equivalence relation on $\gekl{1, \dotsc, N}$. We set out to define a small permutative category $\fB_{G \curvearrowright Z}$ as follows: The objects of $\fBZ$ are given by 
$$
 \menge{\coprod_{i=1}^m (i,U_i) \subseteq \Nz \times Z}{U_i \in \cC \cO}.
$$
For $m=0$, the disjoint union becomes the empty set $\emptyset$. Now given objects $u = \coprod_{i=1}^m (i,U_i)$ and $v = \coprod_{j=1}^n (j,V_j)$, a morphism with target $v$ and domain $u$ is of the form $\coprod_{j,i} (\bbs_{j,i}, \sigma_{j,i}, U_{j,i})$, where $1 \leq i \leq m, 1 \leq j \leq n$, $\bbs_{j,i}$ is the map only defined on $\gekl{i}$ which sends $i$ to $j$, $\sigma_{j,i}$ are compact open bisections of $G$ and $U_{j,i} \in \cC \cO$ such that $\rms(\sigma_{j,i}) = \rho(U_{j,i})$, and $\rho \vert_{U_{j,i}}$ is a homeomorphism $U_{j,i} \cong \rho(U_{j,i})$. Moreover, we require that for all $i$, $U_i = \coprod_j U_{j,i}$, and that for all $j$, $V_j = \coprod_i \sigma_{j,i}.U_{j,i}$. Here we use the notation
$$
 \sigma.U \defeq \menge{g.z}{g \in \sigma, \, z \in U, \, \rms(g) = \rho(z)}
$$
for compact open bisections $\sigma \subseteq G$ and compact open subspaces $U \subseteq Z$. We denote the set of objects of $\fBZ$ by $\obj \fBZ$, the set of morphisms of $\fBZ$ by $\mor \fBZ$ and the set of morphisms with target $v$ and domain $u$ by $\fBZ(v,u)$. Given $\sigma \in \fBZ(v,u)$, we set $\mft(\sigma) \defeq v$ and $\mfd(\sigma) \defeq u$.

The composition of two morphisms $\tau = \coprod_{k,j} (\bbs_{k,j}, \tau_{k,j}, V_{k,j})$ and $\sigma = \coprod_{j,i} (\bbs_{j,i}, \sigma_{j,i}, U_{j,i})$ with $\mfd(\tau) = \mft(\sigma)$ is given by
$$
 \tau \sigma = \coprod_{k,i} (\bbs_{k,i}, \coprod_j \tau_{k,j} \sigma_{j,i}, \coprod_j \sigma^{-1}_{j,i}.(V_{k,j} \cap \sigma_{j,i}.U_{j,i})).
$$
Here the product of two bisections is by definition given by
$$
 \tau_{k,j} \sigma_{j,i} \defeq \menge{hg}{h \in \tau, \, g \in \sigma, \, \rms(h) = \rmr(g)}.
$$

We define the functor $\oplus$. To do so, we introduce the notation that for $m \in \Nz$, $\bbt_m$ denotes the map $\Nz \to \Nz$ given by addition with $m$. Given $u = \coprod_{i=1}^m (i,U_i)$, $u' = \coprod_{i'=1}^{m'} (i',U'_{i'})$, define
$$
 u \oplus u' \defeq \coprod_{i=1}^m (i,U_i) \amalg \coprod_{i'=1}^{m'} (\bbt_m(i'),U'_{i'}).
$$
Furthermore, define $\pi_{u,u'} \in \fBZ(u' \oplus u, u \oplus u')$ by setting
$$
 \pi_{u,u'} \defeq \coprod_{i=1}^m (\bbs_{m'+i,i},\rho(U_i),U_i) \amalg \coprod_{i'=1}^{m'} (\bbs_{i',m+i'},\rho(U'_{i'}),U'_{i'}).
$$
Moreover, given $u = \coprod_{i=1}^m (i,U_i)$, $u' = \coprod_{i'=1}^{m'} (i',U'_{i'})$, $v = \coprod_{j=1}^n (j,V_j)$, $v' = \coprod_{j'=1}^{n'} (j',V'_{j'})$ in $\obj \fBZ$, $\sigma = \coprod_{j,i} (\bbs_{j,i}, \sigma_{j,i}, U_{j,i}) \in \fBZ(v,u)$ and $\sigma' = \coprod_{j',i'} (\bbs_{j',i'}, \sigma_{j',i'}, U_{j',i'}) \in \fBZ(v',u')$, define
$$
 \sigma \oplus \sigma' \defeq \coprod_{j,i} (\bbs_{j,i}, \sigma_{j,i}, U_{j,i}) \amalg \coprod_{j',i'} (\bbs_{\bbt_n(j'),\bbt_m(i')}, \sigma_{j',i'}, U_{j',i'}).
$$
It is now straightforward to check that $\oplus$ indeed defines a functor and that $\gekl{\pi_{u,u'}}$ are natural isomorphisms such that $\fBZ$ becomes a small permutative category, with unit $\emptyset$.

Every element $U \in \cC \cO$ will be viewed as an element of $\obj \fBZ$ by identifying $U$ with $(1,U) \in \obj \fBZ$. With this convention, it is clear that $\coprod_{i=1}^m (i,U_i) = U_1 \oplus \dotso \oplus U_m$.

\bremark
\label{rem:fBZ=Groupoid}
All morphisms in $\fBZ$ are actually invertible, i.e., $\fBZ$ is a groupoid.
\eremark

We will apply this construction to the following special cases: $Z = \Gnu$ as in \S~\ref{sss:C,GPDHom}, viewed as a $G$-space via the anchor map $\rho: \: \Gnu \to \Gn, \, (g_0, \dotsc, g_{\nu - 1}) \ma \rmr(g_0)$ for all $\nu \geq 0$ and action $g.(g_0, \dotsc, g_{\nu - 1}) \defeq (g g_0, \dotsc, g_{\nu - 1})$ if $\nu \geq 1$ and $g.x \defeq \rmr(g)$ if $\nu = 0$. Let $\crO^{(\nu)}$ be as in \S~\ref{sss:C,GPDHom}, i.e., $\crO^{(\nu)}$ denotes the collection of subsets of $\Gnu$ of the form $O_0 \tensor[_{\rms}]{\times}{_{\rmr}} \dotso \tensor[_{\rms}]{\times}{_{\rmr}} O_{\nu-1}$, where $O_{\mu}$ are open bisections with $\rms(O_{\mu + 1}) = \rmr(O_{\mu})$. For $\nu = 0$, we consider $Z = \Gn$ and $\crO^{(0)} = \gekl{\Gn}$. In that case, we write $\fB_G \defeq \fB_{G \curvearrowright \Gn}$. We will also restrict the $G$-action to the (trivial) $\Gn$-action $\Gn \curvearrowright \Gnum$ and consider $\fB_{\Gn \curvearrowright \Gnum}$.

\bremark
\label{rem:B=BisecInRxG}
For $\nu = 0$, objects in $\fB_G = \fB_{G \curvearrowright \Gn}$ are just compact open subspaces $u \subseteq \Nz \times \Gn$ and morphisms in $\fB_G$ are nothing else but compact open bisections of $\cR \times G$, the direct product of the groupoids $\cR$ and $G$. In this case, we may and will reduce the general form $\sigma = \coprod_{j,i} (\bbs_{j,i}, \sigma_{j,i}, U_{j,i})$ of a morphism $\sigma$ to $\sigma = \coprod_{j,i} (\bbs_{j,i}, \sigma_{j,i})$ because the component $U_{j,i}$ is superfluous since $U_{j,i} = \rms(\sigma_{j,i})$. 
\eremark

\bremark
\label{rem:Buu=Fu}
We have $\fB_G(u,u) = \bmF((\cR \times G)_u^u)$. Moreover, given $u = \coprod_i (i,U_i) \in \obj \fB_G$ with pairwise disjoint subspaces $U_i \subseteq \Gn$, write $U \defeq \coprod_i U_i \subseteq \Gn$ and set $\omega \defeq \coprod_i (\bbs_{1,i},U_i) \in \fB_G(U,u)$. Then 
$$
 \fB_G(u,u) = \bmF((R \times G)_u^u) \cong \bmF(G_U^U), \, \sigma \ma \omega \sigma \omega^{-1}.
$$
\eremark

\bremark
Our category of bisections $\fB_G$ coincides with the category in \cite[Definition~2.1]{Li21b}.
\eremark

\subsection{Functoriality of our construction}
\label{ss:FunctPermCat}

Our construction of small permutative categories of bisections is functorial for two types of maps, open embeddings and fibrewise bijective proper surjections. These types of maps also appear in \cite[\S~5]{Li18}. In the following, let $G$ and $\ti{G}$ be ample groupoids with locally compact Hausdorff unit spaces. Denote by $\rms,\rmr$ the source and range maps of $G$ and by $\ti{\rms},\ti{\rmr}$ the source and range maps of $\ti{G}$. Let $Z$ be a $G$-space with anchor map $\rho: \: Z \to \Gn$ and $\ti{Z}$ be a $\ti{G}$-space with anchor map $\ti{\rho}: \: \ti{Z} \to \ti{G}^{(0)}$. Assume that $\rho$ and $\ti{\rho}$ are local homeomorphisms. Let $\crO$ be a family of open Hausdorff subspaces of $Z$ covering $Z$ and $\ti{\crO}$ a family of open Hausdorff subspaces of $\ti{Z}$ covering $\ti{Z}$. As above, construct the small permutative categories $\fB_{G \curvearrowright Z}$ and $\fB_{\ti{G} \curvearrowright \ti{Z}}$.

\subsubsection{The case of open embeddings}
\label{sss:OpenEmb}

Suppose that $\phi: \: G \to \ti{G}$ is a groupoid homomorphism which is an embedding with open image, and that $\psi: \: Z \to \ti{Z}$ is a continuous map such that for all $O \in \crO$ there exists $\ti{O} \in \ti{\crO}$ such that $\psi(O) \subseteq \ti{O}$ and that $\psi$ restricts to a homeomorphism $\psi \vert_O: \: O \cong \psi(O)$. Furthermore, we require that the diagram
$$
 \xymatrix{
 Z \ar[d]_{\rho} \ar[r]^{\psi} & \ti{Z} \ar[d]^{\ti{\rho}}\\
 \Gn \ar[r]_{\phi} & \ti{G}^{(0)}
 }
$$
commutes, and that $\psi(g.z) = \phi(g).\psi(z)$ for all $g \in G$ and $z \in Z$ with $\rms(g) = \rho(z)$. In this situation, we define a functor $F_{\phi,\psi}: \: \fBZ \to \fB_{\ti{G} \curvearrowright \ti{Z}}$ as follows: On objects, set $F_{\phi,\psi}(\coprod_i (i,U_i)) \defeq \coprod_i (i,\psi(U_i))$. Given a morphism $\sigma = \coprod_{j,i} (\bbs_{j,i}, \sigma_{j,i}, U_{j,i})$ in $\fBZ$, set $F_{\phi,\psi}(\sigma) \defeq \coprod_{j,i} (\bbs_{j,i}, \phi(\sigma_{j,i}), \psi(U_{j,i}))$. It is straightforward to check that this is well-defined, i.e., this defines a morphism in $\fB_{\ti{G} \curvearrowright \ti{Z}}$. For instance, we have that 
$$
 \ti{\rms}(\phi(\sigma_{j,i})) = \phi(\rms(\sigma_{j,i})) = \phi(\rho(U_{j,i})) = \ti{\rho}(\psi(U_{j,i})),
$$
and moreover, the restriction
$$
 \ti{\rho} \vert_{\psi(U_{j,i})}: \: \psi(U_{j,i}) \to \ti{\rho}(\psi(U_{j,i})) = \phi(\rho(U_{j,i}))
$$
is a homeomorphism because the composite
$$
 \xymatrix{
 U_{j,i} \ar[r]^{\hspace{-0.25cm} \psi} & \psi(U_{j,i}) \ar[r]^{\hspace{-0.25cm} \ti{\rho}} & \ti{\rho}(\psi(U_{j,i}))
 }
$$
coincides with
$$
 \xymatrix{
 U_{j,i} \ar[r]^{\hspace{-0.25cm} \rho} & \rho(U_{j,i}) \ar[r]^{\hspace{-0.25cm} \phi} & \phi(\rho(U_{j,i})),
 }
$$
and $\psi$ as well as $\phi \circ \rho$ are homeomorphisms.

Given two morphisms $\tau = \coprod_{k,j} (\bbs_{k,j}, \tau_{k,j}, V_{k,j})$ and $\sigma = \coprod_{j,i} (\bbs_{j,i}, \sigma_{j,i}, U_{j,i})$ with $\mfd(\tau) = \mft(\sigma)$, we have
\begin{align*}
 F_{\phi,\psi}(\tau) F_{\phi,\psi}(\sigma) 
 &= \coprod_{k,i} \big( \bbs_{k,i}, \coprod_j \phi(\tau_{k,j}) \phi(\sigma_{j,i}), \coprod_j \phi(\sigma_{j,i})^{-1}.(\psi(V_{k,j}) \cap \phi(\sigma_{j,i}).\psi(U_{j,i})) \big)\\
 &= \coprod_{k,i} \big( \bbs_{k,i}, \coprod_j \phi(\tau_{k,j} \sigma_{j,i}), \coprod_j \psi(\sigma_{j,i}^{-1}.(V_{k,j} \cap \sigma_{j,i}.U_{j,i})) \big)\\
 &= F_{\phi,\psi}(\tau \sigma).
\end{align*}
Hence $F_{\phi,\psi}$ respects composition. Furthermore, $F_{\phi,\psi}$ also respects $\oplus$ by construction. This shows that $F_{\phi,\psi}$ is a permutative functor. Moreover, our construction is (covariantly) functorial in $(\phi,\psi)$, in the sense that $F_{\phi',\psi'} F_{\phi,\psi} = F_{\phi' \phi,\psi' \psi}$.

\subsubsection{The case of fibrewise bijective proper surjections}
\label{sss:FiBiPrSu}

Now suppose that $\phi: \: \ti{G} \to G$ is a groupoid homomorphism which is a fibrewise bijective proper surjection. \an{Fibrewise bijective} means that for all $\ti{x} \in \ti{G}^{(0)}$, $\phi$ restricts to a bijection $\ti{\rmr}^{-1}(\ti{x}) \to \rmr^{-1}(\phi(\ti{x}))$. (Equivalently, we could consider fibres of the source maps.) Moreover, suppose that $\psi: \: \ti{Z} \to Z$ is a continuous map such that for all $O \in \crO$ there exists $\ti{O} \in \ti{\crO}$ such that $\psi^{-1}(O) \subseteq \ti{O}$ and that $\psi$ restricts to a proper map $\psi \vert_{\psi^{-1}(O)}: \: \psi^{-1}(O) \to O$. Furthermore, we require that the diagram
$$
 \xymatrix{
 Z \ar[d]_{\rho} \ar[r]^{\psi} & \ti{Z} \ar[d]^{\ti{\rho}}\\
 \Gn \ar[r]_{\phi} & \ti{G}^{(0)}
 }
$$
commutes, and that $\psi(g.z) = \phi(g).\psi(z)$ for all $g \in G$ and $z \in Z$ with $\rms(g) = \rho(z)$. In addition, assume that for all $\ti{x} \in \ti{G}^{(0)}$, $\psi$ restricts to a bijection
$$
 \ti{\rho}^{-1}(\ti{x}) \cap \psi^{-1}(O) \to \rho^{-1}(\phi(\ti{x})) \cap O.
$$
In this situation, we define a functor $F^{\phi,\psi}: \: \fBZ \to \fB_{\ti{G} \curvearrowright \ti{Z}}$ as follows: On objects, set 
$$
 F^{\phi,\psi}(\coprod_i (i,U_i)) \defeq \coprod_i (i,\psi^{-1}(U_i)).
$$
Given a morphism $\sigma = \coprod_{j,i} (\bbs_{j,i}, \sigma_{j,i}, U_{j,i})$ in $\fBZ$, set 
$$
 F^{\phi,\psi}(\sigma) \defeq \coprod_{j,i} (\bbs_{j,i}, \phi^{-1}(\sigma_{j,i}), \psi^{-1}(U_{j,i})).
$$
Let us now check that this is well-defined, i.e., this defines a morphism in $\fB_{\ti{G} \curvearrowright \ti{Z}}$. 

First observe that for all compact open subspaces $U$ contained in some $O \in \crO$, we have $\phi^{-1}(\rho(U)) = \ti{\rho}(\psi^{-1}(U))$. Indeed, \an{$\supseteq$} is clear. Given $x \in \rho(U)$, let $\ti{x} \in \ti{G}^{(0)}$ with $\phi(\ti{x}) = x$ be arbitrary and choose $z \in \rho^{-1}(x) \cap O$. As $\psi$ restricts to a bijection $\ti{\rho}^{-1}(\ti{x}) \cap \psi^{-1}(O) \to \rho^{-1}(\phi(\ti{x})) \cap O$, there exists $\ti{z} \in \ti{\rho}^{-1}(\ti{x}) \cap \psi^{-1}(O)$ with $\psi(\ti{z}) = z$. It follows that $\ti{z} \in \psi^{-1}(U)$. Hence $\ti{x} = \ti{\rho}(\ti{z}) \in \ti{\rho}(\psi^{-1}(U))$. 

Secondly, observe that for every compact open bisection $\sigma \subseteq G$, $\ti{\rms}(\phi^{-1}(\sigma)) = \phi^{-1}(\rms(\sigma))$. Indeed, \an{$\subseteq$} is clear. Given $x \in \rms(\sigma)$, let $\ti{x} \in \ti{G}^{(0)}$ be arbitrary with $\phi(\ti{x}) = x$. Choose $g \in \rms^{-1}(x) \cap \sigma$. As $\phi$ restricts to a bijection $\ti{\rms}^{-1}(\ti{x}) \to \rms^{-1}(x)$, there exists $\ti{g} \in \ti{\rms}^{-1}(\ti{x})$ with $\phi(\ti{g}) = g$. Hence $\ti{g} \in \phi^{-1}(\sigma)$. Thus $\ti{x} = \ti{\rms}(\ti{g}) \in \ti{\rms}(\phi^{-1}(\sigma))$. This shows \an{$\supseteq$}.

It follows that $\ti{\rms}(\phi^{-1}(\sigma_{j,i})) = \phi^{-1}(\rms(\sigma_{j,i})) = \phi^{-1}(\rho(U_{j,i})) = \ti{\rho}(\psi^{-1}(U_{j,i}))$.

Moreover, the restriction $\ti{\rho} \vert_{\psi^{-1}(U_{j,i})}: \: \psi^{-1}(U_{j,i}) \to \ti{\rho}(\psi^{-1}(U_{j,i})) = \phi^{-1}(\rho(U_{j,i}))$ is a homeomorphism. It suffices to show that this restriction is bijective. Given $\ti{z}_1, \ti{z}_2 \in \psi^{-1}(U_{j,i})$ with $\ti{\rho}(\ti{z}_1) = \ti{\rho}(\ti{z}_2) = \ti{x}$, we have $\rho(\psi(\ti{z}_1)) = \phi(\ti{\rho}(\ti{z}_1)) = \phi(\ti{x}) = \phi(\ti{\rho}(\ti{z}_2)) = \rho(\psi(\ti{z}_2))$. As $\rho$ is injective on $U_{j,i}$, we deduce that $\psi(\ti{z}_1) = \psi(\ti{z}_2)$. But $\psi$ restricts to a bijection $\ti{\rho}^{-1}(\ti{x}) \cap \psi^{-1}(O) \to \rho^{-1}(\phi(\ti{x})) \cap O$, where $O \in \crO$ is such that $U_{j,i} \subseteq O$. Hence we conclude that $\ti{z}_1 = \ti{z}_2$, as desired.

Given two morphisms $\tau = \coprod_{k,j} (\bbs_{k,j}, \tau_{k,j}, V_{k,j})$ and $\sigma = \coprod_{j,i} (\bbs_{j,i}, \sigma_{j,i}, U_{j,i})$ with $\mfd(\tau) = \mft(\sigma)$, we have
\begin{align*}
 F^{\phi,\psi}(\tau) F^{\phi,\psi}(\sigma) 
 &= \coprod_{k,i} \big( \bbs_{k,i}, \coprod_j \phi^{-1}(\tau_{k,j}) \phi^{-1}(\sigma_{j,i}), \coprod_j \phi^{-1}(\sigma_{j,i})^{-1}.(\psi^{-1}(V_{k,j}) \cap \phi^{-1}(\sigma_{j,i}).\psi(U_{j,i})) \big)\\
 &= \coprod_{k,i} \big( \bbs_{k,i}, \coprod_j \phi^{-1}(\tau_{k,j} \sigma_{j,i}), \coprod_j \psi^{-1}(\sigma_{j,i}^{-1}.(V_{k,j} \cap \sigma_{j,i}.U_{j,i})) \big)\\
 &= F^{\phi,\psi}(\tau \sigma).
\end{align*}
Hence $F^{\phi,\psi}$ respects composition. Furthermore, $F^{\phi,\psi}$ also respects $\oplus$ by construction. This shows that $F^{\phi,\psi}$ is a permutative functor. Moreover, our construction is (contravariantly) functorial in $(\phi,\psi)$, in the sense that $F^{\phi',\psi'} F^{\phi,\psi} = F^{\phi \phi',\psi \psi'}$.

\section{Homology for algebraic K-theory spectra of bisections in terms of groupoid homology}
\label{s:HK=HG}

Let $G$ be an ample groupoid with locally compact Hausdorff unit space $\Gn$. Fix an abelian group $\sfC$. Our goal is to identify (reduced) stable homology of $\Kz(\fB_G)$ with groupoid homology of $G$, i.e., $\ti{H}_*(\Kz(\fB_G),\sfC) \cong H_*(G,\sfC)$. Here $\ti{H}_*(\Kz(\fB_G),\sfC) \defeq H_*(\Kz(\fB_G),\sfC)$ for $*>0$ and $H_0(\Kz(\fB_G),\sfC) = \ti{H}_0(\Kz(\fB_G),\sfC) \oplus \sfC$, where the second direct sum comes from the base point corresponding to the unit $\emptyset$ of the small permutative category $\fB_G$.

\subsection{Functors inducing homotopy equivalences of classifying spaces}

First, we establish a criterion for certain functors to induce homotopy equivalences of classifying spaces. Let us start by describing the setting. Let $\Phi: \: \fC \to \fG$ be a functor between small categories. Assume that $\fG$ is a groupoid, and that $\Phi$ is faithful, i.e., the induced maps $\fC(*,\bullet) \to \fG(\Phi(*),\Phi(\bullet))$ are injective for all $\bullet, * \in \obj \fC$. Given an object $u$ of $\fG$, we first recall the definition of the category $u \backslash \Phi$ from \cite[\S~1]{Qui}. Objects of $u \backslash \Phi$ consist of pairs $(v,\sigma)$, where $v \in \obj \fC$ and $\sigma \in \fG(\Phi(v),u)$. A morphism in $u \backslash \Phi$ from $(v,\sigma)$ to $(w,\tau)$ is given by $f \in \fC(w,v)$ such that $\Phi(f) \sigma = \tau$, i.e., the diagram
$$
 \xymatrix{
 u \ar[dr]_{\tau} \ar[r]^{\sigma} & \Phi(v) \ar[d]^{\Phi(f)}\\
  & \Phi(w)
 }
$$
commutes in $\fG$.

We now define a new category $\fF_{u,\Phi}$ as a quotient of $u \backslash \Phi$. Objects of $\fF_{u,\Phi}$ are equivalence classes $[v,\sigma]$ of objects $(v,\sigma)$ of $u \backslash \Phi$, where $(v,\sigma)$ and $(w,\tau)$ are equivalent if there exists an invertible morphism of $u \backslash \Phi$ from $(v,\sigma)$ to $(w,\tau)$, i.e., there exists an invertible element $a \in \fC(w,v)$ such that $\Phi(a) \sigma = \tau$, i.e., the diagram
$$
 \xymatrix{
 u \ar[dr]_{\tau} \ar[r]^{\sigma} & \Phi(v) \ar[d]^{\Phi(a)}\\
  & \Phi(w)
 }
$$
commutes in $\fG$. Morphisms of $\fF_{u,\Phi}$ are equivalence classes $[f]$ of morphisms $f$ of $u \backslash \Phi$, where $f: \: (v,\sigma) \to (w,\tau)$ and $f': \: (v',\sigma') \to (w',\tau')$ are equivalent if there exist invertible elements $b \in \fC(w',w)$, $b' \in \fC(v',v)$, which are invertible morphisms in $u \backslash \Phi$, such that $b f = f' b'$ in $\fC$ (here we view $f$ and $f'$ as morphisms in $\fC$), i.e., the diagram
$$
 \xymatrix{
 v \ar[d]_{b'} \ar[r]^{f} & w \ar[d]^{b}\\
 v' \ar[r]_{f'} & w'
 }
$$
commutes in $\fC$.

Note that a morphism $f$ in $u \backslash \Phi$ from $(v,\sigma)$ to $(w,\tau)$, if it exists, is unique. This is because $\fG$ is a groupoid and $\Phi$ is faithful.

\bprop
\label{prop:uPhi=FuPhi}
The classifying spaces of $u \backslash \Phi$ and $\fF_{u,\Phi}$ are homotopy equivalent.
\eprop
\bproof
For every $c \in \obj \fF_{u,\Phi}$, choose $x_c \in \obj (u \backslash \Phi)$ such that $[x_c] = c$. Given $[f] \in \mor \fF_{u,\Phi}$, let $\mathring{f} \in \mor (u \backslash \Phi)$ be the unique morphism from $x_{\mfd(f)}$ to $x_{\mft(f)}$. This defines a functor $\fF_{u,\Phi} \to u \backslash \Phi$. We also have the canonical functor $u \backslash \Phi \to \fF_{u,\Phi}$ given by forming equivalence classes. By construction, the composite $\fF_{u,\Phi} \to u \backslash \Phi \to \fF_{u,\Phi}$ is the identity on $\fF_{u,\Phi}$. Now let $\Theta$ be the composite $u \backslash \Phi \to \fF_{u,\Phi} \to u \backslash \Phi$. By construction, $\Theta(o) = x_{[o]}$ on objects and $\Theta(f) = \mathring{f}$ on morphisms. We claim that there is a natural transformation $T: \: \id_{u \backslash \Phi} \Rarr \Theta$. Indeed, given an object $o \in \obj (u \backslash \Phi)$, let $T_o \in (u \backslash \Phi)(x_{[o]},o)$ be the unique morphism in $u \backslash \Phi$ from $o$ to $x_{[o]}$. Given a morphism $f \in (u \backslash \Phi)(\ti{o},o)$, the diagram
$$
 \xymatrix{
 \ti{o} \ar[d]_{T_{\ti{o}}} & \ar[l]^{f} \ar[d]^{T_o} o\\
 x_{[\ti{o}]} & \ar[l]_{\mathring{f}} x_{[o]}
 }
$$
commutes because of uniqueness of morphisms in $u \backslash \Phi$. Now our proof is complete because of \cite[Proposition~2]{Qui}.
\eproof

The following is now an immediate consequence of Proposition~\ref{prop:uPhi=FuPhi} and \cite[Theorem~A]{Qui}.
\bcor
\label{cor:Phihomoequi}
If the classifying space of $\fF_{u,\Phi}$ is contractible for all $u \in \obj \fG$, then $\Phi$ induces a homotopy equivalence of classifying spaces.
\ecor

\subsection{Homology for certain free and proper actions}

Now let $\nu \geq 1$ and consider $\fB_{\Gn \acts \Gnum}$ and $\fB_{G \acts \Gnu}$ as defined in \S~\ref{s:smallpermcat}. We set out to define a functor $I: \: \fB_{\Gn \acts \Gnum} \to \fB_{G \acts \Gnu}$. Given $U \in \cC \cO^{(\nu - 1)}$, let 
$$
 I(U) \defeq \menge{(\rho(z),z) \in \Gnu}{z \in U}.
$$
Now define 
$$
 I(\coprod_i (i,U_i)) \defeq \coprod_i (i,I(U_i)).
$$ 
On morphisms, define 
$$
 I(\coprod_{j,i} (\bbs_{j,i},\sigma_{j,i},U_{j,i})) \defeq \coprod_{j,i} (\bbs_{j,i},\sigma_{j,i},I(U_{j,i})).
$$
Note that $\sigma_{j,i} = \rho(U_{j,i})$ because $\coprod_{j,i} (\bbs_{j,i},\sigma_{j,i},U_{j,i})$ is a morphism in $\fB_{\Gn \acts \Gnum}$.

It is straightforward to check that $I$ is a permutative functor, and that $I$ is faithful. Therefore, we are in the setting of Proposition~\ref{prop:uPhi=FuPhi} and Corollary~\ref{cor:Phihomoequi}.
\bprop
\label{prop:FuITrivial}
For all $u \in \obj \fB_{G \acts \Gnu}$, the category $\fF_{u,I}$ is trivial.
\eprop
\bproof
Since $\fB_{\Gn \acts \Gnum}$ is a groupoid (see Remark~\ref{rem:fBZ=Groupoid}), our claim follows from the following observations.
\begin{enumerate}
\item[(i)] For all $u \in \obj \fB_{G \acts \Gnu}$, there exists $v \in \obj \fB_{\Gn \acts \Gnum}$ and a morphism $\sigma \in \fB_{G \acts \Gnu}(I(v),u)$ from $u$ to $I(v)$.
\item[(ii)] Given objects $u \in \obj \fB_{G \acts \Gnu}$, $v, v' \in \obj \fB_{\Gn \acts \Gnum}$ and morphisms $\sigma \in \fB_{G \acts \Gnu}(I(v),u)$, $\sigma' \in \fB_{G \acts \Gnu}(I(v'),u)$, there exists a morphism $\tau \in \fB_{\Gn \acts \Gnum}(v',v)$ from $v$ to $v'$ such that $I(\tau) \sigma = \sigma'$.
\end{enumerate}
To prove (i), write $u = \coprod_i (i,U_i)$. For each $i$, write $U_i = \coprod_{a_i} U_{i,a_i}$, where $U_{i,a_i} = U_{i,a_i}^0 \tensor[_{\rms}]{\times}{_{\rmr}} \dotso \tensor[_{\rms}]{\times}{_{\rmr}} U_{i,a_i}^{\nu - 1}$ and $U_{i,a_i}^{\mu}$ are compact open bisections with $\rms(U_{i,a_i}^{\mu}) = \rmr(U_{i,a_i}^{\mu + 1})$. Define for all $i$ and $a_i$ the compact open subspace $V_{i,a_i} \defeq U_{i,a_i}^1 \tensor[_{\rms}]{\times}{_{\rmr}} \dotso \tensor[_{\rms}]{\times}{_{\rmr}} U_{i,a_i}^{\nu - 1}$. Furthermore, define $v \defeq \coprod_{i, a_i} ((i,a_i),V_{i,a_i})$ and $\sigma \defeq \coprod_{i,a_i} (\bbs_{(i,a_i),i}, (U_{i,a_i}^0)^{-1}, U_{i,a_i})$. Then it is straightforward to check that $\sigma \in \fB_{G \acts \Gnu}(I(v),u)$.

Let us now prove (ii). Write $I(v) = \coprod_j (j,V_j)$, $I(v') = \coprod_k (k,V'_k)$. Let us show that $\sigma' \sigma^{-1} \in \fB_{G \acts \Gnu}(I(v'),I(v))$ is the image of a morphism in $\fB_{\Gn \acts \Gnum}$ under $I$. Write $\sigma' \sigma^{-1} = \coprod_{k,j} (\bbs_{k,j}, \tau_{k,j}, V_{k,j})$. For all $(g_0, \dotsc, g_{\nu - 1}) \in V_{k,j}$, we have $g_0 \in \Gn$. Similarly, for every $(h_0, \dotsc, h_{\nu - 1}) \in V'_k$, we have $h_0 \in \Gn$. It follows that $\tau_{k,j} \subseteq \Gn$. Hence, indeed, $\sigma' \sigma^{-1} = I(\tau)$ for some $\tau \in \fB_{\Gn \acts \Gnum}(v',v)$, as desired.
\eproof

The following is an immediate consequence of Corollary~\ref{cor:Phihomoequi} and Proposition~\ref{prop:FuITrivial}.
\bcor
\label{cor:I=homoequi}
$I$ induces a homotopy equivalence of classifying spaces.
\ecor

Now let $C$ be a compact Hausdorff subspace of $\Gnum$. Let $\fB_C$ be the full subcategory of $\fB_{\Gn \acts \Gnum}$ whose objects are of the form $\coprod_i (i,U_i)$, where $U_i \subseteq C$. It follows that morphisms of $\fB_C$ are of the form $\coprod_{j,i} (\bbs_{j,i},\sigma_{j,i},U_{j,i})$, with $U_{j,i} \subseteq C$ (and hence $\sigma_{j,i}.U_{j,i} = U_{j,i} \subseteq C$). As $C$ is totally disconnected, we can describe $C$ as $C \cong \plim_{l \in \fL} \cU_l$, where $\cU_l$ consist of finitely many compact open subspaces which partition $C$, and $\fL$ is the index set given by all these partitions, partially ordered by refinement. For $l \in \fL$, let $\fB_l$ be the small permutative category with objects of the form $\coprod_i (i,U_i)$ with $U_i \in \cU_l$, and morphisms of the form $\coprod_{j,i} (\bbs_{j(i),i}, \rho(U_i), U_i)$ from $\coprod_i (i,U_i)$ to $\coprod_j (j,U_j)$, where $i \ma j(i)$ is a bijection. In other words, the only morphisms in $\fB_l$ are given by permutations. Note that $\fB_l$ is a subcategory of $\fB_C$.

\blemma
\label{lem:H(K(Bl))}
For all $l \in \fL$, we have a canonical isomorphism
$$
 H_*(\Kz(\fB_l),\sfC)
 \cong
 \bfa
 (\bigoplus_{\cU_l} \sfC) \oplus \sfC & \falls * = 0,\\
 \gekl{0} & \sonst.
 \efa
$$
\elemma
\setlength{\parindent}{0cm} \setlength{\parskip}{0cm}

\bproof
It is easy to see that $\fB_l$ coincides with the free permutative category $P \cU_l$ in the sense of \cite[\S~1]{Tho}, where $\cU_l$ is the category with object set $\cU_l$ and only identity morphisms. Moreover, it is observed in \cite[\S~1]{Tho} that $P \cU_l$ is equivalent to the free symmetric monoidal category $S \cU_l$ on $\cU_l$, so that \cite[Lemma~2.3 and Lemma~2.5]{Tho} imply that $\Kz(P \cU_l)$ is weakly homotopy equivalent to $\Sigma^{\infty}(B \cU_l)^+$. Hence we obtain 
$$
 H_*(\Kz(\fB_l),\sfC) \cong H_*(\Kz(P \cU_l),\sfC) \cong H_*(\Sigma^{\infty}(B \cU_l)^+,\sfC),
$$
and now our claim follows from $H_*(B \cU_l, \sfC) \cong H_*(\cU_l, \sfC)$. Here we view $\cU_l$ as a discrete space. 
\eproof
\setlength{\parindent}{0cm} \setlength{\parskip}{0.5cm}

For $k \leq l$, define a permutative functor $I_{l,k}: \: \fB_k \to \fB_l$ by $U \ma \bigoplus V(U)$ for $U \in \cU_k$, where the sum is taken over all $V(U) \in \cU_l$ contained in $U$, so that $U = \coprod V(U)$ as subspaces in $C$, and by extending this to all objects via $\coprod_i (i,U_i) \ma \bigoplus_i (\bigoplus V(U_i))$. Here we are working with a fixed ordering of $\gekl{V(U)}$ for every $U$. On morphisms, let $I_{l,k}(\coprod_{j,i} (\bbs_{j(i),i}, \rho(U_i), U_i))$ be the morphism
$$
 \bigoplus_i (\bigoplus V(U_i)) \to \bigoplus_j (\bigoplus V(U_j))
$$
induced by the permutation $i \ma j(i)$. Let us now form the homotopy colimit, in the sense of \cite[Construction~3.22]{Tho}, of the following functor $F$ from $\fL$ to permutative categories: We define $F(l) \defeq \fB_l$ on objects $l \in \fL$, and a morphism $k \to l$ (i.e., $k, l \in \fL$ satisfying $k \leq l$) is mapped to $I_{l,k}$ under $F$. We recall the construction of $\fH \defeq {\rm hocolim} \, F$ from \cite[Construction~3.22]{Tho}. Objects of $\fH$ are of the form $(l_1,u_1) \oplus \dotso \oplus (l_n,u_n)$, where $l_i \in \fL$ and $u_i \in F(l_i)$. Morphisms from $(l_1,u_1) \oplus \dotso \oplus (l_n,u_n)$ to $(l'_1,u'_1) \oplus \dotso \oplus (l'_m,u'_m)$ are given by $(\lambda_i, \psi, \chi_j)$, where $\psi: \: \gekl{1, \dotsc, n} \to \gekl{1, \dotsc, m}$ is a surjective map, $\lambda_i: \: l_i \to l'_{\psi(i)}$ are morphisms in $\fL$ (in our case this simply means $l_i \leq l'_{\psi(i)}$), and $\chi_j: \: \bigoplus_{\psi(i) = j} F(\lambda_i)(u_i) \to u'_j$ are morphisms in $F(l'_j)$.

Now assume that $C$ is contained in some $O \in \crO^{(\nu - 1)}$, and define a permutative functor $H: \: \fH \to \fB_C$ as follows: Given $u = \coprod_{\bullet} (\bullet,U_{\bullet})$ with $U_{\bullet} \in \cU_l$, set $H(l,u) \defeq \coprod_{\bullet} (\bullet,U_{\bullet})$ viewed as an object in $\fB_C$. Define $H$ on a morphism given by data $(\lambda_i, \psi, \chi_j)$ as above by sending it to the morphism given by the composition
$$
 \bigoplus_i H(l_i,u_i) \to \bigoplus_i H(l'_{\psi(i)},F(\lambda_i)(u_i)) \to \bigoplus_j H(l'_j,u'_j).
$$
Here, the first map is given as follows: If $u_i = \coprod_{\bullet} (\bullet,U_{\bullet})$, then $F(\lambda_i)(u_i) = \bigoplus_{\bullet} (\bigoplus V(U_{\bullet}))$, and on each component, the first map is given by the morphism $\coprod_{V(U)} (\bbs_{U, V(U)}, \rho(V(U)), V(U))$. The second map is given by 
$$
 \xymatrix{
 \bigoplus_i H(l'_{\psi(i)},F(\lambda_i)(u_i)) \cong \bigoplus_j \bigoplus_{\psi(i) = j} H(l'_{\psi(i)},F(\lambda_i)(u_i)) \ar[r]^{\hspace{3.5cm} \bigoplus_j \chi_j} & \bigoplus_j H(l'_j,u'_j),
 }
$$
where we view $\bigoplus_j \chi_j$ as a morphism in $\fB_C$. 

Our goal is to show that $H$ induces a homotopy equivalence of classifying spaces. To do so, let us first describe $\fF_{u,H}$ for $u \in \obj \fB_C$. Two objects $(v,\sigma)$ and $(w,\tau)$ in $u \backslash H$ are equivalent (with respect to the relation defining $\fF_{u,H}$) if there exists a (necessarily unique) morphism $a$ from $v$ to $w$ which is invertible in $\fH$, i.e., $a$ is given by $(\lambda_i,\psi,\chi_j)$ such that $\lambda_i = \id$ and $l'_{\psi(i)} = l_i$ for all $i$. Hence $\fF_{u,H}$ is a poset, where we define $[v,\sigma] \geq [w,\tau]$ if there exists a morphism from $(v,\sigma)$ to $(w,\tau)$ in $u \backslash H$ given by $f \in \fH(w,v)$ such that $H(f) \sigma = \tau$. We want to show that the classifying space of $\fF_{u,H}$ is contractible. This will follow from the next observation.
\blemma
\label{lem:FuHdirected}
Every two elements of $\fF_{u,H}$ have a common lower bound.
\elemma
\setlength{\parindent}{0cm} \setlength{\parskip}{0cm}

\bproof
We want to show that for all $[v,\sigma]$ and $[w,\tau]$, there exists $[x,\alpha]$ such that $[v,\sigma] \geq [x,\alpha]$ and $[w,\tau] \geq [x,\alpha]$. Indeed, suppose that $u = \coprod_i (i,U_i)$, $v = \bigoplus_{\bullet} (l_{\bullet},v_{\bullet})$, $H(v) = \coprod_j (j,V_j)$, and that $\sigma = \coprod_{j,i} (\bbs_{j,i},\rho(U_{j,i}),U_{j,i})$. Then $U_i = \coprod_j U_{j,i}$ and $V_j = \coprod_i U_{j,i}$. Find an index $l'$ with $l_{\bullet} \leq l'$ for all $\bullet$ such that $U_{j,i}$ can be written as a disjoint union $U_{j,i} = \coprod V_{j,i,\zeta}$ for some $V_{j,i,\zeta} \in \cU_{l'}$. Define $v' \defeq (l', \bigoplus_{j,i,\zeta} V_{j,i,\zeta})$. Construct a morphism $e$ in $\fH$ from $v$ to $v'$ given by data $(\lambda_{\bullet}, \psi, \chi)$ as above (note that $m=1$ for $v'$), with $\psi(\bullet) = 1$ for all $\bullet$, $\lambda_{\bullet}$ is the morphism $l_{\bullet} \to l'$ in $\fL$, which exists because $l_{\bullet} \leq l'$, and
$$
 \chi: \: \bigoplus_j \bigoplus V(V_j) \to \bigoplus_{j,i,\zeta} V_{j,i,\zeta}
$$
given by the obvious permutation, induced by the decomposition $V_j = \coprod_i U_{j,i} = \coprod_i \coprod V_{j,i,\zeta}$. By construction, we have
$$
 H(e) = \coprod_{j,i,\zeta} (\bbs_{(j,i,\zeta),j},\rho(V_{j,i,\zeta}),V_{j,i,\zeta}).
$$
Similarly, suppose that $w = \bigoplus_{\circ} (l_{\circ},w_{\circ})$, $H(w) = \coprod_k (k,W_k)$ and $\tau = \coprod_{k,i} (\bbs_{k,i},\rho(U_{k,i}),U_{k,i})$. Then $W_k = \coprod_i U_{k,i}$. Find an index $\ti{l}$ with $l_{\circ} \leq \ti{l}$ for all $\circ$ such that $U_{k,i}$ can be written as a disjoint union $U_{k,i} = \coprod W_{k,i,\eta}$ for some $W_{k,i,\eta} \in \cU_{\ti{l}}$. Define $\ti{w} \defeq (\ti{l}, \bigoplus_{k,i,\eta} W_{k,i,\eta})$ and construct a morphism $f$ in $\fH$ from $w$ to $\ti{w}$ as above such that
$$
 H(f) = \coprod_{k,i,\eta} (\bbs_{(k,i,\eta),k},\rho(W_{k,i,\eta}),W_{k,i,\eta}).
$$
Now find $l \in \fL$ such that $l', \ti{l} \leq l$. It follows that all $V_{j,i,\zeta}$ and $W_{k,i,\eta}$ can be written as disjoint unions of elements in $\cU_l$. This leads to morphisms $e'$ and $f'$ in $\fH$ such that
\begin{align*}
 H(e') &= \coprod_{k,j,i,\zeta,\eta} (\bbs_{(k,j,i,\zeta,\eta),(j,i,\zeta)},\rho(U_{k,j,i,\zeta,\eta}),U_{k,j,i,\zeta,\eta})\\
 H(f') &= \coprod_{k,j,i,\zeta,\eta} (\bbs_{(k,j,i,\zeta,\eta),(k,i,\eta)},\rho(\ti{U}_{k,j,i,\zeta,\eta}),\ti{U}_{k,j,i,\zeta,\eta})
\end{align*}
satisfying $\mft(e') = \mft(f')$ and $H(e') H(e) \sigma = H(f') H(f) \tau$. Setting $x \defeq \mft(e') = \mft(f')$ and $\alpha \defeq H(e') H(e) \sigma = H(f') H(f) \tau$, we indeed have $[v,\sigma] \geq [x,\alpha]$ and $[w,\tau] \geq [x,\alpha]$, as desired.
\eproof
\setlength{\parindent}{0cm} \setlength{\parskip}{0.5cm}

We obtain the following consequence of Corollary~\ref{cor:Phihomoequi} and Lemma~\ref{lem:FuHdirected} because posets with the property that any two elements have a common lower bound have contractible classifying spaces.
\bcor
\label{cor:H=homoequi}
$H$ induces a homotopy equivalence of classifying spaces.
\ecor

\bprop
\label{prop:H(BC)}
We have
$$
 H_*(\Kz(\fB_C),\sfC) \cong
 \bfa
 \C(C,\sfC) \oplus \sfC & \falls * = 0,\\
 \gekl{0} & \sonst.
 \efa
$$
\eprop
\setlength{\parindent}{0cm} \setlength{\parskip}{0cm}

\bproof
Corollary~\ref{cor:H=homoequi} implies that $H_*(\Kz(\fB_C)) \cong H_*(\Kz(\fH))$ because of \cite[Lemma~2.3]{Tho}. Thus we obtain
$$
 H_*(\Kz(\fB_C), \sfC) \cong H_*(\Kz(\fH), \sfC) \cong H_*({\rm hocolim}_l \, \Kz(\fB_l), \sfC) \cong \ilim_l H_*(\Kz(\fB_l), \sfC).
$$
Here we used \cite[Theorem~4.1]{Tho} for the second isomorphism, and we obtain the third isomorphism because we are taking a filtered colimit. Thus, applying Lemma~\ref{lem:H(K(Bl))}, we derive
$$
 H_*(\Kz(\fB_C), \sfC)
 \cong
 \bfa
 \ilim_l (\bigoplus_{\cU_l} \sfC) \oplus \sfC \cong \C(C,\sfC) \oplus \sfC& \falls * = 0,\\
 \gekl{0} & \sonst,
 \efa
$$
as desired.
\eproof
\setlength{\parindent}{0cm} \setlength{\parskip}{0.5cm}

Now suppose that $O \in \crO^{(\nu - 1)}$. Let $\fB_O$ be the full subcategory of $\fB_{\Gn \acts \Gnum}$ whose objects are of the form $\coprod_i (i,U_i)$, where $U_i \subseteq O$. 
\bprop
\label{prop:H(BO)}
We have
$$
 H_*(\Kz(\fB_O),\sfC) \cong
 \bfa
 \C(O,\sfC) \oplus \sfC & \falls * = 0,\\
 \gekl{0} & \sonst.
 \efa
$$
\eprop
\setlength{\parindent}{0cm} \setlength{\parskip}{0cm}

\bproof
Let us use the same notation as in \S~\ref{ss:AlgKSmallPermCat}. Since $O = \bigcup_U U$, where $U$ runs through all compact open subspaces of $O$, we conclude that for all $p, q$, we have $\fN_p \fB_O(S^n_q) = \bigcup_U \fN_p \fB_U(S^n_q)$. Hence, after taking diagonals, we obtain $\fN \fB_O(S^n) = \bigcup_U \fN \fB_U(S^n)$. Since homology is compatible with inductive limits, we obtain $H_*(X_n) \cong \ilim_U H_*(X_{U,n})$, where $X_n$ is the $n$-th simplicial set of $\Kz(\fB_O)$ and $X_{U,n}$ is the $n$-th simplicial set of $\Kz(\fB_U)$. Therefore, using Proposition~\ref{prop:H(BC)}, we conclude that
$$
 H_*(\Kz(\fB_O),\sfC) \cong \ilim_U H_*(\Kz(\fB_U),\sfC) \cong \ilim_U \C(U,\sfC) \oplus \sfC \cong \C(O,\sfC) \oplus \sfC,
$$
for $* = 0$, and $H_*(\Kz(\fB_O),\sfC) \cong \gekl{0}$ for $* > 0$.
\eproof
\setlength{\parindent}{0cm} \setlength{\parskip}{0.5cm}

Now suppose that we are given $O_1, \dotsc, O_N \in \crO^{(\nu - 1)}$. Let $\fB_{\cap}$ be the full subcategory of $\fB_{\Gn \acts \Gnum}$ whose objects are of the form $\coprod_i (i,U_i)$, where $U_i \subseteq O_1 \cap O_n$ for some $2 \leq n \leq N$. Let $\fB_1$ be the full subcategory of $\fB_{\Gn \acts \Gnum}$ whose objects are of the form $\coprod_i (i,U_i)$, where $U_i \subseteq O_1$. Let $\fB_2$ be the full subcategory of $\fB_{\Gn \acts \Gnum}$ whose objects are of the form $\coprod_i (i,U_i)$, where $U_i \subseteq O_n$ for some $2 \leq n \leq N$. Finally, let $\fB$ be the full subcategory of $\fB_{\Gn \acts \Gnum}$ whose objects are of the form $\coprod_i (i,U_i)$, where $U_i \subseteq O_n$ for some $1 \leq n \leq N$. Let $\fL$ be the category with three objects $\mfl_{\cap}$, $\mfl_1$ and $\mfl_2$, the corresponding identity morphisms, one morphism $\mfl_{\cap} \to \mfl_1$ and another morphism $\mfl_{\cap} \to \mfl_2$. Let $F$ be the functor from $\fL$ to small permutative categories sending $\mfl_{\cap}$ to $\fB_{\cap}$, $\mfl_1$ to $\fB_1$, $\mfl_2$ to $\fB_2$, the morphism $\mfl_{\cap} \to l_1$ to the inclusion $\fB_{\cap} \into \fB_1$ and the morphism $\mfl_{\cap} \to l_2$ to the inclusion $\fB_{\cap} \into \fB_2$.

Let $\fP \defeq {\rm hocolim} \, F$. So objects of $\fP$ are of the form $(l_1,u_1) \oplus \dotso \oplus (l_n,u_n)$, where $l_i \in \obj \fL$ and $u_i \in F(l_i)$. Morphisms from $(l_1,u_1) \oplus \dotso \oplus (l_n,u_n)$ to $(l'_1,u'_1) \oplus \dotso \oplus (l'_m,u'_m)$ are given by $(\lambda_i, \psi, \chi_j)$, where $\psi: \: \gekl{1, \dotsc, n} \to \gekl{1, \dotsc, m}$ is a surjective map, $\lambda_i: \: l_i \to l'_{\psi(i)}$ are morphisms in $\fL$, and $\chi_j: \: \bigoplus_{\psi(i) = j} F(\lambda_i)(u_i) \to u'_j$ are morphisms in $F(l'_j)$. In our case, if $l_i \in \gekl{\mfl_1, \mfl_2}$, then we must have $l'_{\psi(i)} = l_i$ and $\lambda_i = \id$. 

Define a permutative functor $\Pi: \: \fP \to \fB$ by sending $\bigoplus_i (l_i,u_i)$ to $\bigoplus_i u_i$ and the morphism from $(l_1,u_1) \oplus \dotso \oplus (l_n,u_n)$ to $(l'_1,u'_1) \oplus \dotso \oplus (l'_m,u'_m)$ given by $(\lambda_i, \psi, \chi_j)$ to the morphism
$$
 \xymatrix{
 \bigoplus_i u_i = \bigoplus_j \bigoplus_{\psi(i) = j} u_i \ar[r]^{\hspace{1.25cm} \bigoplus_j \chi_j} & \bigoplus_j u'_j,
 }
$$ 
where we view $\chi_j$, which by definition is a morphism in $F(l'_j)$, as a morphism in $\fB$. This is possible because $F(l'_j) \subseteq \fB$.

Let us describe $\fF_{u,\Pi}$ for $u \in \obj \fB$. Two objects $(v,\sigma)$ and $(w,\tau)$ in $u \backslash \Pi$ are equivalent (with respect to the relation defining $\fF_{u,\Pi}$) if there exists a (necessarily unique) morphism $a$ from $v$ to $w$ which is invertible in $\fP$, i.e., $a$ is given by $(\lambda_i,\psi,\chi_j)$ such that $\lambda_i = \id$ and $l'_{\psi(i)} = l_i$ for all $i$. Hence $\fF_{u,\Pi}$ is a poset, where we define $[v,\sigma] \geq [w,\tau]$ if there exists a morphism from $(v,\sigma)$ to $(w,\tau)$ in $u \backslash \Pi$ given by $f \in \fP(w,v)$ such that $\Pi(f) \sigma = \tau$. We want to show that the classifying space of $\fF_{u,\Pi}$ is contractible.

\blemma
\label{lem:FuPidirected}
Every two elements in $\fF_{u,\Pi}$ have a common upper bound.
\elemma
\setlength{\parindent}{0cm} \setlength{\parskip}{0cm}

\bproof
Given $[v,\sigma]$ and $[w,\tau]$, we proceed as in the proof of Lemma~\ref{lem:FuHdirected} to obtain invertible morphisms $e, e', f, f'$ in $\fP$ such that $\Pi(e') \Pi(e) \sigma = \Pi(f') \Pi(f) \tau$. Let $\alpha \defeq \Pi(e') \Pi(e) \sigma = \Pi(f') \Pi(f) \tau$. Set $y \defeq \mft(e')$ and $z \defeq \mft(f')$. Then these objects $y$ and $z$ are of the form $y = (l_1(y),u_1) \oplus \dotso$ and $z = (l_1(z),u_1) \oplus \dotso$. Define another object $x$ of $\fP$ as follows: Set $l_{\bullet}(x) \defeq \mfl_{\cap}$ if $l_{\bullet}(y) = \mfl_{\cap}$ or $l_{\bullet}(z) = \mfl_{\cap}$ or $\mfl_{\cap} \neq l_{\bullet}(y) \neq l_{\bullet}(z) \neq \mfl_{\cap}$. Otherwise, set $l_{\bullet}(x) \defeq l_{\bullet}(y) = l_{\bullet}(z)$. Now define $x \defeq (l_1(x),u_1) \oplus \dotso$. This is well-defined, i.e., if $l_{\bullet}(x) = \mfl_{\cap}$ and $\mfl_{\cap} \neq l_{\bullet}(y) \neq l_{\bullet}(z) \neq \mfl_{\cap}$, then $u_{\bullet}$ must be an object of $\fB_1$ as well as an object of $\fB_n$ for some $2 \leq n \leq N$, so that $u_{\bullet} \in \obj \fB_{\cap}$. Furthermore, define a morphism $\ti{e}$ in $\fP$ from $x$ to $y$ by setting $\psi \defeq \id$, defining $\lambda_i$ as the obvious morphism $l_i(x) \to l_i(y)$ and $\chi_j \defeq \id$. Similarly, define a morphism $\ti{f}$ in $\fP$ from $x$ to $z$. Then $\Pi(\ti{e}) \alpha = \alpha$, $\Pi(\ti{f}) \alpha = \alpha$, so that $[v,\sigma] = [y,\alpha] \leq [x,\alpha]$ and $[w,\tau] = [z,\alpha] \leq [x,\alpha]$, as desired.
\eproof
\setlength{\parindent}{0cm} \setlength{\parskip}{0.5cm}

Since posets with the property that any two elements have a common upper bound have contractible classifying spaces, we obtain the following consequence of Corollary~\ref{cor:Phihomoequi} and Lemma~\ref{lem:FuPidirected}.
\bcor
\label{cor:Pi=homoequi}
$\Pi$ induces a homotopy equivalence of classifying spaces.
\ecor

We need the following observation: Let 
\begin{eqnarray*}
 \C_{\cap} &\defeq& \lspan \menge{c_U}{c \in \sfC, \, U \text{ compact open subspace of } O_1 \cap O_n \text{ for some } 2 \leq n \leq N},\\
 \C_1 &\defeq& \lspan \menge{c_U}{c \in \sfC, \, U \text{ compact open subspace of } O_1},\\
 \C_2 &\defeq& \lspan \menge{c_U}{c \in \sfC, \, U \text{ compact open subspace of } O_n \text{ for some } 2 \leq n \leq N}\\
 \C &\defeq& \lspan \menge{c_U}{c \in \sfC, \, U \text{ compact open subspace of } O_n \text{ for some } 1 \leq n \leq N}.
\end{eqnarray*}
All these are subspaces of $\C(\Gnum,\sfC)$.
\blemma
\label{lem:CC+CC}
The obvious inclusion maps fit into the following short exact sequence:
$$
 0 \to \C_{\cap} \to \C_1 \oplus \C_2 \to \C \to 0,
$$
where the map $\C_{\cap} \to \C_1 \oplus \C_2$ sends $f \in \C_{\cap}$ to $(f,-f) \in \C_1 \oplus \C_2$ and the map $\C_1 \oplus \C_2 \to \C$ sends $(f_1,f_2) \in \C_1 \oplus \C_2$ to $f_1 + f_2 \in \C$.
\elemma
\setlength{\parindent}{0cm} \setlength{\parskip}{0cm}

\bproof
Suppose that $\sum_i (c_i)_{U_i} + \sum_j (\ti{c}_j)_{V_j} = 0$ in $\C$, with $U_i \subseteq O_1$ and $V_j \subseteq O_n$ for some $2 \leq n \leq N$. After disjointifying, we may assume that the $U_i$ are pairwise disjoint. We must have $U_i \subseteq \bigcup_j V_j$, otherwise $c_i = 0$. As in the proof of Lemma~\ref{lem:C=Sum/I}, this allows us to replace $U_i$ by $W_k$, with $W_k \subseteq V_j \subseteq O_n$ for some $2 \leq n \leq N$, so that $W_k \subseteq O_1 \cap O_n$ for some $2 \leq n \leq N$. We conclude that
$$
 (\sum_i (c_i)_{U_i}, \sum_j (\ti{c}_j)_{V_j}) = (\sum_k (\bar{c}_k)_{W_k}, \sum_j (\ti{c}_j)_{V_j}) \equiv (0, \sum_k (\bar{c}_k)_{W_k} + \sum_j (\ti{c}_j)_{V_j})
$$
in $\C_1 \oplus \C_2$, where $\equiv$ is understood modulo the image of $\C_{\cap}$ under the map $\C_{\cap} \to \C_1 \oplus \C_2$ above. It follows that $\sum_k (\bar{c}_k)_{W_k} + \sum_j (\ti{c}_j)_{V_j} = 0$ in $\C$, i.e., $\sum_k (\bar{c}_k)_{W_k} + \sum_j (\ti{c}_j)_{V_j}$ is the constant zero function on $\Gnum$, and hence $\sum_k (\bar{c}_k)_{W_k} + \sum_j (\ti{c}_j)_{V_j} = 0$ in $\C_2$, as desired.
\eproof
\setlength{\parindent}{0cm} \setlength{\parskip}{0.5cm}

Given a finite subset $\underline{O} = \gekl{O_1, \dotsc, O_N} \subseteq \crO^{(\nu - 1)}$, let $\fB_{\underline{O}}$ be the full subcategory of $\fB_{\Gn \acts \Gnum}$ whose objects are of the form $\coprod_i (i,U_i)$, where $U_i \subseteq O_n$ for some $1 \leq n \leq N$. Also define
$$
 \C_{\underline{O}} \defeq \lspan \menge{c_U}{c \in \sfC, \, U \text{ compact open subspace of } O_n \text{ for some } 1 \leq n \leq N}.
$$
\bcor
\label{cor:H(BO1ON)}
We have
$$
 H_*(\Kz(\fB_{\underline{O}}),\sfC) \cong
 \bfa
 \C_{\underline{O}} \oplus \sfC & \falls * = 0,\\
 \gekl{0} & \sonst.
 \efa
$$
\ecor
\setlength{\parindent}{0cm} \setlength{\parskip}{0cm}

\bproof
We proceed inductively on $N$. The case $N=1$ is Proposition~\ref{prop:H(BO)}. Now we turn to the case of general $N$. Corollary~\ref{cor:Pi=homoequi} and \cite[Lemma~2.3, Theorem~4.1]{Tho} imply that $\Kz(\fB_{\underline{O}})$ can be identified with the homotopy pushout of
$$
 \xymatrix{
 \Kz(\fB_{\gekl{O_1 \cap O_2, \dotsc, O_1 \cap O_N}}) \ar[d] \ar[r] & \Kz(\fB_{\gekl{O_1}})\\
 \Kz(\fB_{\gekl{O_2, \dotsc, O_N}})
 }
$$ 
with respect to the maps of spectra induced by the canonical inclusions of small permutative categories (using functoriality, see \S~\ref{ss:AlgKSmallPermCat}). Thus, proceeding inductively, the long exact sequence in homology for pushouts (see for instance \cite[Example~3.7]{Tho}) yields the short exact sequence
$$
 0 \to H_0(\Kz(\fB_{\gekl{O_1 \cap O_2, \dotsc, O_1 \cap O_N}}),\sfC) \to H_0(\Kz(\fB_{\gekl{O_1}}),\sfC) \oplus H_0(\Kz(\fB_{\gekl{O_2, \dotsc, O_N}}),\sfC) \to H_0(\Kz(\fB_{\underline{O}}),\sfC) \to 0,
$$
and that $H_*(\Kz(\fB_{\underline{O}}),\sfC) \cong \gekl{0}$ for all $* > 0$. Comparing this short exact sequence with the one in Lemma~\ref{lem:CC+CC}, our claim follows.
\eproof
\setlength{\parindent}{0cm} \setlength{\parskip}{0.5cm}

\btheo
\label{thm:HBGnu}
For all $\nu \geq 1$, we have
$$
 H_*(\Kz(\fB_{G \acts \Gnu}),\sfC) \cong
 \bfa
 \C(\Gnum,\sfC) \oplus \sfC & \falls * = 0,\\
 \gekl{0} & \sonst.
 \efa
$$
\etheo
\setlength{\parindent}{0cm} \setlength{\parskip}{0cm}

\bproof
Corollary \ref{cor:I=homoequi} implies that $H_*(\Kz(\fB_{G \acts \Gnu}),\sfC) \cong H_*(\Kz(\fB_{\Gn \acts \Gnum}),\sfC)$ because of \cite[Lemma~2.3]{Tho}. A similar argument as for Proposition~\ref{prop:H(BO)} implies that $H_*(\Kz(\fB_{\Gn \acts \Gnum}),\sfC) \cong \ilim_{\underline{O}} H_*(\Kz(\fB_{\underline{O}}),\sfC)$, where we order $\underline{O}$ by inclusion. All in all, using Corollary~\ref{cor:H(BO1ON)}, we conclude that
\begin{equation*}
 H_*(\Kz(\fB_{G \acts \Gnu}),\sfC) 
 \cong \ilim_{\underline{O}} H_*(\Kz(\fB_{\underline{O}}),\sfC)
 \cong 
 \bfa
 \ilim_{\underline{O}} \C_{\underline{O}} \oplus \sfC \cong \C(\Gnum,\sfC) \oplus \sfC & \falls * = 0,\\
 \gekl{0} & \sonst. 
 \efa
 \qedhere
\end{equation*}
\eproof
\setlength{\parindent}{0cm} \setlength{\parskip}{0.5cm}

\subsection{Identifying homology of algebraic K-theory spectra with groupoid homology}

We introduce the notation $\fB^{(\nu)} \defeq \fB_{G \acts \Gnu}$. Let $\Delta_{\nu + 1}^{\mu}$ be the permutative functor $\fBnup \to \fBnu$ induced by the groupoid homomorphism $\id_G$ and the map $G \to \Gn, \, g \ma \rmr(g)$ for $\nu = 0$ and $\mu = 0$, and
$$
 \Gnup \to \Gnu,
 (g_0, \dotsc, g_{\nu}) \ma
 \bfa
 (g_0, \dotsc, g_{\mu} g_{\mu +1}, \dotsc, g_{\nu}) & \falls 0 \leq \mu \leq \nu - 1,\\
 (g_0, \dotsc, g_{\nu - 1}) & \falls \mu = \nu
 \efa
$$
for $\nu > 0$. It is straightforward to check that these maps satisfy the conditions in \S~\ref{sss:OpenEmb}. The functors $\Delta_{\nu + 1}^{\mu}$ induce maps of bisimplicial sets from $(p,q) \ma \fN_p \fBnup (S^n_q)$ to $(p,q) \ma \fN_p \fBnu (S^n_q)$ and hence maps
$$
 \partial_{\nu + 1}^{\mu}: \: C_{p,q} \fN_p \fBnup (S^n_q) \to C_{p,q} \fN_p \fBnu (S^n_q),
$$
for all $p, q$. Here we are using the same notation as in \S~\ref{ss:AlgKSmallPermCat} (and $C_{p,q} \fN_p \fBnu (S^n_q)$ denotes the total complex of $(p,q) \ma \fN_p \fBnu (S^n_q)$). We set $\partial_{\nu + 1} \defeq \sum_{\mu = 0}^{\nu} (-1)^{\mu} \partial_{\nu + 1}^{\mu}$. It is straightforward to check that, for all $p, q$, the sequence
$$
 \xymatrix{
 \dotso \ar[r]^{\hspace{-1.25cm} \partial_{\nu + 2}} & C_{p,q} \fN_p \fBnup (S^n_q) \ar[r]^{\partial_{\nu + 1}} & C_{p,q} \fN_p \fBnu (S^n_q) \ar[r]^{\hspace{1cm} \partial_{\nu}} & \dotso
 }
$$
forms a chain complex. Our goal is to show that this chain complex is exact.

\btheo
\label{thm:CCexact}
The following chain complex is exact:
$$
 \xymatrix{
 \dotso \ar[r]^{\hspace{-1.25cm} \partial_{\nu + 2}} & C_{p,q} \fN_p \fBnup (S^n_q) \ar[r]^{\partial_{\nu + 1}} & C_{p,q} \fN_p \fBnu (S^n_q) \ar[r]^{\hspace{1cm} \partial_{\nu}} & \dotso
 }
$$
\etheo

For the proof, let us construct maps $h_{\nu}: \: \fN_p \fBnu (S^n_q) \to \fN_p \fBnup (S^n_q)$. 

Let us first describe a general method which allows us to modify domains of morphisms. Suppose that $\sigma = \coprod_{j,i} (\bbs_{j,i}, \sigma_{j,i}, U_{j,i})$ is a morphism in $\fBnu$ with domain $\mfd(\sigma) = \coprod_i (i,U_i)$. Suppose that we want to replace $\mfd(\sigma)$ by $\mfd_{\sigma} = \coprod_i (i,U'_i) \in \obj \fBnup$, with $U'_i \subseteq \Gnup$ and $m$ restricting to a homeomorphism $m \vert_{U'_i}: \: U'_i \cong U_i$, where $m$ is the map $m: \: \Gnup \to \Gnu, \, (g_0, \dotsc, g_{\nu}) \ma (g_0 g_1, \dotsc, g_{\nu})$. Then define $h_{\nu}(\sigma, \mfd_{\sigma}) \defeq \coprod_{j,i} (\bbs_{j,i}, \sigma_{j,i}, U'_{j,i})$, where $U'_{j,i} \subseteq U'_i$ such that $m(U'_{j,i})= U_{j,i}$. Now suppose that $\tau$, $\sigma$ are two composable morphisms in $\fBnu$. Given $\mfd_{\sigma}$ as above, construct $h_{\nu}(\sigma, \mfd_{\sigma})$ and set $\mfd_{\tau} \defeq \mft(h_{\nu}(\sigma, \mfd_{\sigma}))$. Then $\mfd_{\tau}$ is of the same form for $\mfd(\tau)$, so that we can form $h_{\nu}(\tau, \mfd_{\tau})$ in the same way as before. By construction, the following holds:
$$
 h_{\nu}(\tau \sigma, \mfd_{\sigma}) = h_{\nu}(\tau, \mfd_{\tau}) h_{\nu}(\sigma, \mfd_{\sigma}).
$$
Now let us apply the procedure above and construct $h_{\nu}$ on $\fBnu(A)$ for some finite based set $A$. Given an object $(u_S, \varphi_{T,T'})$ in $\fBnu(A)$ and an element $a \in A$ which is not the base point, write $u_{\gekl{a}} = \coprod_i (i, U_i)$ and set $h_{\nu}(u_{\gekl{a}}) \defeq \coprod_i (i, U'_i)$, where $U'_i = \menge{(\rho(z),z)}{z \in U_i}$ if $\nu \geq 1$, and $U'_i = U_i$ viewed as a subspace of $G$ via the canonical embedding $\Gn \into G$ if $\nu = 0$. Now proceed inductively on $\# S$. If $S = \gekl{a} \cup S'$, then define $h_{\nu}(u_S)$ as the target of the morphism $h_{\nu}(\varphi_{\gekl{a},S'}, h_{\nu}(u_{\gekl{a}}) \oplus h_{\nu}(u_{S'}))$. Given disjoint $T, T' \subseteq A$ not containing the base point, set $h_{\nu}(\varphi_{T,T'}) \defeq h_{\nu}(\varphi_{T,T'}, h_{\nu}(u_T) \oplus h_{\nu}(u_{T'}))$. It is straightforward to check that $\gekl{h_{\nu}(u_S), h_{\nu}(\varphi_{T,T'})}$ lies in $\fBnup(A)$. Now we extend $h_{\nu}$ to $\fN_p \fBnu(A)$. Suppose we are given $(f_1, \dotsc, f_p) \in \fN_p \fBnu(A)$, with $f_{\sqcup} = \gekl{(f_{\sqcup})_S}$. Write $\mfd(f_p) = \gekl{u_S, \varphi_{T,T'}}$. Then set
$$
 h_{\nu}(f_p)_S \defeq h_{\nu}((f_p)_S,h_{\nu}(u_S)),
$$
and define recursively
$$
 h_{\nu}(f_{p'})_S \defeq h_{\nu}((f_{p'})_S, \mfd(h_{\nu}(f_{p' + 1})_S)).
$$
It is now straightforward to check that $h_{\nu}(f_{p'}) \defeq \gekl{(f_{p'})_S}$ defines a morphism in $\fBnu(A)$, and that
$$
 (h_{\nu}(f_1), \dotsc, h_{\nu}(f_p)) \in \fN_p \fBnup(A).
$$
Apply the above construction to $A = S^n_q$, denote by $h_{\nu}$ the map $h_{\nu}: \: \fN_p \fBnu(S^n_q) \to \fN_p \fBnup(S^n_q)$, and let $(h_{\nu})_*$ be the induced map $C_{p,q} \fN_p \fBnu (S^n_q) \to C_{p,q} \fN_p \fBnup (S^n_q)$. 
\bprop
\label{prop:h}
$(h_{\nu})_*$, $\nu \geq 0$, defines a chain homotopy between the identity map and the zero map, i.e., we have $\partial_1 (h_0)_* = \id$ on $C_{p,q} \fN_p \fB^{(0)}(S^n_q)$ and $\partial_{\nu + 1} (h_{\nu})_* + (h_{\nu + 1})_* \partial_{\nu} = \id$ on $C_{p,q} \fN_p \fBnu(S^n_q)$ for all $\nu \geq 1$.
\eprop
\setlength{\parindent}{0cm} \setlength{\parskip}{0cm}

\bproof
First observe that $\partial_{\nu + 1}^{\mu} (h_{\nu})_*$ and $(h_{\nu + 1})_* \partial_{\nu}^{\mu}$ are determined by how they act on $\Gnu$, in the sense that for both $\partial_{\nu + 1}^{\mu} (h_{\nu})_*$ and $(h_{\nu + 1})_* \partial_{\nu}^{\mu}$, there are maps $\Gnu \to \Gnu$, say $\zeta_{\nu}^{\mu}$ and $\eta_{\nu}^{\mu}$, such that, on the level of $\fBnu(S_q^n)$, for $a \in S_q^n$, $u_{\gekl{a}} = \coprod_i (i,U_i)$ is sent to $\coprod_i (i,\zeta_{\nu}^{\mu}(U_i))$ for $\partial_{\nu + 1}^{\mu} (h_{\nu})_*$ and to $\coprod_i (i,\eta_{\nu}^{\mu}(U_i))$ for $(h_{\nu + 1})_* \partial_{\nu}^{\mu}$, and that $\partial_{\nu + 1}^{\mu} (h_{\nu})_* = (h_{\nu + 1})_* \partial_{\nu}^{\ti{\mu}}$ if $\zeta_{\nu}^{\mu} = \eta_{\nu}^{\ti{\mu}}$.
\setlength{\parindent}{0cm} \setlength{\parskip}{0.5cm}

For $\nu = 0$, it is straightforward to see that $\zeta_0^0 = \id_{\Gn}$, so that $\partial_1 (h_0)_* = \id$. For $\nu \geq 1$, $\zeta_{\nu}^{\mu}$ is given by
$$
 (g_1, \dotsc, g_{\nu}) \ma (\rmr(g_1), g_1, \dotsc, g_{\nu}) 
 \ma
 \bfa
 (g_1, \dotsc, g_{\nu}) & \falls \mu = 0,\\
 (\rmr(g_1), \dotsc, g_{\mu} g_{\mu + 1}, \dotsc, g_{\nu}) & \falls 1 \leq \mu \leq \nu - 1,\\
 (\rmr(g_1), \dotsc, g_{\nu - 1}) & \falls \mu = \nu.
 \efa
$$
At the same time, $\eta_{\nu}^{\ti{\mu}}$ is given by
$$
 (g_1, \dotsc, g_{\nu}) 
 \ma 
 \bfa 
 (\dotsc, g_{\ti{\mu}} g_{\ti{\mu} + 1}, \dotsc. g_{\nu}) & \falls 0 \leq \ti{\mu} \leq \nu - 2,\\
 (g_1, \dotsc, g_{\nu - 1}) & \falls \ti{\mu} = \nu - 1
 \efa
 \ma
 \bfa
 (\rmr(g_1), \dotsc, g_{\ti{\mu}} g_{\ti{\mu} + 1}, \dotsc, g_{\nu}) & \falls 0 \leq \ti{\mu} \leq \nu - 2,\\
 (\rmr(g_1), \dotsc, g_{\nu - 1}) & \falls \ti{\mu} = \nu - 1.
 \efa
$$
Therefore, the computations above show that when we compute $\partial_{\nu + 1} (h_{\nu})_* + (h_{\nu + 1})_* \partial_{\nu}$, then all terms cancel except the identity term (corresponding to $\mu = 0$).
\eproof
\setlength{\parindent}{0cm} \setlength{\parskip}{0.5cm}

The following is an immediate consequence of Proposition~\ref{prop:h}.
\bcor
The chain complex
$$
 \xymatrix{
 \dotso \ar[r]^{\hspace{-1.25cm} \partial_{\nu + 2}} & C_{p,q} \fN_p \fBnup (S^n_q) \ar[r]^{\partial_{\nu + 1}} & C_{p,q} \fN_p \fBnu (S^n_q) \ar[r]^{\hspace{1cm} \partial_{\nu}} & \dotso
 }
$$
is homotopy equivalent to the zero chain complex.
\ecor
In particular, the chain complex
$$
 \xymatrix{
 \dotso \ar[r]^{\hspace{-1.25cm} \partial_{\nu + 2}} & C_{p,q} \fN_p \fBnup (S^n_q) \ar[r]^{\partial_{\nu + 1}} & C_{p,q} \fN_p \fBnu (S^n_q) \ar[r]^{\hspace{1cm} \partial_{\nu}} & \dotso
 }
$$
is exact. This proves Theorem~\ref{thm:CCexact}.

\btheo
\label{thm:HKB=HG}
Let $G$ be an ample groupoid with locally compact Hausdorff unit space $\Gn$ and $\sfC$ an abelian group. We have $\ti{H}_*(\Kz(\fB_G),\sfC) \cong H_*(G,\sfC)$.
\etheo
\setlength{\parindent}{0cm} \setlength{\parskip}{0cm}

\bproof
Theorem~\ref{thm:CCexact} implies that
$$
 \xymatrix{
 \dotso \ar[r]^{\hspace{-1.25cm} \partial_{\nu + 2}} & C_{p,q} \fN_p \fBnup (S^n_q) \ar[r]^{\partial_{\nu + 1}} & C_{p,q} \fN_p \fBnu (S^n_q) \ar[r]^{\hspace{1cm} \partial_{\nu}} & \dotso
 }
$$
is a long exact sequence. After forming diagonals, we obtain the long exact sequence
$$
 \xymatrix{
 \dotso \ar[r]^{\hspace{-1cm} \partial_{\nu + 2}} & C_q \fN \fBnup (S^n) \ar[r]^{\partial_{\nu + 1}} & C_q \fN \fBnu (S^n) \ar[r]^{\hspace{0.75cm} \partial_{\nu}} & \dotso
 }
$$
Let $\ker^n_q(\partial_0) \defeq C_q \fN \fB^{(0)}(S^n) = C_q \fN \fB_G(S^n)$ and, for $\nu \geq 1$, let $\ker^n_q(\partial_{\nu})$ be the kernel of
$$
 \xymatrix{
 C_q \fN \fBnu (S^n) \ar[r]^{\hspace{-0.25cm} \partial_{\nu}} & C_q \fN \fBnum (S^n).
 }
$$
In this way, we obtain short exact sequences
$$
 \xymatrix{
 0 \ar[r] & \ker^n_q(\partial_{\nu + 1}) \ar[r] & C_q \fN \fBnup (S^n) \ar[r]^{\hspace{0.5cm} \partial_{\nu + 1}} & \ker^n_q(\partial_{\nu}) \ar[r] & 0.
 }
$$
These are actually short exact sequences of chain complexes with respect to $q$. Taking homology with respect to $q$, we obtain long exact sequences
$$
 \xymatrix{
 \dotso \ar[r] & H_*(\ker^n_*(\partial_{\nu + 1}),\sfC) \ar[r] & H_*(C_* \fN \fBnup (S^n),\sfC) \ar[r]^{\hspace{0.5cm} (\partial_{\nu + 1})_*} & H_*(\ker^n_*(\partial_{\nu}),\sfC) \ar[r] & \dotso.
 }
$$
Taking the inductive limit for $n \to \infty$, we obtain the long exact sequence
$$
 \xymatrix{
 \dotso \ar[r] & H_*(\ker(\partial_{\nu + 1}),\sfC) \ar[r] & H_*(\Kz(\fBnup),\sfC) \ar[r]^{\hspace{0.25cm} (\partial_{\nu + 1})_*} & H_*(\ker(\partial_{\nu}),\sfC) \ar[r] & \dotso.
 }
$$
Here $H_*(\ker(\partial_{\nu}),\sfC) \defeq \ilim_n H_*(\ker^n_*(\partial_{\nu}),\sfC)$.
\setlength{\parindent}{0cm} \setlength{\parskip}{0.5cm}

By Theorem~\ref{thm:HBGnu}, $H_*(\Kz(\fBnup),\sfC) \cong \gekl{0}$ for all $* \geq 1$. Hence we obtain
$$
 H_*(\ker(\partial_{\nu}),\sfC) \cong H_{*-1}(\ker(\partial_{\nu + 1}),\sfC)
$$
for all $* \geq 2$ and all $\nu \geq 0$. 

This yields, for all $* \geq 2$, 
\begin{equation}
\label{e:H=H=H1}
 H_*(\Kz(\fB_G),\sfC) = H_*(\ker(\partial_0),\sfC) \cong H_{*-1}(\ker(\partial_1),\sfC) \cong \dotso \cong H_1(\ker(\partial_{* - 1}),\sfC).
\end{equation}
Moreover, we obtain that for all $* \geq 1$,
\begin{equation}
\label{e:H1H0H0}
 0 \to H_1(\ker(\partial_{*-1}),\sfC) \to H_0(\ker(\partial_*),\sfC) \to H_0(\Kz(\fB^{(*)}),\sfC)
\end{equation}
is exact.

In addition, because for all $* \geq 0$, the sequence
$$
 \xymatrix{
 C_q \fN \fB^{(* + 2)} (S^n) \ar[r]^{\partial_{* + 2}} & C_q \fN \fB^{(*+1)} (S^n) \ar[r]^{\hspace{0.5cm} \partial_{* + 1}} &  \ker^n_q(\partial_*) \ar[r] & 0
 }
$$
is exact, we obtain the exact sequence
$$
 \xymatrix{
 H_0(C_{\star} \fN \fB^{(*+2)} (S^n),\sfC) \ar[r]^{H_0(\partial_{*+2})} & H_0(C_{\star} \fN \fB^{(*+1)} (S^n),\sfC) \ar[r]^{\hspace{0.5cm} H_0(\partial_{*+1})} & H_0(\ker^n_{\star}(\partial_*),\sfC) \ar[r] & 0,
 }
$$
where we have taken the $0$-th homology with respect to the index $\star$.

Taking the inductive limit for $n \to \infty$, we conclude that the sequence
$$
 \xymatrix{
 H_0(\Kz(\fB^{(*+2)}),\sfC) \ar[r]^{H_0(\partial_{*+2})} & H_0(\Kz(\fB^{(*+1)}),\sfC) \ar[r]^{\hspace{0.25cm} H_0(\partial_{*+1})} & H_0(\ker(\partial_*),\sfC) \ar[r] & 0
 }
$$
is exact. It follows that $H_0(\partial_{*+1})$ induces an isomorphism
\begin{equation}
\label{e:cokerH0=H0ker}
 \coker H_0(\partial_{*+2}) \cong H_0(\ker(\partial_*),\sfC),
\end{equation}
for all $* \geq 0$. 

Plugging \eqref{e:cokerH0=H0ker} into \eqref{e:H1H0H0}, we conclude that for all $* \geq 1$, the sequence
$$
 \xymatrix{
 0 \ar[r] & H_1(\ker(\partial_{*-1}),\sfC) \ar[r] & \coker H_0(\partial_{*+2}) \ar[r]^{H_0(\partial_{*+1})} & H_0(\Kz(\fB^{(*)}),\sfC)
 }
$$
is exact. Hence
\begin{equation}
\label{e:H1=H}
 H_1(\ker(\partial_{*-1}),\sfC) \cong \rukl{\ker H_0(\partial_{*+1})} / \rukl{\coker H_0(\partial_{*+2})} 
 \cong H_*(\gekl{H_0(\Kz(\fB^{(*+1)}),\sfC), H_0(\partial_{*+1})}),
\end{equation}
for all $* \geq 1$. 

Combining \eqref{e:H=H=H1}, \eqref{e:H1=H}, \eqref{e:cokerH0=H0ker} and applying Theorem~\ref{thm:HBGnu}, we conclude that
$$
 H_*(\Kz(\fB_G),\sfC) = H_*(\ker(\partial_0),\sfC) \cong H_*(\gekl{\C(G^{(*)},\sfC) \oplus \sfC, \, \dot{\partial}_{*+1}}),
$$
i.e., $H_*(\Kz(\fB_G),\sfC)$ is given by the homology of the chain complex
$$
 \xymatrix{
 \dotso \ar[r]^{\hspace{-1cm} \dot{\partial}_4} & \C(G^{(2)},\sfC) \oplus \sfC \ar[r]^{\dot{\partial}_3} & \C(G^{(1)},\sfC) \oplus \sfC \ar[r]^{\dot{\partial}_2} & \C(G^{(0)},\sfC) \oplus \sfC \ar[r] & 0,
 }
$$
where $\dot{\partial}_{*+1}$ denotes the map induced by $\partial_{*+1}$ on $H_0(\Kz(\fB^{(*+1)}),\sfC) \cong \C(G^{(*)},\sfC) \oplus \sfC$. Now $\dot{\partial}_{*+1}$ is of the form $\bar{\partial}_* \oplus 0$ for all $* \geq 1$ odd and $\bar{\partial}_* \oplus \id_{\sfC}$ for all $* \geq 1$ even, and it is straightforward to identify $\bar{\partial}_*$ with the map $\ti{\partial}_*$ defined in \eqref{e:tipartial}. Hence, all in all, we conclude that
\begin{equation*}
 \ti{H}_*(\Kz(\fB_G),\sfC) \cong H_*(\gekl{\C(G^{(*)},\sfC), \, \ti{\partial}_*}) \cong H_*(G,\sfC). \qedhere
\end{equation*}
\eproof
\setlength{\parindent}{0cm} \setlength{\parskip}{0.5cm}

\section{Homology of topological full groups in terms of homology of infinite loop spaces}

Let $G$ be an ample groupoid with locally compact Hausdorff unit space $\Gn$. In the following, we write $\fB \defeq \fB_G$. Our goal is to describe the homology of the topological full group $\bmF(G)$ in terms of the homology of the infinite loop space attached to $\Kz(\fB)$.

\subsection{Simplicial complexes attached to small permutative categories of bisections}

We will work in the setting of \cite{RWW} (see also \cite{SW}). Let us first introduce some relevant notions (Remark~\ref{rem:V=V,Q=Q} explains the motivation for the following subcategories).

Define the subcategory $\textbf{V}_{\fB}$ of $\fB$ by setting $\obj \textbf{V}_{\fB} \defeq \obj \fB$ and, for $u = \coprod_{i \in I} (i,U_i), \, v = \coprod_{j \in J} (j,V_j) \in \obj \textbf{V}_{\fB}$, $\textbf{V}_{\fB}(v,u) \defeq \menge{\pi \sigma}{\pi \in \fP(v,u), \, \sigma \in \fB(u,u)}$, where $\fP(v,u)$ are morphisms of the form $\pi = \coprod_i (\bbs_{j(i),i}, U_i)$ for a bijection $i \ma j(i)$ between $I$ and $J$. In other words, morphisms in $\fP(v,u)$ are just given by permutations, and for objects $u$, $v$, there is a morphism from $u$ to $v$ in $\textbf{V}_{\fB}$ only if $u$ and $v$ are equal up to permutation. $\textbf{V}_{\fB}$ inherits the structure of a small permutative category from $\fB$.

Moreover, we define a category $\textbf{Q}_{\fB}$ as follows. Set $\obj \textbf{Q}_{\fB} \defeq \obj \textbf{V}_{\fB} = \obj \fB$. Given $u, v \in \obj \textbf{Q}_{\fB}$, define the set of morphisms $\textbf{Q}_{\fB}(v,u)$ to be equivalence classes of pairs $(u',\sigma)$, where $u' \in \obj \textbf{Q}_{\fB}$ and $\sigma \in \textbf{V}_{\fB}(v,u' \oplus u)$, with respect to the equivalence relation that $(u'_1,\sigma_1) \sim (u'_2,\sigma_2)$ if there exists $\tau \in \textbf{V}_{\fB}(u'_2,u'_1)$ such that $\sigma_2 (\tau \oplus u) = \sigma_1$, i.e., the diagram
$$
 \xymatrix{
 u'_1 \oplus u \ar[d]_{\tau \oplus u} \ar[r]^{\sigma_1} & v\\
 u'_2 \oplus u \ar[ur]_{\sigma_2} &
 }
$$
commutes, where $u$ is the identity morphism at $u$. Let us now define composition of morphisms in $\textbf{Q}_{\fB}$. Given $[u',\sigma] \in \textbf{Q}_{\fB}(v,u)$ and $[v',\tau] \in \textbf{Q}_{\fB}(w,v)$, define $[v',\tau] [u',\sigma] \in \textbf{Q}_{\fB}(w,u)$ as $[v' \oplus u', \tau (v' \oplus \sigma)]$, where $\tau (v' \oplus \sigma)$ is the composite
$$
 \xymatrix{
 v' \oplus u' \oplus u \ar[r]^{\quad v' \oplus \sigma} & v' \oplus v \ar[r]^{\tau} & w.
 }
$$
We want to define the structure of a small permutative category on $\textbf{Q}_{\fB}$. To this end, we define a functor $\oplus$. On objects, $\oplus$ acts just in the same way as in $\fB$ or $\textbf{V}_{\fB}$. Given $[u'_1,\sigma_1] \in \textbf{Q}_{\fB}(v_1,u_1)$ and $[u'_2,\sigma_2] \in \textbf{Q}_{\fB}(v_2,u_2)$, define $[u'_1,\sigma_1] \oplus [u'_2,\sigma_2] \in \textbf{Q}_{\fB}(v_1 \oplus v_2, u_1 \oplus u_2)$ as the morphism $[u'_1 \oplus u'_2, \tau]$, where $\tau$ is given by the composite
$$
 \xymatrix{
 u'_1 \oplus u'_2 \oplus u_1 \oplus u_2 \ar[rr]^{u'_1 \oplus \pi_{u'_2,u_1} \oplus u_2} & & u'_1 \oplus u_1 \oplus u'_2 \oplus u_2 \ar[r]^{\hspace{0.75cm} \sigma_1 \oplus \sigma_2} & v_1 \oplus v_2.
 }
$$
It is straightforward to check that $(\textbf{Q}_{\fB},\oplus)$ is a small permutative category and that $\emptyset$ is the unit with respect to $\oplus$.

\bremark
\label{rem:V=V,Q=Q}
In the special case of the groupoid $G = \cR_r \times G_k$ whose topological full group is isomorphic to the Higman-Thompson group $V_{k,r}$ (see \S~\ref{ss:ExTopFullGrp}), the categories we just introduced already appeared in \cite{SW}, where they are described using the language of Cantor algebras: $\textbf{Cantor}\reg$ in \cite{SW} is the restriction of $\fB$ to objects of the form $(\Gn)^{\oplus m}$, $\textbf{V}$ in \cite{SW} is the restriction of $\textbf{V}_{\fB}$ to objects of the form $(\Gn)^{\oplus m}$, and $\textbf{Q}$ in \cite{SW} is the restriction of $\textbf{Q}_{\fB}$ to objects of the form $(\Gn)^{\oplus m}$, for the groupoid $G = \cR_r \times G_k$. Also note that $\textbf{V}_{\fB}$ plays the role of $\cG$ in \cite{RWW} and $\textbf{Q}_{\fB}$ plays the role of $U \cG$ in \cite{RWW}.
\eremark

Let us now verify several conditions from \cite{RWW}. 
\blemma
\label{lem:i,ii,C}
The condition in \cite[Proposition~1.7~(i)]{RWW} is satisfied in $\textbf{V}_{\fB}$. Moreover, the condition in \cite[Proposition~1.7~(ii)]{RWW} is satisfied in $\textbf{V}_{\fB}$. Furthermore, condition $\textbf{C}$ in \cite[Definition~1.9]{RWW} is satisfied in $\textbf{V}_{\fB}$.
\elemma
\setlength{\parindent}{0cm} \setlength{\parskip}{0cm}

\bproof
The condition in \cite[Proposition~1.7~(i)]{RWW} holds because $\textbf{V}_{\fB}(\emptyset,\emptyset) = \gekl{\id_{\emptyset}}$. The condition in \cite[Proposition~1.7~(ii)]{RWW} holds because, for all $u, v \in \obj \textbf{V}_{\fB}$, $u \oplus v \cong \emptyset$ in $\textbf{V}_{\fB}$ implies that $u = \emptyset$ and $v = \emptyset$. Finally, condition $\textbf{C}$ in \cite[Definition~1.9]{RWW} holds because, for all $u, v, w \in \obj \textbf{V}_{\fB}$, $u \oplus w \cong v \oplus w$ in $\textbf{V}_{\fB}$ implies that $u \cong v$ in $\textbf{V}_{\fB}$.
\eproof
\setlength{\parindent}{0cm} \setlength{\parskip}{0.5cm}

\blemma
\label{lem:LS}
Conditions $\textbf{LS1}$ and $\textbf{LS2}$ in \cite[Definition~2.5]{RWW} are satisfied for all $a \in \obj \textbf{Q}_{\fB}$ and $\emptyset \neq u \in \obj \textbf{Q}_{\fB}$.
\elemma
\setlength{\parindent}{0cm} \setlength{\parskip}{0cm}

\bproof
First we need the following terminology: Given $o \in \obj \textbf{Q}_{\fB}$, define $\iota_o$ as the unique element of $\textbf{Q}_{\fB}(o,\emptyset)$ given by $\iota_o = [o,o]$. 
\setlength{\parindent}{0cm} \setlength{\parskip}{0.5cm}

To verify condition $\textbf{LS1}$, we have to show that $\iota_a \oplus u \oplus \iota_u \neq \iota_{a \oplus u} \oplus u$ in $\textbf{Q}_{\fB}(a \oplus u \oplus u, u)$. Indeed, $\iota_a \oplus u \oplus \iota_u = [a \oplus u, \sigma]$, where $\sigma$ is the morphism $a \oplus \pi_{u,u}$ from $(a \oplus u) \oplus u$ to $a \oplus u \oplus u$, while $\iota_{a \oplus u} \oplus u = [a \oplus u, \tau]$, where $\tau$ is the identity morphism from $(a \oplus u) \oplus u = a \oplus u \oplus u$ to $a \oplus u \oplus u$. If $\iota_a \oplus u \oplus \iota_u = \iota_{a \oplus u} \oplus u$ would hold in $\textbf{Q}_{\fB}(a \oplus u \oplus u, u)$, then there would be $\rho \in \textbf{V}_{\fB}(a \oplus u, a \oplus u)$ such that $\tau = \sigma (\rho \oplus u)$, i.e., the diagram
$$
 \xymatrix{
 a \oplus u \oplus u \ar[d]_{\rho \oplus u} \ar[r]^{\tau} & a \oplus u \oplus u\\
 a \oplus u \oplus u \ar[ur]_{\sigma} &
 }
$$
commutes. But that is impossible.

To verify condition $\textbf{LS2}$, we show that for all $u, v, w 	\in \obj \textbf{Q}_{\fB}$ with $u \neq \emptyset$, the map 
$$
 \textbf{Q}_{\fB}(w,v) \to \textbf{Q}_{\fB}(w \oplus u, v), \, [v', \sigma] \ma [v', \sigma] \oplus \iota_u
$$
is injective. Note that $[v', \sigma] \oplus \iota_u = [v' \oplus u, \ti{\sigma}]$, where $\ti{\sigma}$ is the composite
$$
 \xymatrix{
 (v' \oplus u) \oplus v \ar[r]^{\hspace{0.25cm} v' \oplus \pi_{u,v}} & v' \oplus v \oplus u \ar[r]^{\hspace{0.25cm} \sigma \oplus u} & w \oplus u.
 }
$$
Now assume that $[v'_1, \sigma_1], [v'_2, \sigma_2] \in \textbf{Q}_{\fB}(w,v)$ satisfy $[v'_1, \sigma_1] \oplus \iota_u = [v'_2, \sigma_2] \oplus \iota_u$. As above, let $[v'_1, \sigma_1] \oplus \iota_u = [v'_1 \oplus u, \ti{\sigma}_1]$ and $[v'_2, \sigma_2] \oplus \iota_u = [v'_2 \oplus u, \ti{\sigma}_2]$. Then there must exist $\rho \in \textbf{V}_{\fB}(v'_2 \oplus u, v'_1 \oplus u)$ such that
$$
 \ti{\sigma}_1 = \ti{\sigma}_2 (\rho \oplus v),
$$
i.e, the diagram
$$
 \xymatrix{
 v'_1 \oplus u \oplus v \ar[d]_{\rho \oplus v} \ar[r]^{\ti{\sigma}_1} & w \oplus u\\
 v'_1 \oplus u \oplus u \ar[ur]_{\ti{\sigma}_2} &
 }
$$
commutes. But then $\rho$ must be of the form $\tau \oplus u$ for some $\tau \in \textbf{V}_{\fB}(v'_2,v'_1)$, with $\sigma_2 (\tau \oplus v) = \sigma_1$. Hence $[v'_1, \sigma_1] = [v'_2, \sigma_2]$ in $\textbf{Q}_{\fB}(w,v)$, as desired.
\eproof
\setlength{\parindent}{0cm} \setlength{\parskip}{0.5cm}

Now let $a, u \in \obj \textbf{Q}_{\fB}$ with $u \neq \emptyset$, and fix $r \geq 1$. Define a semi-simplicial set $W$ by setting $W_p \defeq \textbf{Q}_{\fB}(a \oplus u^{\oplus r}, u^{\oplus (p+1)})$ for $0 \leq p \leq r-1$ and defining the face maps as $d_i: \: W_p \to W_{p-1}, \, [u', \sigma] \ma [u', \sigma] (u^{\oplus i} \oplus \iota_u \oplus u^{\oplus(p-i)})$ for $0 \leq i \leq p$. 

Moreover, let $S$ be the simplicial complex with vertices given by $W_0 = \textbf{Q}_{\fB}(a \oplus u^{\oplus r}, u)$, and $\gekl{[u'_i,\sigma_i]}_i$ forms a simplex of $S$ if there exists a simplex of $W$ with vertices $\gekl{[u'_i,\sigma_i]}_i$.

Given $[u',\sigma] \in W_p$, then the vertices of $[u',\sigma]$ are given by
$$
 \xymatrix{
 u \ar[r]^{\hspace{-0.5cm} \iota_i} & u^{\oplus(p+1)} \ar[r]^{[u',\sigma]} & a \oplus u^{\oplus r}
 }
$$
in $\textbf{Q}_{\fB}$, for $1 \leq i \leq p+1$. Here $\iota_i$ is the element $[u^{\oplus p}, u^{\oplus(i-1)} \oplus \pi_{u^{\oplus(p-i+1)}, u}] \in \textbf{Q}_{\fB}(u^{\oplus(p+1)},u)$. Therefore, vertices of $[u',\sigma]$ are determined by $\sigma \epsilon_i$, where $\epsilon_i \subseteq \cR \times \Gn$ is the compact open bisection with $\rms(\epsilon_i) = u$ and $\rmr(\epsilon_i)$ is the $i$-th summand of $u$ in $u' \oplus u^{\oplus (p+1)}$.

We have the following result from \cite{RWW}.
\btheo
\label{thm:SConnWConn}
For all $1 \leq q \leq r$, $S$ is $(r-q)$-connected if and only if $W$ is $(r-q)$-connected.
\etheo
\setlength{\parindent}{0cm} \setlength{\parskip}{0cm}

Note that we will only need the implication $\Rarr$.
\bproof
This follows from \cite[Proposition~2.9]{RWW} and \cite[Theorem~2.10]{RWW} because $\textbf{V}_{\fB}$ satisfies conditions $\textbf{LS1}$ and $\textbf{LS2}$ from \cite[Definition~2.5]{RWW} by Lemma~\ref{lem:LS} and $\textbf{Q}_{\fB}$ satisfies conditions $\textbf{H1}$ and $\textbf{H2}$ from \cite[Definition~1.3]{RWW}. The latter follows from parts (c) and (d) of \cite[Theorem~1.10]{RWW}, because condition $\textbf{C}$ in \cite[Definition~1.9]{RWW} is satisfied in $\textbf{V}_{\fB}$ by Lemma~\ref{lem:i,ii,C} and the condition in part (d) is easy to check in our setting.
\eproof
\setlength{\parindent}{0cm} \setlength{\parskip}{0.5cm}

To proceed, we need the following observation. Recall that $M(G)$ denotes the set of all non-zero Radon measures $\mu$ on $\Gn$ which are invariant, and that an ample groupoid $G$ is said to have comparison if for all non-empty compact open sets $U, V \subseteq \Gn$ with $\mu(U) < \mu(V)$ for all $\mu \in M(G)$, there exists a compact open bisection $\sigma \subseteq G$ with $\rms(\sigma) = U$ and $\rmr(\sigma) \subseteq V$ (see \S~\ref{ss:GPD}). 
\blemma
\label{lem:ExtendBisec}
Assume that $G$ is minimal and has comparison. Let $Y \subseteq \Gn$ be a compact open subspace. Suppose that $A \subsetneq Y$ is given, together with a compact open bisection $\sigma$ with $\rms(\sigma) = A$ and $B \defeq \rmr(\sigma) \subsetneq Y$. If $2 \mu(A) < \mu(Y)$ for all $\mu \in M(G)$, then there exists a compact open bisection $\tau$ with $\rmr(\tau) = Y = \rms(\tau)$ and $\tau A = \sigma$.
\elemma
\setlength{\parindent}{0cm} \setlength{\parskip}{0cm}

Here $\tau A$ is the restriction of $\tau$ to $A$, i.e.,
$$
 \tau A = \menge{g \in \tau}{\rms(g) \in A},
$$
which is also the product of $\tau$ and $A$ as bisections.
\bproof
$2 \mu(A) < \mu(Y)$ for all $\mu \in M(G)$ implies that $\mu(B) = \mu(A) < \mu(Y \setminus A)$ for all $\mu \in M(G)$. Set $C \defeq (Y \setminus A) \cap B$. Then $\mu(B) < \mu(Y \setminus A)$ implies $\mu(B \setminus C) = \mu(B) - \mu(C) < \mu(Y \setminus A) - \mu(C) = \mu((Y \setminus A) \setminus C)$ for all $\mu \in M(G)$. As $G$ has comparison, there exists a compact open bisection $\ti{\rho}$ with $\rms(\ti{\rho}) = B \setminus C$ and $\rmr(\ti{\rho}) \subseteq (Y \setminus A) \setminus C$. As $(B \setminus C) \cap ((Y \setminus A) \setminus C) = \emptyset$, there exists a compact open bisection $\bar{\rho}$ with $\rmr(\bar{\rho}) = (B \setminus C) \amalg ((Y \setminus A) \setminus C) = \rms(\bar{\rho})$ extending $\ti{\rho}$, i.e., $\bar{\rho} (B \setminus C) = \ti{\rho}$. Now 
$$
 \rho \defeq \bar{\rho} \amalg (Y \setminus (B \cup (Y \setminus A))) \amalg C
$$
defines a compact open bisection with $\rmr(\rho) = Y = \rms(\rho)$, with the property that
$$
 \rho.B = \rho.(B \setminus C) \amalg \rho.C \subseteq ((Y \setminus A) \setminus C) \amalg C = Y \setminus A,
$$
i.e., $(\rho.B) \cap A = \emptyset$. Here $\rho.B = \menge{\rmr(g)}{g \in \rho, \, \rms(g) \in B}$ (see \S~\ref{s:smallpermcat} for the definition of the $G$-action on $\Gn$).
\setlength{\parindent}{0cm} \setlength{\parskip}{0.5cm}

Now consider the decomposition
$$
 Y \setminus A = (Y \setminus (A \amalg \rho.B)) \amalg \rho.B
$$
and define 
$$
 \tau' \defeq \rho^{-1} \rukl{ (Y \setminus (A \amalg \rho.B)) \amalg \sigma^{-1} (\rho B)^{-1} }.
$$
This is well-defined because $(\sigma^{-1} (\rho B)^{-1}).(\rho.B) = \sigma^{-1}.B = A$. Then, by construction $\tau'.(Y \setminus A) = \rho^{-1}.(Y \setminus \rho.B) = Y \setminus B$. So we obtain the compact open bisection with the desired properties by setting 
$\tau \defeq \sigma \amalg (\tau' (Y \setminus A))$.
\eproof
\setlength{\parindent}{0cm} \setlength{\parskip}{0.5cm}

\bremark
In the special case that $G$ is purely infinite, the condition $2 \mu(A) < \mu(Y)$ for all $\mu \in M(G)$ in Lemma~\ref{lem:ExtendBisec} is empty. Indeed, if $G$ is minimal and purely infinite, we can always find a compact open bisection with $\rmr(\rho) = Y = \rms(\rho)$ such that $(\rho.B) \cap A = \emptyset$.
\eremark

\bremark
It is straightforward to see that if $G$ has comparison, then so does $\cR \times G$. In the following, we will frequently apply Lemma~\ref{lem:ExtendBisec} to the groupoid $\cR \times G$, using the observation that $\fB$ consists of compact open bisections in $\cR \times G$ (see Remark~\ref{rem:B=BisecInRxG}). 
\eremark

\bprop
\label{prop:DescVertBisec}
Assume that $a \neq \emptyset$, and that $G$ is minimal and has comparison.
\setlength{\parindent}{0cm} \setlength{\parskip}{0cm}

\begin{enumerate}
\item[(i)] If $r \geq 2$, then vertices of $S$ (which coincide with the vertices of $W$) are in bijection with compact open bisections $\sigma$ with $\rms(\sigma) = u$ and $\rmr(\sigma) \subsetneq a \oplus u^{\oplus r}$, via the map $\cV$ defined by $\cV(\sigma) \defeq [u',\tau]$, where $\tau \in \textbf{V}_{\fB}(a \oplus u^{\oplus r}, a \oplus u^{\oplus r})$ satisfies $\tau \epsilon_r = \sigma$. Here $\epsilon_r \subseteq \cR \times \Gn$ is the compact open bisection with $\rms(\epsilon_r) = u$ and $\rmr(\epsilon_r)$ is the $r$-th summand of $u$ in $a \oplus u^{\oplus r}$.
\item[(ii)] Given $[u',\tau] \in W_p$, $[u'',\sigma] \in \textbf{Q}_{\fB}(a \oplus u^{\oplus r},u)$ is a vertex of $[u',\tau]$ if and only if there exists $i$ such that $\tau \epsilon_i = \sigma \epsilon_u$, where $\epsilon_i \subseteq \cR \times \Gn$ is the compact open bisection with $\rms(\epsilon_i) = u$ and $\rmr(\epsilon_i)$ is the $i$-th summand of $u$ in $u' \oplus u^{\oplus (p+1)}$, and $\epsilon_u \subseteq \cR \times \Gn$ is the compact open bisection with $\rms(\epsilon_u) = u$ and $\rmr(\epsilon_u)$ is the summand $u$ in $u'' \oplus u$.
\item[(iii)] Assume that $G$ is purely infinite and $0 \leq p \leq r-1$ is arbitrary or that $M(G) \neq \emptyset$ and $p$ satisfies $2(p+1) \leq r$. Given $p+1$ bisections $\sigma_i$ as in (i), $\cV(\sigma_i)$ form a $p$-simplex of $S$ if and only if $\rmr(\sigma_i)$ are pairwise disjoint and $\coprod_i \rmr(\sigma_i) \subsetneq a \oplus u^{\oplus r}$.
\end{enumerate}
\eprop
\setlength{\parindent}{0cm} \setlength{\parskip}{0cm}
Note that proper inclusions $\subsetneq$ are needed in (i) and (iii) to leave space for other vertices in higher-dimensional simplices.

\bproof
(i) The map $\cV$ is well-defined because given $\sigma$ as in (i), Lemma~\ref{lem:ExtendBisec} implies that there exists $\tau \in \textbf{V}_{\fB}(a \oplus u^{\oplus r}, a \oplus u^{\oplus r})$ with $\tau \epsilon_r = \sigma$. Here we use the assumptions $r \geq 2$ and $a \neq \emptyset$. It is straightforward to see that $\cV$ is bijective.
\setlength{\parindent}{0cm} \setlength{\parskip}{0.5cm}

(ii) is clear by definition.

Let us prove (iii). By assumption, there exists a compact open bisection $\tau'$ with $\rmr(\tau') \subseteq a \oplus u^{\oplus r}$ and $\rms(\tau') = u^{\oplus (p+1)}$ such that $\tau' \epsilon_i = \sigma_i$. By Lemma~\ref{lem:ExtendBisec}, there exists a compact open bisection $\tau$ with $\rmr(\tau) = a \oplus u^{\oplus r} = \rms(\tau)$ such that $\tau \epsilon_{u^{\oplus (p+1)}} = \tau'$. Here $\epsilon_{u^{\oplus (p+1)}} \subseteq \cR \times \Gn$ is the compact open bisection with $\rms(\epsilon_{u^{\oplus (p+1)}}) = u^{\oplus (p+1)}$ and $\rmr(\epsilon_{u^{\oplus (p+1)}})$ is the summand $u^{\oplus (p+1)}$ in $a \oplus u^{\oplus r} = a \oplus u^{\oplus (r-p-1)} \oplus u^{\oplus (p+1)}$. Now (ii) implies our claim.
\eproof
\setlength{\parindent}{0cm} \setlength{\parskip}{0.5cm}

Let us summarize the description of $S$ we obtain based on Proposition~\ref{prop:DescVertBisec}.
\bcor
\label{cor:DescS}
Assume that $a \neq \emptyset$, and that $G$ is minimal and has comparison. 

If $G$ is purely infinite, i.e., when $M(G) = \emptyset$, then $S$ can be described as the simplicial complex with vertices given by compact open bisections $\sigma$ with $\rms(\sigma) = u$ and $\rmr(\sigma) \subsetneq a \oplus u^{\oplus r}$, with the property that for all $0 \leq p \leq r-1$, $p+1$ vertices $\sigma_i$ form a $p$-simplex if and only if $\rmr(\sigma_i)$ are pairwise disjoint and $\coprod_i \rmr(\sigma_i) \subsetneq a \oplus u^{\oplus r}$.

If $M(G) \neq \emptyset$, then $S$ is a simplicial complex whose vertices are given by compact open bisections $\sigma$ with $\rms(\sigma) = u$ and $\rmr(\sigma) \subsetneq a \oplus u^{\oplus r}$, with the property that for all $p$ with $2(p+1) \leq r$, $p+1$ vertices $\sigma_i$ form a $p$-simplex if and only if $\rmr(\sigma_i)$ are pairwise disjoint and $\coprod_i \rmr(\sigma_i) \subsetneq a \oplus u^{\oplus r}$.
\ecor

\subsection{Connectivity of simplicial complexes}

\btheo
\label{thm:PInfSConn}
Assume that $a \neq \emptyset$, and that $G$ is minimal and purely infinite. Then $S$ is $(r-1)$-connected.
\etheo
\setlength{\parindent}{0cm} \setlength{\parskip}{0cm}

\bproof
We work with the description of $S$ from Corollary~\ref{cor:DescS}. Let $f: \: S^k \to S$ be a continuous map from the $k$-dimensional sphere $S^k$ to $S$, where $k \leq r-2$. Find a triangulation of $S^k$ such that $f$ is simplicial. Let $\nu_i$ be the number of $i$-simplices in the triangulation of $S^k$, and set $\nu \defeq \sum_{i=0}^k \nu_i$. Given a simplex $\bm{\sigma} = \gekl{\sigma_1, \dotsc, \sigma_i}$ of $S$, define $\rmr(\bm{\sigma}) \defeq \coprod_i \rmr(\sigma_i)$. Let $\cP$ be a partition of $a \oplus u^{\oplus r}$ into compact open subspaces refining all compact open subspaces of the form $\rmr(f(v))$, where $v$ is a vertex of $S^k$, such that for all simplices $\Delta$ of $S^k$, $\rmr(f(\Delta))^c$ contains at least $\nu+2$ elements of $\cP$.
\setlength{\parindent}{0cm} \setlength{\parskip}{0.5cm}

Our goal is to show that $f$ is homotopic to another simplicial map whose image only contains vertices $\sigma \in S$ which are small, in the sense that there exists $V \in \cP$ such that $\rmr(\sigma) \subseteq V$. In the process, we will re-triangulate $S^k$ such that there are always at most $\nu$ vertices. We will modify $f$ such that we keep the property that for all simplices $\Delta$ of $S^k$, $\rmr(f(\Delta))^c$ contains at least $\nu+2$ elements of $\cP$. In the following, we call a simplex $\Delta$ of $S^k$ bad if all vertices of $f(\Delta)$ are not small. In other words, a simplex $\Delta$ is not bad if and only if at least one vertex of $f(\Delta)$ is small. Let us now go through the bad simplices, removing them one by one, proceeding inductively on $\dim \Delta$.

We start with the case $\dim \Delta = k$. Since $\rmr(f(\Delta))^c$ contains at least $\nu+2$ elements of $\cP$, we can choose $V \in \cP$ with $V \subsetneq \rmr(f(\Delta))^c$. Furthermore, choose a compact open bisection $\tau$ with $\rmr(\tau) \subseteq V$ and $\rms(\tau) = u$. Then $f(\Delta) \cup \gekl{\tau}$ is a simplex in $S$. Add a vertex $a$ in the centre of $\Delta$, replace $\Delta$ by $\partial \Delta * a$ and replace $f$ by $f \vert_{\partial \Delta} * (a \ma \tau)$. This new map is homotopic to $f$ through the simplex $f(\Delta) \cup \gekl{\tau}$ because the two maps are contiguous (see for instance \cite[Chapter~3, \S~5]{Spa}).

We have added the vertex $a$, which is mapped to $\tau$ and hence is small. In this way we removed $\Delta$. Hence the number of bad simplices decreased. Moreover, we only increased the number of vertices by at most $1$. In addition, we still keep the property that for all simplices $\ti{\Delta}$ of $S^k$, $\rmr(f(\ti{\Delta}))^c$ contains at least $\nu+2$ elements of $\cP$. This is clear if $a \notin \ti{\Delta}$. If $a \in \ti{\Delta}$, then $\ti{\Delta} \setminus \gekl{a} \subseteq \Delta$. We replaced a vertex $v$ of $\Delta$ by $a$. $v$ must have been an original vertex, and hence $\rmr(f(v))$ covers at least one element of $\cP$. This is the reason why we keep the property, as claimed.

Now let $\Delta$ be a bad simplex of maximal dimension $\dim \Delta = j < k$. Then, by maximality, $f({\rm Link} \, \Delta)$ only contains small vertices. $\rmr(f(\Delta))^c$ contains at least $\nu+2$ elements of $\cP$, say $\gekl{V_i}$. Choose compact open bisections with $\rmr(\tau_i) \subseteq V_i$ and $\rms(\tau_i) = u$. Then, for all $i$, $f(\Delta) \cup \gekl{\tau_i}$ is a simplex in $S$. By the pigeonhole principle there exist $\omega, \omega' \in \gekl{\tau_i}$ such that no vertex of ${\rm Link} \, \Delta$ is mapped to $\omega$ or $\omega'$. So for all simplices $\Delta'$ of ${\rm Link} \, \Delta$, $f(\Delta') \cup f(\Delta) \cup \gekl{\omega}$ is a simplex in $S$. Add a vertex $a$ in the centre of $\Delta$, replace $\Delta$ by $a * \partial \Delta$, and replace $f$ in ${\rm Star} \, \Delta = ({\rm Link} \, \Delta) * \Delta \cong S^{k-j-1} * D^j$ by $(f \vert_{{\rm Link} \, \Delta}) * (a \ma \omega) * (f \vert_{\partial \Delta})$. We obtain a new map which is homotopic to $f$ via $(f \vert_{{\rm Link} \, \Delta}) * (a \ma \omega) * (f \vert_{\Delta})$ on $({\rm Link} \, \Delta) * a * \Delta \cong S^{k-j-1} * D^0 * D^j$ because they are contiguous, as above.

The number of vertices increased by at most $1$ (it only increases if $\partial \Delta \neq \emptyset$). Moreover, $a$ is mapped to $\omega$ and hence is small. Therefore, we have not added any new bad simplices. As we removed $\Delta$, the number of bad simplices decreased. In addition, we still keep the property that for all simplices $\ti{\Delta}$ of $S^k$, $\rmr(f(\ti{\Delta}))^c$ contains at least $\nu+2$ elements of $\cP$. This is clear if $a \notin \ti{\Delta}$. If $a \in \ti{\Delta}$, then $\ti{\Delta} \setminus \gekl{a} \subseteq \Delta * (\ti{\Delta} \cap {\rm Link} \, \Delta)$, which is an original simplex, with at least one original vertex not in $\ti{\Delta}$. The range of the image under $f$ of this original vertex covers at least one element of $\cP$. This is the reason why we keep the property, as claimed.

After this process, we obtain a map, again denoted by $f$, together with a triangulation of $S^k$ with at most $\nu$ vertices such that all vertices in the image of $f$ are small. Our new triangulation has at most $\nu$ vertices (where $\nu$ is the number of simplices in the original triangulation) because the number of vertices increases by at most $1$ for each bad simplex $\Delta$ in the original triangulation with $\dim \Delta > 0$. As $\# \cP \geq \nu+2$, there exist $V, V' \in \cP$ such that the image of every vertex of $S^k$ is disjoint from $V$ and $V'$. Choose compact open bisections $\tau, \tau'$ with $\rmr(\tau) \subseteq V$, $\rmr(\tau') \subseteq V'$ and $\rms(\tau) = u = \rms(\tau')$. It follows that for every simplex $\Delta$ of $S^k$, $f(\Delta) \cup \gekl{\tau}$ is a simplex in $S$.

It follows that $f$ is contiguous to the simplicial map sending all vertices to $\tau$. Hence $f$ is homotopic to a constant map, as desired.
\eproof
\setlength{\parindent}{0cm} \setlength{\parskip}{0.5cm}

\btheo
\label{thm:FiniteSConn}
Assume that $a \neq \emptyset$, and that $G$ is minimal and has comparison with $M(G) \neq \emptyset$. Then $S$ is $(l-1)$-connected if $\frac{(l+2)(l+3)}{2} < r$.
\etheo
\setlength{\parindent}{0cm} \setlength{\parskip}{0cm}

\bproof
Given $k \leq l$, $\frac{(l+2)(l+3)}{2} < r$ implies $2(p+1) < r$ for all $p$ with $p \leq k$. Hence we may and will work with the description of the $p$-simplices of $S$ for $2(p+1) < r$ from Corollary~\ref{cor:DescS}. 
\setlength{\parindent}{0cm} \setlength{\parskip}{0.5cm}

Let $f: \: S^{k-1} \to S$ be a continuous map from the $(k-1)$-dimensional sphere $S^{k-1}$ to $S$, where $k \leq l$. Find a triangulation of $S^{k-1}$ such that $f$ is simplicial.

Set $u_j \defeq u^{\oplus (j+3)}$. By our assumption that $\frac{(l+2)(l+3)}{2} < r$, we have
$$
 u_{k-1} \oplus u_{k-2} \oplus \dotso \oplus u_0 \oplus u_{-1} \subseteq a \oplus u^{\oplus r}.
$$
Set $u^{(i)} \defeq u_{k-1} \oplus u_{k-2} \oplus \dotso \oplus u_i$. We call an $i$-simplex of $S^{k-1}$ bad if all its vertices $v$ satisfy $\rmr(f(v)) \not\subseteq u^{(i)}$. Our goal is to show that $f$ is homotopic to another simplicial map $\ti{f}$ such that $\rmr(\ti{f}(v)) \subseteq u^{(0)}$ for all $v$ in a possibly new triangulation of $S^{k-1}$. As $G$ has comparison, there exists a compact open bisection $\tau$ with $\rmr(\tau) \subseteq u_{-1}$ and $\rms(\tau) = u$. Hence $\ti{f}(\Delta) \cup \gekl{\tau}$ is a simplex in $S$ for all simplices $\Delta$ of $S^{k-1}$. So the same argument as in the final step of the proof of Theorem~\ref{thm:PInfSConn} shows that $\ti{f}$ is homotopic to a constant map.

Let us now explain the procedure to remove bad simplices. Again, we start with bad simplices $\Delta$ with $\dim \Delta = k-1$. By comparison, there exists a compact open bisection $\tau$ with $\rms(\tau) = u$ and $\rmr(\tau) \subsetneq \rmr(f(\Delta))^c \cap u_{k-1}$. Then $f(\Delta) \cup \gekl{\tau}$ is a simplex in $S$. Add a vertex $a$ in the centre of $\Delta$, replace $\Delta$ by $\partial \Delta * a$ and replace $f$ by $f \vert_{\partial \Delta} * (a \ma \tau)$. This new map is homotopic to $f$ through the simplex $f(\Delta) \cup \gekl{\tau}$ because the two maps are contiguous.

In this way, we decreased the number of bad $(k-1)$-simplices.

Now assume that $\Delta$ is a bad simplex of maximal dimension $j = \dim \Delta < k-1$. By maximality, all vertices in $f({\rm Link} \, \Delta)$ have range in $u^{(j+1)}$. Otherwise, if there exists a vertex $v \in {\rm Link} \, \Delta$ with $\rmr(f(v)) \not\subseteq u^{(j+1)}$, then $\Delta \cup \gekl{v}$ would be a bad simplex (here we use that $u^{(j+1)} \subseteq u^{(j)}$) of dimension $j+1 = \dim \Delta + 1$, i.e., of higher dimension than $\Delta$. By comparison, there exists a compact open bisection $\tau$ with $\rms(\tau) = u$ and $\rmr(\tau) \subsetneq \rmr(f(\Delta))^c \cap u_j$. In particular, $\rmr(\tau) \cap u^{(j+1)} = \emptyset$. Hence, for all simplices $\Delta'$ of ${\rm Link} \, \Delta$, $f(\Delta') \cup f(\Delta) \cup \gekl{\tau}$ is a simplex in $S$. Add a vertex $a$ in the centre of $\Delta$, replace $\Delta$ by $a * \partial \Delta$, and replace $f$ in ${\rm Star} \, \Delta = ({\rm Link} \, \Delta) * \Delta \cong S^{k-j-1} * D^j$ by $(f \vert_{{\rm Link} \, \Delta}) * (a \ma \tau) * (f \vert_{\partial \Delta})$. We obtain a new map which is homotopic to $f$ via $(f \vert_{{\rm Link} \, \Delta}) * (a \ma \tau) * (f \vert_{\Delta})$ on $({\rm Link} \, \Delta) * a * \Delta \cong S^{k-j-1} * D^0 * D^j$. Again, we succeeded in decreasing the number of bad simplices. Indeed, after this modification, a simplex $\ti{\Delta}$ containing $a$ is not a bad simplex of dimension $\leq j$ because $\rmr(f(a)) = \rmr(\tau) \subseteq u_j \subseteq u^{(j)}$. If $\dim \ti{\Delta} > j$, then $({\rm Link} \, \Delta) \cap \ti{\Delta} \neq \emptyset$, so that $\ti{\Delta}$ is not bad. And if $a \notin \ti{\Delta}$, then $\ti{\Delta}$ is a simplex in the original triangulation but with $\ti{\Delta} \neq \Delta$.
\eproof
\setlength{\parindent}{0cm} \setlength{\parskip}{0.5cm}

\subsection{Homological stability and Morita invariance}
\label{ss:Morita}

As before, let $a, u \in \obj \textbf{Q}_{\fB}$ with $u \neq \emptyset$. Now let $W^r$ be the semi-simplicial set as defined above given by $W^r_p \defeq \textbf{Q}_{\fB}(a \oplus u^{\oplus r}, u^{\oplus (p+1)})$ for $0 \leq p \leq r-1$. We add the superscript $r$ to keep track of the number of summands of $u$ in $a \oplus u^{\oplus r}$ because we now want to vary $r$. 

First, we establish the following consequence of Theorems~\ref{thm:SConnWConn}, \ref{thm:PInfSConn}, \ref{thm:FiniteSConn} as well as \cite[Theorems~3.1 and 3.4]{RWW}.
\bcor
\label{cor:i(r)grows}
Assume that $a \neq \emptyset$, and that $G$ is minimal and has comparison.

For all $r$ there exists $i(r)$, which grows monotonously with $r$ such that $i(r) \to \infty$ if $r \to \infty$, with the property that the canonical map $\textbf{V}_{\fB}( a \oplus u^{\oplus r},  a \oplus u^{\oplus r}) \to \textbf{V}_{\fB}( a \oplus u^{\oplus (r+1)},  a \oplus u^{\oplus (r+1)})$ induces isomorphisms
$$
 H_i(\textbf{V}_{\fB}( a \oplus u^{\oplus r},  a \oplus u^{\oplus r}),\sfC) \to H_i(\textbf{V}_{\fB}( a \oplus u^{\oplus (r+1)},  a \oplus u^{\oplus (r+1)}),\sfC)
$$
for all $i \leq i(r)$.
\ecor
\setlength{\parindent}{0cm} \setlength{\parskip}{0cm}

\bproof
First of all, the proof of Theorem~3.1 and Theorem~3.4 in \cite{RWW} yields the following statement: 
\setlength{\parindent}{0cm} \setlength{\parskip}{0.5cm}

Suppose that $N$ is an integer such that $W^r$ is $(\frac{r-2}{k})$-connected for all $r$ with $r+1 \leq N$. Then, for all $r$ with $r+1 \leq N$, the canonical map
$$
 H_i(\textbf{V}_{\fB}( a \oplus u^{\oplus r},  a \oplus u^{\oplus r}),\sfC) \to H_i(\textbf{V}_{\fB}( a \oplus u^{\oplus (r+1)},  a \oplus u^{\oplus (r+1)}),\sfC)
$$
is an epimorphism for all $i$ with $i \leq \frac{r}{k}$ and an isomorphism for all $i$ with $i \leq \frac{r-1}{k}$. 

Here we have introduced an upper bound $N$ as the upper bound for $l$ in Theorem~\ref{thm:FiniteSConn} does not grow linearly with $r$.

Indeed, examining the proof of Theorem~3.1 and Theorem~3.4 in \cite{RWW}, we see that, in the notation of the proofs of \cite[Theorems~3.1 and 3.4]{RWW}, the proofs of (a), (b), (c) and (d) work for $r$ fixed. Moreover, the proofs of ($\textbf{E}_I 1$) and ($\textbf{E}_I 2$) work for all $r$ with $r+1 \leq N$. Similarly, the proofs of ($\textbf{I}_I 1$) and ($\textbf{I}_I 2$) work for all $r$ with $r+1 \leq N$. The proof of ($\textbf{I}_I 3$) works anyway.

It is now straightforward to derive the desired statement using Theorems~\ref{thm:SConnWConn}, \ref{thm:PInfSConn} and \ref{thm:FiniteSConn}.
\eproof
\setlength{\parindent}{0cm} \setlength{\parskip}{0.5cm}

\btheo
\label{thm:BvvBvzvz}
Suppose that $G$ is an ample groupoid, with locally compact Hausdorff unit space, and assume that $G$ is minimal and has comparison. Moreover assume that $\Gn$ has no isolated points. Then for all $v, z \in \obj \fB$ with $v \neq \emptyset$, the canonical map $\fB(v,v) \to \fB(v \oplus z,v \oplus z)$ induces an $H_*(\sqcup,\sfC)$-isomorphism, i.e., an isomorphism $H_*(\fB(v,v),\sfC) \cong H_*(\fB(v \oplus z,v \oplus z),\sfC)$ for all $* \geq 0$.
\etheo
\setlength{\parindent}{0cm} \setlength{\parskip}{0cm}

\bproof
In the following, we prove the statement for the case where $v = V$ for some non-empty compact open subspace $V \subseteq \Gn$. The general case is similar. 
\setlength{\parindent}{0cm} \setlength{\parskip}{0.5cm}

First of all, fix an index $* \geq 0$. 

Since $G$ is minimal, $\Gn$ is totally disconnected, locally compact, Hausdorff, without isolated points, given an arbitrary (big) natural number $r$, there exist a non-empty compact open subspace $u \subseteq v$ and $r$ compact open bisections $\sigma_i \subseteq G$ with $\rms(\sigma_i) = u$ such that $\rmr(\sigma_i)$ are pairwise disjoint with $\rmr(\sigma_i) \subseteq v$ and $v \setminus \coprod_i \rmr(\sigma_i) \neq \emptyset$. Set $a \defeq v \setminus \coprod_i \rmr(\sigma_i)$. Let $\sigma$ be a compact open bisection with $\rms(\sigma) = a \oplus u^{\oplus r}$ and $\rmr(\sigma) = v$. Then conjugation with $\sigma$, i.e., $\sigma^{-1} \sqcup \sigma$ yields an isomorphism
$$
 \fB(v,v) \cong \textbf{V}_{\fB}(a \oplus u^{\oplus r},a \oplus u^{\oplus r}).
$$
By Corollary~\ref{cor:i(r)grows}, for sufficiently big $r$ (more precisely, for all $r$ such that $i(r) \geq *$), we have that
$$
 \sqcup \oplus u: \: \fB(v,v) \to \fB(v \oplus u,v \oplus u)
$$
induces an $H_*(\sqcup,\sfC)$-isomorphism. Hence, for all $r$ sufficiently big and all $m \geq 0$, 
$$
 \sqcup \oplus u^{\oplus m}: \: \fB(v,v) \to \fB(v \oplus u^{\oplus m},v \oplus u^{\oplus m})
$$
induces an $H_*(\sqcup,\sfC)$-isomorphism. Note that in these two statements, $u$ depends on $r$.

Since $G$ is minimal, there exists $m$ large enough so that there exists a compact open bisection $\tau$ with $\rms(\tau) = z$ and $\rmr(\tau) \subseteq u^{\oplus m}$. Now the composite
\begin{equation}
\label{e:emb-conj-emb}
 \xymatrix{
 \fB(v,v) \ar@{^{(}->}[r] & \fB(v \oplus z,v \oplus z) \ar[rr]^{\hspace{-0.5cm} v \amalg (\tau^{-1} \sqcup \tau)} & & \fB(v \oplus \rmr(\tau),v \oplus \rmr(\tau)) \ar@{^{(}->}[r] & \fB(v \oplus u^{\oplus m},v \oplus u^{\oplus m})
 }
\end{equation}
is the canonical embedding $\sqcup \oplus u^{\oplus m}$. 

It follows that $\fB(v,v) \into \fB(v \oplus z,v \oplus z)$ induces an injective map in $H_*(\sqcup,\sfC)$. But this holds for arbitrary non-empty $v$ and $z$. Hence also $\fB(v \oplus \rmr(\tau),v \oplus \rmr(\tau)) \into \fB(v \oplus u^{\oplus m},v \oplus u^{\oplus m})$ induces an injective map in $H_*(\sqcup,\sfC)$. But because $\sqcup \oplus u^{\oplus m}$ induces an $H_*(\sqcup,\sfC)$-isomorphism and coincides with the composition in \eqref{e:emb-conj-emb}, the map in $H_*(\sqcup,\sfC)$ induced by $\fB(v \oplus \rmr(\tau),v \oplus \rmr(\tau)) \into \fB(v \oplus u^{\oplus m},v \oplus u^{\oplus m})$ is also surjective. Hence $\fB(v \oplus \rmr(\tau),v \oplus \rmr(\tau)) \into \fB(v \oplus u^{\oplus m},v \oplus u^{\oplus m})$ induces an $H_*(\sqcup,\sfC)$-isomorphism. Therefore, since $v \amalg (\tau^{-1} \sqcup \tau)$ is an isomorphism, it follows that $\fB(v,v) \into \fB(v \oplus z,v \oplus z)$ induces an $H_*(\sqcup,\sfC)$-isomorphism, as desired.
\eproof
\setlength{\parindent}{0cm} \setlength{\parskip}{0.5cm}

The following is now an immediate consequence because of Remark~\ref{rem:Buu=Fu} (apply Theorem~\ref{thm:BvvBvzvz} to $v = U$ and $z = V \setminus U$).
\btheo
\label{thm:GUUGVV}
Suppose that $G$ is an ample groupoid, with locally compact Hausdorff unit space, and assume that $G$ is minimal and has comparison. Moreover assume that $\Gn$ has no isolated points. Then for all non-empty compact open subspaces $U \subseteq V$ of $\Gn$, the canonical map $\bmF(G_U^U) \to \bmF(G_V^V)$ induces an $H_*(\sqcup,\sfC)$-isomorphism for all abelian groups $\sfC$ and all $* \geq 0$.
\etheo

We obtain the following consequence.
\bcor
\label{cor:Morita}
Suppose that $G$ is an ample groupoid, with locally compact Hausdorff unit space, and assume that $G$ is minimal and has comparison. Moreover assume that $\Gn$ has no isolated points. Let $H$ be an ample groupoid, with locally compact Hausdorff unit space, which is equivalent to $G$. Then $H_*(\bmF(G),\sfC) \cong H_*(\bmF(H),\sfC)$ for all abelian groups $\sfC$ and all $* \geq 0$.
\ecor
\setlength{\parindent}{0cm} \setlength{\parskip}{0cm}

\bproof
Given a $(G,H)$-equivalence, let $L$ be the corresponding linking groupoid as in \cite[\S~4.1]{CS}. Then $L$ is an ample groupoid, with locally compact Hausdorff unit space without isolated points, which is minimal and has comparison, because $G$ has these properties. Similarly, $H$ is minimal, has comparison, and the unit space of $H$ has no isolated points. Hence Theorem~\ref{thm:GUUGVV} applies, and our claim follows because $G$ and $H$ are isomorphic to reductions of $L$.
\eproof
\setlength{\parindent}{0cm} \setlength{\parskip}{0.5cm}

\bremark
\label{rem:DMorita}
The same arguments as above, using \cite[Corollary~3.9]{RWW}, show that the analogues of Theorem~\ref{thm:GUUGVV} and Corollary~\ref{cor:Morita} are also true for the commutator subgroup in place of the topological full group, i.e., in the setting of Theorem~\ref{thm:GUUGVV} and Corollary~\ref{cor:Morita}, homology of commutator subgroups of topological full groups is also Morita invariant.
\eremark

\subsection{Identifying homology of infinite loop spaces with homology of topological full groups}

Let $\Omega^{\infty}_0 \Kz(\fB)$ denote the connected component of the base point of $\Omega^{\infty} \Kz(\fB)$ (see \S~\ref{ss:AlgKSmallPermCat}).

\btheo
\label{thm:GroupCompl}
Let $G$ be an ample groupoid, with locally compact Hausdorff unit space $\Gn$, and $\sfC$ an abelian group. Then there exists a map $B \bmF(\cR \times G) \to \Omega^{\infty}_0 \Kz(\fB)$ which induces an $H_*(\sqcup,\sfC)$-isomorphism, i.e.,
$$
 H_*(\bmF(\cR \times G),\sfC) \cong H_*(\Omega^{\infty}_0 \Kz(\fB),\sfC).
$$
\etheo
\setlength{\parindent}{0cm} \setlength{\parskip}{0cm}

\bproof
The proof is as in \cite{SW} (compare \cite[Theorem~5.4]{SW}, which is based on \cite{McDS,R-W}). Let $M = \abs{\fB}$ be the nerve or classifying space of $\fB$ as in \S~\ref{ss:AlgKSmallPermCat}. Let $M_{\infty}$ be the homotopy colimit of $M$ with respect to the maps given by $\sqcup \oplus v$, for $v \in \obj \fB$. The group completion theorem \cite{McDS} (see also \cite{R-W}) implies that there exists a map 
\begin{equation}
\label{MinfToOmega}
 M_{\infty} \to \Omega^{\infty} \Kz(\fB)
\end{equation}
which induces an $H_*(\sqcup,\sfC)$-isomorphism. The component of $\emptyset$ of $M_{\infty}$ can be identified with $B \bmF(\cR \times G)$ as
$$
 \bmF(\cR \times G) = \ilim_u \bmF((\cR \times G)_u^u) \cong \ilim_u \fB(u,u)
$$
by definition (see also Remark~\ref{rem:Buu=Fu}) and because $\fB$, being a groupoid, is equivalent to $\coprod_u \fB(u,u)$, where $u$ runs through a system of representatives for the components of $\fB$.
\setlength{\parindent}{0cm} \setlength{\parskip}{0.5cm}

Hence, restricting \eqref{MinfToOmega} to the component of $\emptyset$ of $M_{\infty}$ and the component of the base point $\Omega^{\infty}_0 \Kz(\fB)$ of $\Omega^{\infty} \Kz(\fB)$, our claim follows.
\eproof
\setlength{\parindent}{0cm} \setlength{\parskip}{0.5cm}

Let us now combine Theorem~\ref{thm:GroupCompl} with Morita invariance from \S~\ref{ss:Morita}.
\btheo
\label{thm:HFG=HOKB}
Let $G$ be an ample groupoid, with locally compact Hausdorff $\Gn$ without isolated points. Assume that $G$ is minimal and has comparison. Let $\sfC$ be an abelian group. Then $H_*(\bmF(\cR \times G),\sfC) \cong H_*(\bmF(G),\sfC)$. Hence $H_*(\bmF(G),\sfC) \cong H_*(\Omega^{\infty}_0 \Kz(\fB),\sfC)$.
\etheo
\setlength{\parindent}{0cm} \setlength{\parskip}{0cm}

\bproof
By definition, $\bmF(G) = \ilim_U \bmF(G_U^U)$, where the limit is taken over all compact open subspaces $U \subseteq \Gn$. Furthermore, $\bmF(\cR \times G) \cong \ilim_u \fB(u,u)$, using Remark~\ref{rem:Buu=Fu}. Now Theorem~\ref{thm:GUUGVV} implies that for all $\emptyset \neq U \subseteq \Gn$, $H_*(\bmF(G),\sfC) \cong H_*(\bmF(G_U^U),\sfC)$. Moreover, for all $z \in \obj \fB$, $\fB(U,U) \into \fB(U \oplus z,U \oplus z)$ induces an $H_*(\sqcup,\sfC)$-isomorphism, for all $* \geq 0$, by Theorem~\ref{thm:BvvBvzvz}. It follows that $H_*(\bmF(\cR \times G),\sfC) \cong H_*(\bmF(G_U^U),\sfC) \cong H_*(\bmF(G),\sfC)$, for all $* \geq 0$. Now apply Theorem~\ref{thm:GroupCompl}.
\eproof
\setlength{\parindent}{0cm} \setlength{\parskip}{0.5cm}

As explained in \S~\ref{ss:TFG} and \S~\ref{ss:ExTopFullGrp}, interesting examples of infinite simple groups do not arise directly from the construction of topological full groups; rather, they are given by alternating full groups $\bmA(G)$ for special ample groupoids $G$. For almost finite or purely infinite groupoids $G$ which are minimal, effective and Hausdorff, with unit space $\Gn$ homeomorphic to the Cantor space, the alternating full group coincides with the commutator subgroup $\bmD(G)$ of $\bmF(G)$ (see \cite{Mat15,Nek19}). The commutator subgroup $\bmD(G)$ (also called derived subgroup) is the subgroup of $\bmF(G)$ generated by commutators of the form $\sigma \tau \sigma^{-1} \tau^{-1}$ for $\sigma, \tau \in \bmF(G)$. Let us now explain how our approach allows us to study homology of $\bmD(G)$ as well.

In the following, let $\ti{\Omega}^{\infty}_0 \Kz(\fB_G)$ be the universal cover of $\Omega^{\infty}_0 \Kz(\fB_G)$.
\btheo
\label{thm:HDG=HtiOKB}
Let $G$ be an ample groupoid, with locally compact Hausdorff unit space. Then
$$
 H_*(\bmD(\cR \times G),\sfC) \cong H_*(\ti{\Omega}^{\infty}_0 \Kz(\fB_G),\sfC)
$$
for all abelian groups $\sfC$ and all $* \geq 0$.
\etheo
\setlength{\parindent}{0cm} \setlength{\parskip}{0cm}

\bproof
We have
\begin{align*}
 H_*(\bmD(\cR \times G),\sfC) &\cong H_*(\bmF(\cR \times G),\Zz[\bmF(\cR \times G)/\bmD(\cR \times G)] \otimes \sfC) \cong H_*(\bmF(\cR \times G),\Zz [H_1(\bmF(\cR \times G))] \otimes \sfC)\\
 &\cong H_*(\Omega^{\infty}_0 \Kz(\fB_G), \Zz [H_1(\Omega^{\infty}_0 \Kz(\fB_G))] \otimes \sfC).
\end{align*}
For the first isomorphism, we used Shapiro's Lemma (see for instance \cite[Chapter~III, Proposition~(6.2)]{Bro}). The second isomorphism is induced by the canonical identification $H_1(\bmF(\cR \times G)) \cong \bmF(\cR \times G)/\bmD(\cR \times G)$. The third isomorphism follows from the group completion theorem \cite{McDS} (see also \cite{R-W}) in the same way as Theorem~\ref{thm:GroupCompl} does (see also \cite[Remark~2.5]{R-W}).
\setlength{\parindent}{0cm} \setlength{\parskip}{0.5cm}

In addition, we have
$$
 H_*(\Omega^{\infty}_0 \Kz(\fB_G), \Zz [H_1(\Omega^{\infty}_0 \Kz(\fB_G))] \otimes \sfC) \cong H_*(\Omega^{\infty}_0 \Kz(\fB_G), \Zz [\pi_1(\Omega^{\infty}_0 \Kz(\fB_G))] \otimes \sfC) \cong H_*(\ti{\Omega}^{\infty}_0 \Kz(\fB_G), \sfC).
$$
The first isomorphism is induced by the Hurewicz isomorphism $\pi_1(\Omega^{\infty}_0 \Kz(\fB_G)) \cong H_1(\Omega^{\infty}_0 \Kz(\fB_G))$ (using that $\pi_1(\Omega^{\infty}_0 \Kz(\fB_G))$ is abelian because $\Omega^{\infty}_0 \Kz(\fB_G)$ is an infinite loop space). For the second isomorphism,  we refer the reader for instance to \cite[Example~3H.2]{Hat}. This proves our claim.
\eproof
\setlength{\parindent}{0cm} \setlength{\parskip}{0.5cm}

\bcor
\label{cor:HDG=HtiOKB}
Let $G$ be an ample groupoid, with locally compact Hausdorff unit space $\Gn$ without isolated points. Assume that $G$ is minimal and has comparison. Then
$$
 H_*(\bmD(G),\sfC) \cong H_*(\bmD(\cR \times G),\sfC) \cong H_*(\ti{\Omega}^{\infty}_0 \Kz(\fB_G),\sfC)
$$
for all abelian groups $\sfC$ and all $* \geq 0$.
\ecor
\setlength{\parindent}{0cm} \setlength{\parskip}{0cm}

\bproof
As noted in Remark~\ref{rem:DMorita}, using \cite[Corollary~3.9]{RWW}, the same argument as for Theorem~\ref{thm:HFG=HOKB} implies that
$$
 H_*(\bmD(\cR \times G), \sfC) \cong H_*(\bmD(G_U^U), \sfC) \cong H_*(\bmD(G), \sfC),
$$
for all non-empty compact open subspaces $U \subseteq \Gn$ (see also the argument for \cite[Corollary~6.7]{SW}, which is similar).
\eproof
\setlength{\parindent}{0cm} \setlength{\parskip}{0.5cm}

\section{Applications}
\label{s:App}

Let us now derive consequences from our main results. Let $G$ be an ample groupoid with locally compact Hausdorff unit space. 

First, note that $\Omega_0^{\infty} \Kz(\fB_G)$ inherits the structure of an infinite loop space from $\Omega^{\infty} \Kz(\fB_G)$. In particular, $\Omega^{\infty}_0 \Kz(\fB_G)$ is up to weak homotopy equivalence a homotopy-associative $H$-space. It can also be derived directly from definitions (see for instance \cite{Swi}) that $\Omega^{\infty}_0 \Kz(\fB_G)$ inherits the structure of a homotopy-associative $H$-space from $\Omega^{\infty} \Kz(\fB_G)$.


Moreover, observe that the $H$-space structure on $\Omega^{\infty}_0 \Kz(\fB_G)$ can be lifted to the universal cover $\ti{\Omega}^{\infty}_0 \Kz(\fB_G)$ (see for instance \cite[Section~3.C, Exercise~4]{Hat}). Moreover, because $\ti{\Omega}^{\infty}_0 \Kz(\fB_G)$ is the universal cover of $\Omega^{\infty}_0 \Kz(\fB_G)$, we have
\begin{equation}
\label{e:pitiOmega}
 \pi_*(\ti{\Omega}^{\infty}_0 \Kz(\fB_G)) \cong 
 \bfa
 \pi_*(\Omega^{\infty}_0 \Kz(\fB_G)) & \falls * \geq 2,\\
 \gekl{0} & \falls * = 0, 1.
 \efa
\end{equation}

\subsection{Rational homology}

Let us start with rational computations. We need the following notation. Given an ample groupoid $G$ with locally compact Hausdorff unit space, we denote by $H^{\rm odd}_*(G,\Qz)$ the groupoid homology of $G$ with rational coefficients in odd degree, i.e., 
$$
 H^{\rm odd}_*(G,\Qz) \defeq
 \bfa
 H_*(G,\Qz) & \falls * > 0 \text{ odd},\\
 \gekl{0} & \sonst.
 \efa
$$
Similarly, let $H^{\rm even}_*(G,\Qz)$ be the groupoid homology of $G$ with rational coefficients in (positive) even degree, i.e., 
$$
 H^{\rm even}_*(G,\Qz) \defeq
 \bfa
 H_*(G,\Qz) & \falls * > 0 \text{ even},\\
 \gekl{0} & \sonst.
 \efa
$$
\bcor
\label{cor:RatHom}
Let $G$ be an ample groupoid with locally compact Hausdorff unit space $\Gn$. Then
$$
 H_*(\bmF(\cR \times G),\Qz) \cong {\rm Ext}(H^{\rm odd}_*(G,\Qz)) \otimes {\rm Sym}(H^{\rm even}_*(G,\Qz))
$$
as graded vector spaces over $\Qz$.

Suppose, in addition, that $\Gn$ does not have isolated points, and that $G$ is minimal and has comparison. Then
$$
 H_*(\bmF(G),\Qz) \cong {\rm Ext}(H^{\rm odd}_*(G,\Qz)) \otimes {\rm Sym}(H^{\rm even}_*(G,\Qz))
$$
as graded vector spaces over $\Qz$.
\ecor
\setlength{\parindent}{0cm} \setlength{\parskip}{0cm}

\bproof
\cite[Chapter~II, Proposition~6.30~(iii)]{Sch}) implies that the Hurewicz maps induce isomorphisms
$$
 \pi_*(\Kz(\fB_G)) \otimes \Qz \cong H_*(\Kz(\fB_G),\Qz)
$$
for all $* \geq 0$. Hence
$$
 \pi_*(\Omega^{\infty}_0 \Kz(\fB_G)) \otimes \Qz \cong \pi_*(\Omega^{\infty} \Kz(\fB_G)) \otimes \Qz \cong \pi_*(\Kz(\fB_G)) \otimes \Qz \cong H_*(\Kz(\fB_G),\Qz) \cong H_*(G,\Qz)
$$
for all $*>0$, where we applied Theorem~\ref{thm:HKB=HG}. As explained above, $\Omega^{\infty}_0 \Kz(\fB_G)$ is a homotopy-associative $H$-space. Hence the theorem in the appendix of \cite{MM} implies that $H_*(\Omega^{\infty}_0 \Kz(\fB_G), \Qz)$ is isomorphic to the universal enveloping algebra $U(\pi_*(\Omega^{\infty}_0 \Kz(\fB_G)) \otimes \Qz)$ of $\pi_*(\Omega^{\infty}_0 \Kz(\fB_G)) \otimes \Qz$, which is a Lie algebra with respect to the Samelson product (see for instance the appendix of \cite{MM}), i.e., 
$$
 H_*(\Omega^{\infty}_0 \Kz(\fB_G), \Qz) \cong U(\pi_*(\Omega^{\infty}_0 \Kz(\fB_G)) \otimes \Qz).
$$
Now the Poincar{\'e}-Birkhoff-Witt Theorem (see for instance \cite[Theorem~21.1]{FHT}) implies that
$$
 U(\pi_*(\Omega^{\infty}_0 \Kz(\fB_G)) \otimes \Qz) \cong \Lambda(\pi_*(\Omega^{\infty}_0 \Kz(\fB_G)) \otimes \Qz),
$$
where $\Lambda$ stands for the free commutative graded algebra (see for instance \cite[Chapter~I, \S~3, Example~6]{FHT}). Note that $\Lambda(\pi_*(\Omega^{\infty}_0 \Kz(\fB_G)) \otimes \Qz)$ is constructed from the vector space $\pi_*(\Omega^{\infty}_0 \Kz(\fB_G)) \otimes \Qz$ and does not use the Lie algebra structure of $\pi_*(\Omega^{\infty}_0 \Kz(\fB_G)) \otimes \Qz$ anymore. Finally, it follows from \cite[Chapter~II, \S~12~(a)]{FHT} that 
$$
 \Lambda(\pi_*(\Omega^{\infty}_0 \Kz(\fB_G)) \otimes \Qz) \cong {\rm Ext}(H^{\rm odd}_*(G,\Qz)) \otimes {\rm Sym}(H^{\rm even}_*(G,\Qz)).
$$
All in all, using Theorem~\ref{thm:GroupCompl}, we obtain 
$$
 H_*(\bmF(\cR \times G),\Qz) \cong H_*(\Omega^{\infty}_0 \Kz(\fB_G), \Qz) \cong {\rm Ext}(H^{\rm odd}_*(G,\Qz)) \otimes {\rm Sym}(H^{\rm even}_*(G,\Qz)),
$$
and, if $\Gn$ does not have isolated points and $G$ is minimal and has comparison, Theorem~\ref{thm:HFG=HOKB} implies
$$
 H_*(\bmF(G),\Qz) \cong H_*(\Omega^{\infty}_0 \Kz(\fB_G), \Qz) \cong {\rm Ext}(H^{\rm odd}_*(G,\Qz)) \otimes {\rm Sym}(H^{\rm even}_*(G,\Qz)),
$$
as desired.
\eproof
\setlength{\parindent}{0cm} \setlength{\parskip}{0.5cm}

We record the following immediate consequence.
\bcor
\label{cor:RatAcyclic}
Let $G$ be an ample groupoid with locally compact Hausdorff unit space. $\bmF(\cR \times G)$ is rationally acyclic (i.e., $H_*(\bmF(\cR \times G),\Qz) \cong \gekl{0}$ for all $*>0$) if and only if $H_*(G,\Qz) \cong \gekl{0}$ for all $*>0$.

Suppose, in addition, that $\Gn$ does not have isolated points, and that $G$ is minimal and has comparison. Then $\bmF(G)$ is rationally acyclic if and only if $H_*(G,\Qz) \cong \gekl{0}$ for all $*>0$.
\ecor

Moreover, specializing to degree $1$, we obtain the following consequence of Corollary~\ref{cor:RatHom}.
\bcor
Let $G$ be an ample groupoid, with locally compact Hausdorff unit space $\Gn$. Then
$$
 H_1(\bmF(\cR \times G),\Qz) \cong H_1(G,\Qz).
$$
If, in addition, $\Gn$ has no isolated points, and $G$ is minimal and has comparison, then
$$
 H_1(\bmF(G),\Qz) \cong H_1(G,\Qz).
$$
\ecor

Furthermore, we obtain a formula for the Poincar{\'e} series 
$$
 P_{\bmF(G)}(t) \defeq \sum_{j=0}^{\infty} \dim H_j(\bmF(G),\Qz) t^j.
$$
For that purpose, we define, for $j>0$,
$$
 d_j \defeq \dim H_j(G,\Qz).
$$
\bcor
\label{cor:P_FG}
Let $G$ be an ample groupoid with locally compact Hausdorff unit space $\Gn$. Then
$$
 P_{\bmF(\cR \times G)}(t) = \prod_{j=1}^{\infty} (1 + (-1)^{j+1} t^j)^{(-1)^{j+1} d_j} = (1 + t^1)^{d_1} (1 + t^3)^{d_3} \dotsm (1 - t^2)^{-d_2} (1 - t^4)^{-d_4} \dotsm.
$$
If, in addition, $\Gn$ has no isolated points, and $G$ is minimal and has comparison, then
$$
 P_{\bmF(G)}(t) = \prod_{j=1}^{\infty} (1 + (-1)^{j+1} t^j)^{(-1)^{j+1} d_j} = (1 + t^1)^{d_1} (1 + t^3)^{d_3} \dotsm (1 - t^2)^{-d_2} (1 - t^4)^{-d_4} \dotsm.
$$
\ecor
\setlength{\parindent}{0cm} \setlength{\parskip}{0cm}

\bproof
This is an immediate consequence of Corollary~\ref{cor:RatHom}, together with known formulas for the Poincar{\'e} series for exterior and symmetric algebras (see for instance \cite[\S~5.7, equation~(5.15)]{Gre} and \cite[\S~9.11]{Gre}).
\eproof
\setlength{\parindent}{0cm} \setlength{\parskip}{0.5cm}

Let us now turn to rational computations for commutator subgroups. Given an ample groupoid $G$ with locally compact Hausdorff unit space, we denote by $H^{\rm odd}_{*>1}(G,\Qz)$ the groupoid homology of $G$ with rational coefficients in odd degree $>1$, i.e., 
$$
 H^{\rm odd}_{*>1}(G,\Qz) \defeq
 \bfa
 H_*(G,\Qz) & \falls * > 1 \text{ odd},\\
 \gekl{0} & \sonst.
 \efa
$$
Let $H^{\rm even}_*(G,\Qz)$ be defined as above.
\bcor
\label{cor:DRatHom}
Let $G$ be an ample groupoid with locally compact Hausdorff unit space $\Gn$. Then
$$
 H_*(\bmD(\cR \times G),\Qz) \cong {\rm Ext}(H^{\rm odd}_{*>1}(G,\Qz)) \otimes {\rm Sym}(H^{\rm even}_*(G,\Qz))
$$
as graded vector spaces over $\Qz$.

If, in addition, $\Gn$ has no isolated points, and $G$ is minimal and has comparison, then
$$
 H_*(\bmD(G),\Qz) \cong {\rm Ext}(H^{\rm odd}_{*>1}(G,\Qz)) \otimes {\rm Sym}(H^{\rm even}_*(G,\Qz))
$$
as graded vector spaces over $\Qz$.
\ecor
\bproof
The proof is similar to the one for Corollary~\ref{cor:RatHom}. We have
$$
 \pi_*(\ti{\Omega}^{\infty}_0 \Kz(\fB_G)) \otimes \Qz \cong \pi_*(\Omega^{\infty}_0 \Kz(\fB_G)) \otimes \Qz \cong \pi_*(\Omega^{\infty} \Kz(\fB_G)) \otimes \Qz \cong H_*(\Kz(\fB_G),\Qz) \cong H_*(G,\Qz)
$$
for all $*>1$, where we applied Theorem~\ref{thm:HKB=HG}. As explained above, $\ti{\Omega}^{\infty}_0 \Kz(\fB_G)$ is a homotopy-associative $H$-space. Hence the theorem in the appendix of \cite{MM} implies that $H_*(\ti{\Omega}^{\infty}_0 \Kz(\fB_G), \Qz)$ is isomorphic to the universal enveloping algebra $U(\pi_*(\ti{\Omega}^{\infty}_0 \Kz(\fB_G)) \otimes \Qz)$ of $\pi_*(\ti{\Omega}^{\infty}_0 \Kz(\fB_G)) \otimes \Qz$, which is a Lie algebra with respect to the Samelson product (see for instance the appendix of \cite{MM}). Now the Poincar{\'e}-Birkhoff-Witt Theorem (see for instance \cite[Theorem~21.1]{FHT}) and \cite[Chapter~II, \S~12~(a)]{FHT} imply that
$$
 U(\pi_*(\ti{\Omega}^{\infty}_0 \Kz(\fB_G)) \otimes \Qz) \cong \Lambda(\pi_*(\ti{\Omega}^{\infty}_0 \Kz(\fB_G)) \otimes \Qz) 
 \cong {\rm Ext}(H^{\rm odd}_{*>1}(G,\Qz)) \otimes {\rm Sym}(H^{\rm even}_*(G,\Qz)).
$$
All in all, using \eqref{e:pitiOmega} and Theorem~\ref{thm:HDG=HtiOKB}, we obtain
$$
 H_*(\bmD(\cR \times G),\Qz) \cong H_*(\ti{\Omega}^{\infty}_0 \Kz(\fB_G), \Qz) \cong {\rm Ext}(H^{\rm odd}_{*>1}(G,\Qz)) \otimes {\rm Sym}(H^{\rm even}_*(G,\Qz)),
$$
and, if $\Gn$ does not have isolated points and $G$ is minimal and has comparison, \eqref{e:pitiOmega} and Corollary~\ref{cor:HDG=HtiOKB} imply
$$
 H_*(\bmD(G),\Qz) \cong H_*(\ti{\Omega}^{\infty}_0 \Kz(\fB_G), \Qz) \cong {\rm Ext}(H^{\rm odd}_{*>1}(G,\Qz)) \otimes {\rm Sym}(H^{\rm even}_*(G,\Qz)),
$$
as desired.
\eproof

We record the following immediate consequence, as above.
\bcor
\label{cor:DRatAcyclic}
Let $G$ be an ample groupoid with locally compact Hausdorff unit space. $\bmD(\cR \times G)$ is rationally acyclic if and only if $H_*(G,\Qz) \cong \gekl{0}$ for all $*>1$.

Suppose, in addition, that $\Gn$ does not have isolated points, and that $G$ is minimal and has comparison. Then $\bmD(G)$ is rationally acyclic if and only if $H_*(G,\Qz) \cong \gekl{0}$ for all $*>1$.
\ecor

As before, we obtain the following formula for the Poincar{\'e} series 
$$
 P_{\bmD(G)}(t) \defeq \sum_{j=0}^{\infty} \dim H_j(\bmD(G),\Qz) t^j.
$$
Recall that we defined $d_j \defeq \dim H_j(G,\Qz)$ for $j>0$.
\bcor
\label{cor:P_DG}
Let $G$ be an ample groupoid with locally compact Hausdorff unit space $\Gn$. Then
$$
 P_{\bmD(\cR \times G)}(t) = \prod_{j=2}^{\infty} (1 + (-1)^{j+1} t^j)^{(-1)^{j+1} d_j} = (1 + t^3)^{d_3} (1 + t^5)^{d_5} \dotsm (1 - t^2)^{-d_2} (1 - t^4)^{-d_4} \dotsm.
$$
If, in addition, $\Gn$ has no isolated points, and $G$ is minimal and has comparison, then
$$
 P_{\bmD(G)}(t) = \prod_{j=2}^{\infty} (1 + (-1)^{j+1} t^j)^{(-1)^{j+1} d_j} = (1 + t^3)^{d_3} (1 + t^5)^{d_5} \dotsm (1 - t^2)^{-d_2} (1 - t^4)^{-d_4} \dotsm.
$$
\ecor

\subsection{Vanishing results}

In the following, we write $H_*(G) \defeq H_*(G,\Zz)$ and $H_*(\bmF(G)) \defeq H_*(\bmF(G),\Zz)$.
\bcor
\label{cor:HVanish}
Let $G$ be an ample groupoid with locally compact Hausdorff unit space $\Gn$. Given $k \in \Zz$ with $k > 0$, suppose that $H_*(G) \cong \gekl{0}$ for all $* < k$. Then $H_*(\bmF(\cR \times G)) \cong \gekl{0}$ for all $0 < * < k$ and $H_k(\bmF(\cR \times G)) \cong H_k(G)$. If, in addition, $\Gn$ has no isolated points and $G$ is minimal and has comparison, then $H_*(\bmF(G)) \cong \gekl{0}$ for all $0 < * < k$ and $H_k(\bmF(G)) \cong H_k(G)$.
\ecor
\setlength{\parindent}{0cm} \setlength{\parskip}{0cm}

\bproof
Assume that $H_*(G) \cong \gekl{0}$ for all $* < k$. Theorem~\ref{thm:HKB=HG} implies that $\ti{H}_*(\Kz(\fB_G),\Zz) \cong \gekl{0}$ for all $* < k$. Thus we obtain, for all $1 \leq * \leq k$, that
$$
 \pi_*(\Omega^{\infty}_0 \Kz(\fB_G)) \cong \pi_*(\Omega^{\infty} \Kz(\fB_G)) \cong \pi_*(\Kz(\fB_G)) \cong \ti{H}_*(\Kz(\fB_G),\Zz)
$$
by applying the stable Hurewicz Theorem (see for instance \cite[Chapter~II, Proposition~6.30~(i)]{Sch}) inductively. Now the usual Hurewicz Theorem for spaces implies that, for all $1 \leq * \leq k$,
$$
 \pi_*(\Omega^{\infty}_0 \Kz(\fB_G)) \cong H_*(\Omega^{\infty}_0 \Kz(\fB_G),\Zz).
$$
Now Theorem~\ref{thm:GroupCompl} implies that $H_*(\bmF(\cR \times G)) \cong \gekl{0}$ for all $0 < * < k$ and $H_k(\bmF(\cR \times G)) \cong H_k(G)$. And if, in addition, $\Gn$ has no isolated points and $G$ is minimal and has comparison, then Theorem~\ref{thm:HFG=HOKB} implies that $H_*(\bmF(G)) \cong \gekl{0}$ for all $0 < * < k$ and $H_k(\bmF(G)) \cong H_k(G)$, as desired.
\eproof
\setlength{\parindent}{0cm} \setlength{\parskip}{0.5cm}

We record the following immediate consequence.
\bcor
\label{cor:Acyclic}
Let $G$ be an ample groupoid with locally compact Hausdorff unit space $\Gn$. If $H_*(G) \cong \gekl{0}$ for all $* \geq 0$, then $\bmF(\cR \times G)$ is integrally acyclic, i.e., $H_*(\bmF(\cR \times G)) \cong \gekl{0}$ for all $* > 0$, and $\bmF(\cR \times G) = \bmD(\cR \times G)$. If, in addition, $\Gn$ has no isolated points and $G$ is minimal and has comparison, then $\bmF(G)$ is integrally acyclic, and $\bmF(G) = \bmD(G)$. 
\ecor

For commutator subgroups, we obtain the following consequences of Theorem~\ref{thm:HDG=HtiOKB}, Corollary~\ref{cor:HDG=HtiOKB} and \eqref{e:pitiOmega} as well as Corollary~\ref{cor:HVanish}.
\bcor
\label{cor:HDVanish}
Let $G$ be an ample groupoid with locally compact Hausdorff unit space $\Gn$. We always have $H_1(\bmD(\cR \times G)) \cong \gekl{0}$, and if, in addition, $\Gn$ has no isolated points and $G$ is minimal and has comparison, then $H_1(\bmD(G)) \cong \gekl{0}$. 

Now suppose that $k$ is an integer with $k \geq 2$, and that $H_*(G) \cong \gekl{0}$ for all $* < k$. Then $\bmD(\cR \times G) = \bmF(\cR \times G)$ and $H_*(\bmD(\cR \times G)) \cong \gekl{0}$ for all $0 < * < k$ and $H_k(\bmD(\cR \times G)) \cong H_k(G)$, and if, in addition, $\Gn$ has no isolated points and $G$ is minimal and has comparison, then $\bmD(G) = \bmF(G)$ and $H_*(\bmD(G)) \cong \gekl{0}$ for all $0 < * < k$ and $H_k(\bmD(G)) \cong H_k(G)$. 
\ecor

\bremark
In particular, this means that $\bmD(\cR \times G)$ is always perfect, and that $\bmD(G)$ is perfect if $\Gn$ has no isolated points and $G$ is minimal and has comparison. Perfection is also discussed in \cite{McDS,R-W}. Note that Matui proved in \cite{Mat15} that for second countable, locally compact, Hausdorff, minimal groupoids $G$ which are almost finite or purely infinite, $\bmD(G)$ is even simple.
\eremark

\subsection{Low degree exact sequences}

\btheo
\label{thm:H2H0H1H1}
Let $G$ be an ample groupoid, with locally compact Hausdorff unit space $\Gn$. There exists an exact sequence
$$
 \xymatrix{
 H_2(\bmD(\cR \times G)) \ar[r] & H_2(G) \ar[r] & H_0(G,\Zz/2) \ar[r]^{\hspace{-0.25cm} \zeta} & H_1(\bmF(\cR \times G)) \ar[r]^{\hspace{0.75cm} \eta} & H_1(G) \ar[r] & 0.
 }
$$
The map $\eta$ is determined by the property that the composition
$
 \xymatrix{
 H_1(\bmF(\cR \times G)) \ar[r]^{\hspace{0.75cm} \eta} & H_1(G) \ar[r] & H_1(\cR \times G)
 }
$
sends the class of an element $\sigma \in \bmF(\cR \times G)$ to the class of $1_{\sigma}$ in $H_1(\cR \times G)$, where the second map is induced by the canonical inclusion $G \into \cR \times G$.

The map $\zeta$ sends the class of $1_U$ in $H_0(G, \Zz/2)$ for a compact open subspace $U \subseteq \Gn$ to the class of the element $\tau \amalg \tau^{-1} \in \bmF((\cR \times G)^{U \oplus U}_{U \oplus U}) \subseteq \bmF(\cR \times G)$, where $\tau = \gekl{\bbs_{2,1}} \times U$.
\etheo

\bremark
Morita invariance of groupoid homology (see for instance \cite[\S~3]{Mat12} or \cite[\S~4]{Mil}) implies that the canonical inclusion $G \into \cR \times G$ induces isomorphisms $H_*(G) \cong H_*(\cR \times G)$ and $H_0(G,\Zz/2) \cong H_0(\cR \times G,\Zz/2)$. This is the reason why $\eta$ is determined by the composition
$$
 \xymatrix{
 H_1(\bmF(\cR \times G)) \ar[r]^{\hspace{0.75cm} \eta} & H_1(G) \ar[r] & H_1(\cR \times G).
 }
$$
\eremark

\bproof[Proof of Theorem~\ref{thm:H2H0H1H1}]
Let $\Sz$ be the sphere spectrum. The Atiyah-Hirzebruch spectral sequence (see for instance \cite[Part~III, \S~7]{Ada74}) has $E^2_{p,q} = \ti{H}_p(\Kz(\fB_G),\pi_q(\Sz))$ and converges to $\pi_{p+q}(\Kz(\fB_G))$. Since the Atiyah-Hirzebruch spectral sequence is a first quadrant spectral sequence (i.e., it satisfies $E^2_{p,q} \cong \gekl{0}$ whenever $p<0$ or $q<0$), we obtain a low degree exact sequence
\begin{equation}
\label{e:AHSS_original}
 \xymatrix{
 \pi_2(\Kz(\fB_G)) \ar[r] & \ti{H}_2(\Kz(\fB_G),\pi_0(\Sz)) \ar[r] & \ti{H}_0(\Kz(\fB_G),\pi_1(\Sz)) \ar[r] & \pi_1(\Kz(\fB_G)) \ar[r] & \ti{H}_1(\Kz(\fB_G)) \ar[r] & 0.
 }
\end{equation}
Now we plug in $\pi_0(\Sz) \cong \Zz$ and $\pi_1(\Sz) \cong \Zz/2$ (see \cite[\S~4.2]{Hat}), $\ti{H}_*(\Kz(\fB_G)) \cong H_*(G)$, $\ti{H}_0(\Kz(\fB_G),\Zz/2) \cong H_0(G,\Zz/2)$ from Theorem~\ref{thm:HKB=HG}, as well as
$$
 \pi_1(\Kz(\fB_G)) \cong \pi_1(\Omega^{\infty} \Kz(\fB_G)) \cong \pi_1(\Omega_0^{\infty} \Kz(\fB_G)) \cong H_1(\Omega_0^{\infty} \Kz(\fB_G)) \cong H_1(\bmF(\cR \times G)),
$$
where the third isomorphism is given by the Hurewicz homomorphism, the last isomorphism is from Theorem~\ref{thm:GroupCompl}, and
$$
 \pi_2(\Kz(\fB_G)) \cong \pi_2(\Omega^{\infty} \Kz(\fB_G)) \cong \pi_2(\Omega_0^{\infty} \Kz(\fB_G)) \cong \pi_2(\ti{\Omega}_0^{\infty} \Kz(\fB_G)) \cong H_2(\ti{\Omega}_0^{\infty} \Kz(\fB_G)) \cong H_2(\bmD(\cR \times G)),
$$
where the third isomorphism is from \eqref{e:pitiOmega}, the fourth isomorphism is given by the Hurewicz homomorphism and the last isomorphism is from Theorem~\ref{thm:HDG=HtiOKB}. We obtain the exact sequence
$$
 \xymatrix{
 H_2(\bmD(\cR \times G)) \ar[r] & H_2(G) \ar[r] & H_0(G,\Zz/2) \ar[r]^{\hspace{-0.25cm} \zeta} & H_1(\bmF(\cR \times G)) \ar[r]^{\hspace{0.75cm} \eta} & H_1(G) \ar[r] & 0,
 }
$$
as desired. It remains to determine the maps $\zeta$ and $\eta$.

Take a compact open subspace $u \subseteq \Nz \times \Gn = (\cR \times G)^{(0)}$ and $\sigma \in \bmF((\cR \times G)_u^u)$. Consider the commutative diagram
$$
 \xymatrix{
 \bmF((\cR \times G)_u^u) \ar[d] \ar[r] & H_1(\bmF((\cR \times G)_u^u)) \ar[d]\\
 \pi_1(\abs{\fB_G}) \ar[d] \ar[r] & H_1(\abs{\fB_G}) \ar[d]\\
 \pi_1(\Omega^{\infty}_0 \Kz(\fB_G)) \ar[d] \ar[r] & H_1(\Omega^{\infty}_0 \Kz(\fB_G)) \ar[d]\\
 \pi_1(\Kz(\fB_G)) \ar[r] & H_1(\Kz(\fB_G))
 }
$$
Here the first horizontal map is the canonical quotient map, all other horizontal maps are given by Hurewicz homomorphisms, and all vertical maps are the canonical ones. The diagram commutes by naturality of Hurewicz homomorphisms (see for instance \cite[Chapter~7, \S~4, Theorem~3]{Spa}. Moreover, note that the map $\pi_1(\Kz(\fB_G)) \to \ti{H}_1(\Kz(\fB_G))$ in \eqref{e:AHSS_original} is given by the stable Hurewicz homomorphism. Hence, in order to determine $\eta([\sigma])$, view $\sigma$ as an element of $C_{1,0} \fN_1 \fB_G$ (where we identify $\fB_G(S^1_0)$ with $\fB_G$) and determine the image of $[\sigma]$ under the isomorphism $H_1(\Kz(\fB_G)) \cong H_1(G)$ from Theorem~\ref{thm:HKB=HG}. 

We follow the proof of Theorem~\ref{thm:HKB=HG} and use the same notation. First, consider the exact sequence
$$
 \xymatrix{
 0 \ar[r] & \ker \partial_1 \ar[r] & C_{*,0} \fN_* \fB_G^{(1)} \ar[r]^{\partial_1} & C_{*,0} \fN_* \fB_G \ar[r] & 0.
 }
$$
In homology, we obtain the exact sequence
$$
 \xymatrix{
 0 \ar[r] & H_1(C_{*,0} \fN_* \fB_G) \ar[r] & H_0(\ker \partial_1) \ar[r] & H_0(C_{*,0} \fN_* \fB_G^{(1)}) \ar[r] & H_0(C_{*,0} \fN_* \fB_G) \ar[r] & 0.
 }
$$
We need the image of $[\sigma] \in H_1(C_{*,0} \fN_* \fB_G)$ under the map $H_1(C_{*,0} \fN_* \fB_G) \to H_0(\ker \partial_1)$. Observe that $(\sigma,u) \in \fN_1 \fB_G^{(1)}$ maps to $\sigma \in \fN_1 \fB_G$ under $\partial_1$. Under the map $C_{1,0} \fN_1 \fB_G^{(1)} \to C_{0,0} \fN_0 \fB_G^{(1)}$, $(\sigma,u)$ is mapped to $\sigma - u$. Hence $[\sigma] \in H_1(C_{*,0} \fN_* \fB_G)$ is mapped to $[\sigma] - [u]$ under the map $H_1(C_{*,0} \fN_* \fB_G) \to H_0(\ker \partial_1)$. 

Next, consider the exact sequence
$$
 \xymatrix{
 C_{0,0} \fN_0 \fB_G^{(3)} \ar[r]^{\partial_3} & C_{0,0} \fN_0 \fB_G^{(2)} \ar[r]^{\quad \partial_2} & \ker \partial_1 \ar[r] & 0.
 }
$$
In homology, $H_0(\partial_2)$ induces the identification $\coker H_0(\partial_3) \cong H_0(\ker \partial_1)$. Now $u \tensor[_{\rms}]{\times}{_{\rmr}} \sigma \in C_{0,0} \fN_0 \fB_G^{(2)}$ is mapped to $[\sigma] - [u]$ under $\partial_2$. Hence $H_0(\partial_2)$ sends $[u \tensor[_{\rms}]{\times}{_{\rmr}} \sigma]$ to $[\sigma] - [u]$. 

Now consider the commutative diagram
$$
 \xymatrix{
 \ti{H}_0(C_* \fN \fB_G^{(2)}) \ar[r]^{\cong} \ar[d] & \C(G^{(1)}) \ar[d]\\
 \ti{H}_0(C_* \fN \fB_{\cR \times G}^{(2)}) \ar[r]^{\cong} & \C((\cR \times G)^{(1)})
 }
$$
The horizontal maps are the isomorphisms $\ti{H}_0(C_* \fN \fB_G^{(2)}) \cong \C(G^{(1)})$ and $\ti{H}_0(C_* \fN \fB_{\cR \times G}^{(2)}) \cong \C((\cR \times G)^{(1)})$ from Theorem~\ref{thm:HBGnu}, and the vertical arrows are the maps induced by the canonical inclusion $G \into \cR \times G$. The left vertical arrow sends $[u \tensor[_{\rms}]{\times}{_{\rmr}} \sigma]$ to $[u \tensor[_{\rms}]{\times}{_{\rmr}} \sigma]$ (now viewing $u \tensor[_{\rms}]{\times}{_{\rmr}} \sigma$ as an element of $C_{0,0} \fN_0 \fB_{\cR \times G}^{(2)}$), which then is sent to $[1_{\sigma}]$ by the lower horizontal arrow.

So all in all, we conclude that the composition 
$
\xymatrix{
 H_1(\bmF(\cR \times G)) \ar[r]^{\hspace{0.75cm} \eta} & H_1(G) \ar[r] & H_1(\cR \times G)
 }
$
sends $[\sigma]$ to $[1_{\sigma}] \in H_1(G)$, as desired.

To determine $\zeta$, we use naturality of the Atiyah-Hirzebruch spectral sequence. Take a non-empty compact open subspace $U \subseteq \Gn$. Let $\fB_{\gekl{U}}$ be the category consisting of direct sums of $\emptyset$ and copies of $U$ and morphisms given by permutations of summands of $U$. In other words, $\fB_{\gekl{U}}$ is the small permutative category constructed in \S~\ref{s:smallpermcat} for the discrete groupoid $\cR \times \gekl{U}$. We have a canonical embedding $\fB_{\gekl{U}} \into \fB_G$.

Since $\fB_{\gekl{U}}$ is the free permutative category on $\gekl{\emptyset, U}$, the category with two objects $\emptyset$ and $U$ and only identity morphisms, we have $\ti{H}_0(\Kz \fB_{\gekl{U}},\Zz/2) \cong \Zz/2$, with generator $[U]$, $\ti{H}_2(\Kz \fB_{\gekl{U}}) \cong \gekl{0}$ and $\ti{H}_1(\Kz \fB_{\gekl{U}}) \cong \gekl{0}$. So the low degree exact sequence obtained from the Atiyah-Hirzebruch spectral sequence for $\Kz \fB_{\gekl{U}}$ degenerates to the isomorphism
$$
 \xymatrix{
 \ti{H}_0(\Kz \fB_{\gekl{U}},\Zz/2) \ar[r]^{\hspace{-0.75cm} \cong} & \pi(\Kz \fB_{\gekl{U}}) \cong H_1(\bmF(\cR))
 }
$$
sending the generator $[U]$ to the class of the non-trivial permutation $\pi: \: U \oplus U \cong U \oplus U$. Here we have applied Theorem~\ref{thm:GroupCompl} to the discrete groupoid $\cR \cong \cR \times \gekl{U}$ and used that $\bmF(\cR) \cong S_{\infty} = \bigcup_N S_N$, where $S_N$ is the symmetric group on a finite set of size $N$.

Hence by naturality of the Atiyah-Hirzebruch spectral sequence, $\zeta([1_U]) = [\tau \amalg \tau^{-1}]$ because under the maps induced by the canonical embedding $\fB_{\gekl{U}} \into \fB_G$, $[U]$ is mapped to $[1_U]$ and $[\pi]$ is mapped to $[\tau \amalg \tau^{-1}]$.
\eproof

\bcor
\label{cor:AHConj}
Let $G$ be an ample groupoid whose unit space $\Gn$ is locally compact Hausdorff and has no isolated points. Assume that $G$ is minimal and has comparison. Then there is an exact sequence
$$
 \xymatrix{
 H_2(\bmD(G)) \ar[r] & H_2(G) \ar[r] & H_0(G,\Zz/2) \ar[r]^{\zeta} & H_1(\bmF(G)) \ar[r]^{\hspace{0.25cm} \eta} & H_1(G) \ar[r] & 0.
 }
$$
The maps $\eta$ and $\zeta$ coincide with the ones in \cite[\S~2.3]{Mat16} and \cite[\S~7]{Nek19}. 

In particular, Matui's AH-conjecture is true for every ample groupoid $G$ which is minimal, has comparison and whose unit space is locally compact Hausdorff without isolated points.
\ecor
\setlength{\parindent}{0cm} \setlength{\parskip}{0cm}

In particular, this proves Matui's AH-conjecture for all purely infinite minimal ample groupoids, which was not known before. Corollary~\ref{cor:AHConj} also verifies the AH-conjecture for all minimal ample groupoids which are $\sigma$-compact, Hausdorff and almost finite, and whose unit spaces are compact without isolated points. Previously, this was only known under the additional assumption of principality \cite{Mat12}.
\bproof
We obtain the desired exact sequence by plugging in the isomorphism $H_*(\bmF(G)) \cong H_*(\bmF(\cR \times G))$ given by Theorem~\ref{thm:HFG=HOKB} into the exact sequence obtained in Theorem~\ref{thm:H2H0H1H1}.
\setlength{\parindent}{0cm} \setlength{\parskip}{0.5cm}

Let us explain why our maps $\eta$ and $\zeta$ coincide with the maps in Matui's AH-conjecture as in \cite[\S~2.3]{Mat16} and \cite[\S~7]{Nek19}. This is straightforward to see for $\eta$, so it remains to consider $\zeta$. The corresponding map in \cite[\S~2.3]{Mat16} and \cite[\S~7]{Nek19} is defined as follows: Given a compact open subspace $U \subseteq \Gn$ together with a compact open bisection $\sigma$ with $\rms(\sigma) = U$ and $\rmr(\sigma) \cap U = \emptyset$, the map in \cite[\S~2.3]{Mat16} and \cite[\S~7]{Nek19} sends $[1_U]$ to $\sigma \amalg \sigma^{-1}$. Let $\tau \amalg \tau^{-1}$ be as above. Let $\ti{\sigma}$ be the composition
$$
 \xymatrix{
 V & \ar[l]_{\hspace{-0.25cm} \sigma} (1,U) & \ar[l] (2,U),
 }
$$
where $V = \rmr(\sigma)$. Then $\ti{\sigma}^{-1}$ is given by the composition
$$
 \xymatrix{
 (2,U) & \ar[l] (1,U) & \ar[l]_{\hspace{0.5cm} \sigma^{-1}} V.
 }
$$
It follows that
$$
 (\ti{\sigma} \amalg \ti{\sigma}^{-1}) (\sigma \amalg \sigma^{-1}) (\ti{\sigma} \amalg \ti{\sigma}^{-1}) = \tau \amalg \tau^{-1}.
$$
Hence $[\tau \amalg \tau^{-1}] = [\sigma \amalg \sigma^{-1}]$ in $H_1(\bmF(\cR \times G))$. So our map $\zeta$ indeed coincides with the corresponding one in \cite[\S~2.3]{Mat16} and \cite[\S~7]{Nek19}.

Our claim about Matui's AH-conjecture follows immediately. (Note that the AH-conjecture is sometimes formulated with $H_0(G) \otimes (\Zz / 2)$ instead of $H_0(G,\Zz/2)$. However, these two groups are canonically isomorphic by the K{\"u}nneth formula.)
\eproof
\setlength{\parindent}{0cm} \setlength{\parskip}{0.5cm}

\bremark
\label{rem:AmplifiedAHConj}
An immediate consequence of Theorem~\ref{thm:H2H0H1H1} and Corollary~\ref{cor:AHConj} is that the stable version of Matui's AH-conjecture (with $H_1(\bmF(\cR \times G))$ in place of $H_1(\bmF(G))$) is always true for all ample groupoids. Equivalently, Matui's AH-conjecture is always true for groupoids of the form $\cR \times G$, where $G$ is an arbitrary ample groupoid.
\eremark

\bremark
Theorem~\ref{thm:H2H0H1H1} also implies that the strong AH-conjecture holds if $H_2(G) \cong \gekl{0}$. More precisely, by exactness, the strong AH-conjecture holds (i.e., the map $H_0(G,\Zz/2) \to H_1(\bmF(G))$ is injective) if and only if the map $H_2(G) \to H_0(G,\Zz/2)$ in Theorem~\ref{thm:H2H0H1H1} is the zero map.
\eremark

In addition to the alternating full group $\bmA(G)$, Nekrashevych also introduced the subgroup $\bmS(G)$ of $\bmF(G)$, which is an analogue of the symmetric group. By definition, $\bmA(G) \subseteq \bmS(G)$. Let $\bmK(G)$ be the kernel of the index map (see \cite[\S~2.3]{Mat16}, the index map coincides with $\eta$ in Theorem~\ref{thm:H2H0H1H1} and Corollary~\ref{cor:AHConj}). Clearly, $\bmD(G) \subseteq \bmK(G)$. As observed in \cite{Nek19}, it is easy to see that $\bmS(G) \subseteq \bmK(G)$ and $\bmA(G) \subseteq \bmD(G)$. Nekrashevych points out in \cite{Nek19} that \an{it would be interesting to understand when the equalities $\bmA(G) = \bmD(G)$ and $\bmS(G) = \bmK(G)$ hold}. Our work yields the following result about the relation between the subgroups $\bmA(G)$, $\bmD(G)$, $\bmS(G)$ and $\bmK(G)$ of $\bmF(G)$.
\bcor
Let $G$ be an ample groupoid whose unit space $\Gn$ is locally compact Hausdorff and has no isolated points. Assume that $G$ is minimal and has comparison. Then $\bmK(G)$ is generated by $\bmS(G)$ and $\bmD(G)$. Moreover, the following are equivalent:
\setlength{\parindent}{0cm} \setlength{\parskip}{0cm}

\begin{enumerate}
\item[(i)] $\bmD(G) \subseteq \bmS(G)$,
\item[(ii)] $\bmS(G) = \bmK(G)$,
\item[(iii)] $\bmA(G) = \bmD(G)$.
\end{enumerate}
\ecor
\setlength{\parindent}{0cm} \setlength{\parskip}{0cm}

\bproof
That $\bmK(G)$ is generated by $\bmS(G)$ and $\bmD(G)$ follows by combining the exact sequence in Corollary~\ref{cor:AHConj} with \cite[Theorem~7.2]{Nek19}.
\setlength{\parindent}{0cm} \setlength{\parskip}{0.5cm}

(i) $\Rarr$ (ii) is clear because $\bmK(G)$ is generated by $\bmS(G)$ and $\bmD(G)$. To see (ii) $\Rarr$ (iii), observe that (ii) produces an embedding $\bmD(G) / \bmA(G) \into \bmK(G) / \bmA(G) = \bmS(G) / \bmA(G)$. It follows that $\bmD(G) / \bmA(G)$ is abelian because \cite[Theorem~7.2]{Nek19} implies that $\bmS(G) / \bmA(G)$ is abelian. At the same time, we know that $H_1(\bmD(G)) \cong \gekl{0}$ by Corollary~\ref{cor:HDVanish}. Hence $\bmD(G) / \bmA(G) \cong \gekl{0}$, i.e., $\bmA(G) = \bmD(G)$.

(iii) $\Rarr$ (i) is clear because $\bmA(G) \subseteq \bmS(G)$.
\eproof
\setlength{\parindent}{0cm} \setlength{\parskip}{0.5cm}

\subsection{Examples}
\label{ss:Examples}

In the following, we present a few examples to illustrate our main results.

Let $\Zz \curvearrowright X$ be a Cantor minimal system and $G \defeq \Zz \ltimes X$ the corresponding transformation groupoid. As mentioned in \S~\ref{sss:GPDHomEx}, $H_1(G) \cong \Zz$ and $H_*(G) \cong \gekl{0}$ for all $*>1$. Thus, by Corollary~\ref{cor:RatHom},
$$
 H_*(\bmF(G),\Qz)
 \cong
 \bfa
 \Qz & \falls *=0 \text{ or } *=1,\\
 \gekl{0} & \sonst.
 \efa
$$
Moreover, Corollary~\ref{cor:DRatAcyclic} implies that $\bmD(G)$ is rationally acyclic, i.e., $H_*(\bmD(G),\Qz) \cong \gekl{0}$ for all $*>0$. In particular, this computes the rational homology of the examples of finitely generated infinite simple amenable groups found in \cite{JM}.

If $G$ is the transformation groupoid of a Cantor minimal $\Zz^d$-system with $d>1$, we have $H_d(G) \cong \Zz$, so that $H_d(\bmD(G),\Qz) \not\cong \gekl{0}$ by Corollary~\ref{cor:DRatHom}, and hence $\bmD(G)$ is not rationally acyclic.

For tiling groupoids as in \S~\ref{sss:Tiling}, explicit groupoid homology computations in for instance \cite{GK,FHK01,FHK02,GHK13} and Corollaries~\ref{cor:RatHom} and \ref{cor:DRatHom} lead to rational group homology computations for the corresponding topological full groups and their commutator subgroups. For instance, let $G$ be the groupoid attached to the classical Penrose tiling. Then 
$$
 H_*(G) \cong 
 \bfa 
 \Zz^8 & \falls *=0,\\
 \Zz^5 & \falls *=1,\\
 \Zz & \falls *=2,\\
 \gekl{0} & \sonst.
 \efa
$$
Hence, by Corollary~\ref{cor:P_FG}, the Poincar{\'e} series for $H_*(\bmF(G),\Qz)$ is given by 
$(1+t^1)^5 \, (1 - t^2)^{-1}$.
For $H_*(\bmD(G),\Qz)$, the Poincar{\'e} series is given by
$(1 - t^2)^{-1}$
by Corollary~\ref{cor:P_DG}, so that
$$
 H_*(\bmD(G),\Qz) \cong 
 \bfa
 \Qz & \falls * \text{ is even},\\
 \gekl{0} & \falls * \text{ is odd}.
 \efa
$$

Let $G_A$ be an SFT groupoid as in \S~\ref{sss:GraphGPD}, where the transition matrix $A$ is irreducible and not a permutation matrix. Let $d$ be the rank of $\ker(id - A^t)$. Then, using the groupoid homology results in \S~\ref{sss:GPDHomEx}, Corollary~\ref{cor:RatHom} implies that $H_*(\bmF(G_A),\Qz) \cong \Qz^{\binom{d}{*}}$, and Corollary~\ref{cor:DRatAcyclic} implies that $\bmD(G_A)$ is rationally acyclic.

Given a one vertex $k$-graph $\Lambda$ as in \S~\ref{sss:GPDHomEx}, the groupoid homology results in \S~\ref{sss:GPDHomEx}, Corollary~\ref{cor:RatAcyclic} and Corollary~\ref{cor:HVanish} imply that $\bmF(G_{\Lambda})$ is always rationally acyclic, and that $\bmF(G_{\Lambda})$ is even integrally acyclic if $\gcd(N_1, \dotsc, N_k) = 0$. In particular, the Brin-Higman-Thompson groups $n V_{k,r}$ are always rationally acyclic, and $n V_{k,r}$ are even integrally acyclic if $k=2$. Note that Brin's groups $nV$ from \cite{Bri04} coincide with $nV_{2,1}$ and hence are integrally acyclic. Moreover, Theorem~\ref{thm:GUUGVV} implies that $H_*(nV_{k,r},\sfC)$ does not depend on $r$, for all abelian groups $\sfC$ and $* \geq 0$. More precisely, for all $r \leq s$, the canonical embedding $nV_{k,r} \into nV_{k,s}$ induces isomorphisms $H_*(nV_{k,r},\sfC) \cong H_*(nV_{k,s},\sfC)$ for all abelian groups $\sfC$ and $* \geq 0$.

Consider a Katsura-Exel-Pardo groupoid $G_{A,B}$ as in \S~\ref{sss:SelfSimGPD}, where $A$ and $B$ are row-finite matrices with integer entries, and all entries of $A$ are non-negative. Suppose that for all $1 \leq i, j \leq N$, $B_{i,j} = 0$ if and only if $A_{i,j} = 0$. Further assume that $A$ is irreducible and not a permutation matrix. Let $d_A$ be the rank of $\ker(id - A^t)$ and $d_B$ the rank of $\ker(id - B^t)$. Then, using the groupoid homology results in \S~\ref{sss:GPDHomEx}, Corollary~\ref{cor:P_FG} implies that the Poincar{\'e} series of $H_*(\bmF(G_{A,B}),\Qz)$ is given by
$(1 + t)^{d_A + d_B} (1 - t^2)^{-d_B}$,
and Corollary~\ref{cor:P_DG} implies that the Poincar{\'e} series of $H_*(\bmD(G_{A,B}),\Qz)$ is given by
$(1 - t^2)^{-d_B}$.

Let us now discuss groupoids arising from piecewise affine transformations on the unit interval as in \S~\ref{sss:AffineGPD}. First let $\lambda$ be an algebraic integer with $1 \neq \lambda \in (0,\infty)$ whose minimal polynomial is given by $f(T) = T^d + a_{d-1} T^{d-1} + \dotso + a_1 T + a_0$. Let $G$ be the groupoid from \S~\ref{sss:AffineGPD} for parameter $\lambda$ (see also \cite{Li15}). The concrete computations of groupoid homology in \cite[\S~5.2]{Li15} and Corollaries~\ref{cor:RatHom}, \ref{cor:RatAcyclic}, \ref{cor:DRatHom} and \ref{cor:DRatAcyclic} imply the following:
\setlength{\parindent}{0cm} \setlength{\parskip}{0cm}

\begin{itemize}
\item If $d=2$ and $a_0 \neq 1$, then $\bmF(G)$ and $\bmD(G)$ are rationally acyclic.
\item If $d=2$ and $a_0 = 1$, then $H_*(\bmF(G),\Qz) \cong \Qz$ for all $* \geq 0$, and $H_*(\bmD(G),\Qz) \cong \Qz$ for all even $* \geq 0$ and $H_*(\bmD(G),\Qz) \cong \gekl{0}$ for all odd $* \geq 1$.
\item If $d=3$ and $a_0 \neq -1$, then $\bmF(G)$ and $\bmD(G)$ are rationally acyclic.
\item If $d=3$ and $a_0 = -1$, then $H_*(\bmF(G),\Qz) \cong \Qz$ for all $* \geq 0$ with $* \neq 1$ and $H_1(\bmF(G),\Qz) \cong \gekl{0}$, and $H_*(\bmD(G),\Qz) \cong \Qz$ for all $* \geq 0$ with $* \neq 1$ and $H_1(\bmD(G)) \cong \gekl{0}$.
\end{itemize}
If $\lambda$ is transcendental and $G$ is the groupoid from \S~\ref{sss:AffineGPD} for parameter $\lambda$, then the computations mentioned in \cite[\S~6]{Li15} and Corollaries~\ref{cor:RatHom}, \ref{cor:DRatHom} imply that $H_*(\bmF(G),\Qz) \cong \bigoplus_{i=0}^{\infty} \Qz$ for all $* \geq 0$ and $H_*(\bmD(G),\Qz) \cong \bigoplus_{i=0}^{\infty} \Qz$ for all $* \geq 0$ with $* \neq 1$, and $H_1(\bmD(G)) \cong \gekl{0}$.
\setlength{\parindent}{0cm} \setlength{\parskip}{0.5cm}

Finally, we discuss examples where we can apply our vanishing and acyclicity results (Corollaries~\ref{cor:HVanish}, \ref{cor:Acyclic}, \ref{cor:HDVanish}). Suppose $\sfA$ is a countably generated abelian group. Using Katsura-Exel-Pardo groupoids and combining results in \cite[\S~4]{Kat08a}, \cite[\S~3]{Kat08b} and \cite{Ort}, we can find purely infinite minimal groupoids $G(0,\sfA)$ and $G(1,\Zz)$ such that
$$
 H_*(G(0,\sfA)) \cong
 \bfa
 \sfA & \falls *=0,\\
 \gekl{0} & \sonst,
 \efa
 \qquad
 \text{and }
 H_*(G(1,\Zz)) \cong
 \bfa
 \Zz & \falls *=1,\\
 \gekl{0} & \sonst.
 \efa
$$
Given $k \in \Zz$ with $k > 0$, the K{\"u}nneth formula (see \cite[Theorem~2.4]{Mat16}) implies that
$$
 H_*(G(1,\Zz)^k \times G(0,\sfA)) \cong
 \bfa
 \gekl{0} & \falls * < k,\\
 \sfA & \falls * = k.
 \efa
$$
Hence Corollary~\ref{cor:HVanish} implies that
$$
 H_*(\bmF(G(1,\Zz)^k \times G(0,\sfA))) \cong
 \bfa
 \Zz & \falls * = 0,\\
 \gekl{0} & \falls * < k,\\
 \sfA & \falls * = k,
 \efa
$$
and, if $k \geq 2$, Corollary~\ref{cor:HDVanish} implies that $\bmD(G(1,\Zz)^k \times G(0,\sfA)) = \bmF(G(1,\Zz)^k \times G(0,\sfA))$.

Let us now turn to acyclicity results. Let $G_2$ be the Deaconu-Renault groupoid for the one-sided full shift on two symbols (this is a special case of an SFT groupoid as in \S~\ref{sss:GraphGPD}, where $A$ is the $1 \times 1$ matrix with entry $2$). Let $G$ be an arbitrary minimal ample groupoid. Then $G_2 \times G$ is purely infinite minimal. Moreover, the K{\"u}nneth formula (see \cite[Theorem~2.4]{Mat16}) implies that $H_*(G_2 \times G) \cong \gekl{0}$ for all $* \geq 0$. Hence Corollary~\ref{cor:Acyclic} implies that $\bmF(G_2 \times G)$ is integrally acyclic and $\bmF(G_2 \times G) = \bmD(G_2 \times G)$.

\bremark
\label{rem:ManyAcyclic}
In combination with the groupoids constructed in \cite[\S~9.2]{CL21}, we obtain continuum many pairwise non-isomorphic infinite simple groups which are all integrally acyclic. Indeed, consider the groupoids of the form $G_2 \times \cG_{\Gamma}$ from \cite[\S~9.2]{CL21}, where $\Gamma$ is an abelian, torsion-free, finite rank group which is not free abelian. By construction, these groupoids are ample, locally compact, Hausdorff, purely infinite, topologically free, with unit space homeomorphic to the Cantor space. Moreover, as observed above, $\bmF(G_2 \times \cG_{\Gamma})$ is integrally acyclic, and $\bmF(G_2 \times \cG_{\Gamma}) = \bmD(G_2 \times \cG_{\Gamma})$. Hence \cite[Theorem~4.16]{Mat15} implies that $\bmF(G_2 \times \cG_{\Gamma})$ is simple. Now the rigidity results in \cite[Theorem~3.10]{Mat15} and \cite[Theorem~3.11]{Nek19} together with the argument for \cite[Theorem~9.3]{CL21} imply that, for two abelian, torsion-free, finite rank groups $\Gamma$ and $\Lambda$ which are not free abelian, $\bmF(G_2 \times \cG_{\Gamma}) \cong \bmF(G_2 \times \cG_{\Lambda})$ if and only if $G_2 \times \cG_{\Gamma} \cong G_2 \times \cG_{\Lambda}$ if and only if $\Gamma \cong \Lambda$. Thus we obtain continuum many pairwise non-isomorphic infinite simple, integrally acyclic groups because there are continuum many pairwise non-isomorphic abelian, torsion-free, finite rank groups which are not free abelian.
\eremark

\end{document}